\magnification=\magstep 1
\input amssym.def

\baselineskip=18pt \vsize = 9truein \hsize=5.5truein
\hoffset=.5in

\font\title=cmbx10 scaled \magstep 2
\font\sc=cmcsc10
 \font\smalli=cmti8
 1
\def\m@th{\mathsurround=0pt}
\def\recase#1{\phantom\{ \vcenter{\normalbaselines\m@th
     \ialign{##\hfil & \quad ## \crcr#1\crcr}}\Biggr\} }
\def\recasebig#1{\left. \vcenter{\normalbaselines\m@th
     \ialign{##\hfil & \quad ## \crcr#1\crcr}}\right\} }

\def \R{\bf R}
\def \N{\bf N}
\def \Z{\bf Z} 
\def \cA{{\cal A}}
\def \cB{{\cal B}}
\def \cC{{\cal C}}
\def \cD{{\cal D}}
\def \cE{{\cal E}}
\def \cG{{\cal G}}
\def \cH{{\cal H}}
\def \cI{{\cal I}}

\def \cL{{\cal L}}
\def \cJ{{\cal J}}
\def \cM{{\cal M}}
\def \cP{{\cal P}}
\def \cF{{\cal F}}
\def \cQ{{\cal Q}}
\def \cR{{\cal R}}
\def \cS{{\cal S}}
\def \cT{{\cal T}}

\def \basic{{\rm basic}}
\def \card{{\rm card}}
\def \cen{{\rm cen}}

\def \double{{\rm double}}
\def \Fri{{\rm fri}}
\def \fri{{\rm fri}}
\def \good{{\rm good}}

\def \key{{\rm key}}

\def \ord{{\rm ord}}
\def \pos{{\rm pos}}

\def \spaced{{\rm spaced}}
\def \sxtp{{\rm sxtp}}
\def \union{{\rm union}}
\def \var{\hbox{\rm Var}}
\hyphenation{argu-ment}

\def\ref{\bigskip\hangindent=25pt\hangafter=1\noindent}
\def\refs{\medskip\hangindent=25pt\hangafter=1\noindent}
\overfullrule= 0 pt 
\parskip =0pt

\centerline{\bf A STRICTLY STATIONARY, ``CAUSAL,'' 
5--TUPLEWISE}
\centerline{\bf INDEPENDENT COUNTEREXAMPLE TO THE} 
\centerline{\bf CENTRAL LIMIT THEOREM}
\bigskip

\centerline{\sc Richard C.\ Bradley} \centerline{\smalli Department
of Mathematics, Indiana University, Bloomington, Indiana 47405,
U.S.A.} \centerline{\smalli E-mail: bradleyr@indiana.edu}
\bigskip

\par{\bf Abstract.} A strictly stationary sequence of random variables is
constructed with the following properties: (i)~the random variables
take the values $-1$ and $+1$ with probability 1/2 each, (ii)~every
five of the random variables are independent of each other,
(iii)~the sequence is ``causal'' in a certain sense,
(iv)~the sequence has a trivial double tail $\sigma$-field,
 and
(v)~regardless of the normalization used, the partial sums 
do not
converge to a (nondegenerate) normal law. 
The example has some features in common with a recent
construction (for an arbitrary fixed positive
integer $N$), by Alexander Pruss and the author, of a
strictly stationary $N$-tuplewise independent
counterexample to the central limit theorem.
\medskip

\par {\it AMS} 2010 {\it Mathematics Subject Classification.} 60G10, 60F05.
\medskip

\hfill\break

\vskip .2truein \vfill\eject

\centerline {\bf 1. Introduction} 
\bigskip
For a given integer $N \geq 2$ and a given
sequence $X := (X_k, k \in {\bf Z})$ of random variables defined on
a probability space $(\Omega, {\cal F}, P)$, the random variables
$X_k, k \in {\bf Z}$ are said to be ``$N$-tuplewise independent'' if
for every choice of $N$ distinct integers $k(1), k(2), \dots, k(N)$,
the random variables $X_{k(1)}, X_{k(2)},\dots, X_{k(N)}$, are
independent. For $N=2$ (resp.\ $N=3$), the word ``$N$-tuplewise'' is
also expressed as ``pairwise'' (resp.\ ``triplewise'').

  Etemadi [14] proved a strong law of large numbers
for sequences of pairwise independent, identically distributed
random variables with finite absolute first moment. 
Janson [19] showed with several classes of counterexamples that for strictly stationary sequences of pairwise independent, nondegenerate, square-integrable random variables, the Central Limit Theorem (henceforth abbreviated CLT) need not hold. 
Subsequently, the author [3, Theorem 1] constructed another such counterexample, a 3-state one that has the additional property of satisfying the absolute regularity (weak Bernoulli) condition. (The definition of that condition will be given later in this Introduction.)\ \
Yet another counterexample was constructed by Cuesta and Matran [10, Section 2.3], a construction based on elementary 
number-theoretic properties of addition on 
$\{0, 1, \dots, p-1\}$ mod $p$, where $p$ is a prime number.
    
For an arbitrary fixed integer $N \geq 3$, Pruss [28] constructed a (not strictly stationary) sequence of bounded, nondegenerate, $N$-tuplewise independent, identically distributed random variables for which the CLT fails to hold. In that paper, Pruss left open the question whether, for any integer $N \geq 3$, a strictly stationary counterexample exists.     
For $N = 3$, the author [4, Theorem 1] answered that question affirmatively by showing that the counterexample in [3,
Theorem 1] alluded to above is in fact triplewise independent.
Recently, the author and Pruss [6] have answered that question
affirmatively for arbitrary $N \geq 3$, with a (strictly stationary) counterexample adapted from the (nonstationary) one in Pruss [28].

  In a similar spirit, for an arbitrary integer
$N \geq 2$, Flaminio [16] constructed a 
(nondegenerate) strictly stationary,
finite-state, $N$-tuplewise independent random sequence $X := (X_k, k \in {\bf Z})$ which also has zero entropy and is mixing (in the ergodic-theoretic sense). 
That paper explicitly left open the
question of whether those examples satisfy the CLT.

  In this paper here, a (nondegenerate) strictly stationary,
5-tuplewise independent counterexample to the CLT will be
constructed which is finite-state and has the extra property of being ``causal'' and therefore ``Bernoulli.'' 
The interest in the property of ``Bernoulli,'' for (finite-state) $N$-tuplewise independent counterexamples to the CLT, was suggested to the author by Jon Aaronson and Benjamin Weiss. 
There does not seem to be a visible way of constructing such a (strictly stationary) finite-state Bernoulli counterexample which is $N$-tuplewise independent for any given $N \geq 6$.  The techniques in the example given here do not appear to adapt effectively to $N \geq 6$; and there is no visible way of adapting the recent example of the author and Pruss [6] alluded to above, into one that is (finite-state and)
Bernoulli.  
The example given here will also have the further
property of possessing a trivial double tail $\sigma$-field.  (The terms ``causal,'' ``Bernoulli,'' and ``double tail $\sigma$-field'' will be defined below.)
\medskip

  {\it The main result and comments on it.}\  Before the
result is stated, some notations will be needed.

  Let ${\bf N}$ denote the set of all positive integers.
Let ${\cal R}$ denote the Borel $\sigma$-field on the real number line ${\bf R}$.  
The notation $\Rightarrow$ will mean convergence in
distribution.

  Suppose $X := (X_k, k \in {\bf Z})$ is a sequence of
random variables on a probability space $(\Omega, {\cal F}, P)$.
For each positive integer $n$, define the partial sum
$$ S_n := S(X,n) := X_1+ X_2 + \dots + X_n.  \leqno (1.1)$$
Also, for $-\infty \leq J \leq L \leq \infty$, let ${\cal F}_J^L$ denote the $\sigma$-field $\subset {\cal F}$ generated by the random variables $X_k,\ J \leq k \leq L\ (k \in {\bf Z})$. 
The ``double tail $\sigma$-field'' of the sequence $X$ is
$$ {\cal T}_{\rm double}(X) :=
\bigcap_{n \in {\bf N}} ({\cal F}_{-\infty}^{-n} \vee {\cal
F}_n^\infty).  \leqno (1.2) $$ A $\sigma$-field ${\cal A} \subset {\cal F}$ is said to be ``trivial'' if $P(A) = 0$ or 1 for every $A \in {\cal A}$.

  Here is our main result:
\bigskip

  {\sc Theorem 1.1.} {\it There exists a strictly stationary
sequence $X := (X_k$, $ k \in {\bf Z})$ of random variables  (on some probability space $(\Omega, {\cal F}, P)$) with the following six properties:

  (A)  The random variables $X_k$ take just the values
$-1$ and $1$, with $P(X_0 = -1) = P(X_0 = 1) = 1/2$ (and hence $EX_0 = 0$ and $EX_0^2 = 1$).

  (B)  For every five distinct integers
$k(1),\ k(2),\ k(3),\ k(4),\ {\rm and}\ k(5)$, the five random
variables $X_{k(1)},\ X_{k(2)},\ X_{k(3)},\ X_{k(4)},\ {\rm and}\ X_{k(5)}$ are independent.

  (C)  There exist a sequence $(\eta_k, k \in {\bf Z})$ of independent, identically distributed real-valued random variables (on $(\Omega, {\cal F}, P)$), and a Borel function $f :
{\bf R} \times {\bf R} \times {\bf R} \times \dots \to \{-1,1\}$,
such that for every $k \in {\bf Z}$,
$$ X_k = f(\eta_k, \eta_{k-1}, \eta_{k-2}, \dots)\ \ 
{\rm a.s.} \leqno (1.3) $$

  (D)  The double tail $\sigma$-field of $X$ is trivial
(that is, $P(A) = 0$ or $1$ for every $A \in {\cal T}_{\rm
double}(X)$).

  (E) One has that $\limsup_{n \to \infty}
E(S_n/\sqrt n \thinspace )^6 < 15$.

  (F) For every infinite set $Q \subset {\bf N}$,
there exist an infinite set $T \subset Q$ and a nondegenerate,
non-normal probability measure $\mu$ on $({\bf R}, {\cal R})$ such that $S_n/\sqrt n \Rightarrow \mu$ as $n \to \infty,\ n \in T$.}
\bigskip

  Here are some comments on the various properties in
this theorem --- starting with properties (A), (B), (E), and (F), the ones most closely tied to the central limit question.

  First, by property (A), in our look at the central limit
question in connection with this example, the natural normalization of the partial sums is $S_n/\sqrt n$.

  Property (B) in Theorem 1.1 is of course 5-tuplewise
independence. For $N \geq 6$, the question of possible existence of a similar, strictly stationary, $N$-tuplewise independent counterexample --- including properties (C) and (D) --- remains open. 
The techniques in the construction for Theorem 1.1 do not appear to extend effectively to $N \geq 6$.
As a comparison, for a given arbitrary fixed positive
integer $N$, the paper of the author and Pruss [6] gives
a construction of a strictly stationary sequence in which
(i)~the random variables are uniformly distributed on the
interval $[-\sqrt 3, \sqrt 3\, ]$ (and hence have mean 0 and variance 1),
(ii)~the sequence satisfies $N$-tuplewise independence,
and (iii)~the sequence satisfies property (F) 
(as well as a variant of property (E)) in Theorem 1.1;
that sequence satisfies ergodicity (as was shown in [6]),
but does not satisfy property (C) (suitably
reformulated) or property (D) in Theorem 1.1.
     
  Property (E) may seem rather pointless at first.
However, property (F) is an elementary consequence of properties (A), (B), and (E) together with the fact that a $N(0,1)$ random variable $Z$ satisfies $EZ^6 = 15$. (The argument will be given in detail in section 9 below.
An analogous argument, involving moments of a high even
order, was used by the author and Pruss [6].)

    In property (F), the probability measure $\mu$ may
depend on the set $Q$.

(As a comparison, in a couple of the pairwise independent
counterexamples alluded to above --- one of 
those in Janson [19] and the one in Cuesta and Matran [10, Section 2.3] --- the partial sums, appropriately normalized,
converge in distribution to a nondegenerate, non-normal
law as $n \to \infty$ along the entire sequence of positive integers.)

  The formulation of property (F) may seem somewhat
awkward. However, properties (E) and (F) indirectly give the
following information:

  (i) The family of distributions of the random variables
$(S_n/\sqrt n,\ n \in {\bf N})$ is tight.

  (ii) There does not exist an infinite set
$Q \subset {\bf N}$ such that $S_n/\sqrt n$ converges to 0 (or to
any other constant) in probability as $n \to \infty,\ n \in Q$.

  (iii) There does not exist an infinite set
$Q \subset {\bf N}$ such that $S_n/\sqrt n$ converges in
distribution to a (nondegenerate) normal law as $n \to \infty,\ n
\in Q$.

  (iv) By (F) and the Theorem of Types (see e.g.\
[1, Theorem 14.2]), there do not exist an infinite set $Q \subset
{\bf N}$ and real numbers $a_n,\ b_n,\ n \in Q$, with $b_n \to
\infty$ as $n \to \infty,\ n \in Q$, such that $(S_n-a_n)/b_n$
converges in distribution to a nondegenerate normal law as $n \to
\infty,\ n \in Q$.

  One can describe property (C) by saying that the
random sequence $X$ is ``causal.''  
Such uses of that term are well known in the literature.  
See e.g.\ Brockwell and Davis [7] for its
use in the context of linear models in time series analysis.

  Properties (C) and (D) are motivated partly by the
general question of what properties from ergodic theory can help insure that a CLT holds. 
In order to elaborate on that, we will need
to give some more definitions and background information.
\medskip

  {\it The ``Bernoulli'' property.}\ \ Suppose $A$ is
nonempty finite set, and $X := (X_k, k \in {\bf Z})$ is a
nondegenerate strictly stationary sequence of random variables
taking their values in $A$. 
This sequence $X$ is said to be ``Bernoulli'' if, without changing its distribution (on $A^{\bf Z}$), it can be represented as a stationary coding of a finite-state
i.i.d.\ sequence --- that is, if $X$ can be represented in the form
$$   X_k = h\Bigl((\dots, Y_{k-1}, Y_k),
(Y_{k+1}, Y_{k+2}, \dots)\Bigl)  \leqno (1.4)  $$ for $k \in {\bf
Z}$, where $Y := (Y_k, k \in {\bf Z})$ is a sequence of independent,
identically distributed random variables taking their values in a finite set $B$, and $h: B^{\bf Z} \to A$ is a Borel function. 
(Here, a given ``two sided'' sequence $b := (b_k, k \in {\bf Z})$ of elements of $B$ is written as $((\dots, b_{-1}, b_0), (b_1, b_2, \dots))$ as a convenient way to avoid ambiguity.) This is one of
numerous equivalent ways of formulating the class of (strictly
stationary, finite-state) random sequences that are ``Bernoulli.''
For a list of some others, see e.g.\ Shields [31, p.\ 235, lines 12-15].

  By a special case of a classic theorem of Ornstein 
(see [24, p.\ 350, line 11], with reference to [23]),
properties ((A) and) (C) in Theorem 1.1 imply the Bernoulli property.
That is:
\medskip

  {\sc Remark 1.2.} Automatically, the random sequence
$X$ in Theorem 1.1 is Bernoulli.
\medskip

  In ergodic theory, it is well known that a nondegenerate strictly stationary, finite-state sequence with zero
entropy is not Bernoulli; see e.g.\ [26, section 6.4]. 
It follows that the examples of Flaminio [16] alluded to above are not Bernoulli. 
Consequently, even if those examples of Flaminio turn out
to be counterexamples to the CLT (apparently still an open
question), Theorem 1.1 (with Remark 1.2) still gives new information (beyond the examples given by Flaminio [16] and by the author and Pruss [6]) in that it provides a 5-tuplewise independent counterexample which has the additional property of being Bernoulli.
\medskip

  {\it Two strong mixing conditions.}\ \ Next, we would like
to use Theorem 1.1 to obtain some perspective on a classic CLT
(stated in Theorem 1.4 below) involving Rosenblatt's [29] ``strong
mixing'' condition.   This will require some more definitions and
background information.

  Again suppose $X := (X_k, k \in {\bf Z})$ is a
strictly stationary sequence of (real-valued, not necessarily
finite-state) random variables on a probability space $(\Omega,
{\cal F}, P)$.

For any two $\sigma$-fields ${\cal A}$ and ${\cal B} \subset {\cal F}$, define the measures of dependence
$$ \eqalign{
\alpha({\cal A}, {\cal B}) &:= \sup_{A \in {\cal A},\ B \in {\cal
B}} |P(A \cap B) - P(A)P(B)| \indent {\rm and} \cr \beta({\cal A},
{\cal B}) &:=\  \sup\ (1/2) \sum_{i=1}^I \sum_{j=1}^J |P(A_i \cap
B_j) - P(A_i)P(B_j)| \cr }$$ where the latter supremum is taken over all pairs of finite partitions $\{A_1, \dots, 
\allowbreak A_I\}$ and $\{B_1,
\dots, B_J\}$ of $\Omega$ such that $A_i \in {\cal A}$ for each $i$ and $B_j \in {\cal B}$ for each $j$. 
(The factor of $1/2$ in the definition of $\beta({\cal A}, {\cal B})$ is of no significance, but has become standard in the literature.)\ \  
For each positive integer $n$, define the dependence coefficients
$$ \eqalign{
\alpha(n) := \alpha(X,n) &:= \alpha({\cal F}_{-\infty}^0, {\cal
F}_n^\infty) \indent {\rm and} \cr \beta(n) := \beta(X,n) &:=
\beta({\cal F}_{-\infty}^0, {\cal F}_n^\infty).  \cr }$$ By strict
stationarity, $\alpha({\cal F}_{-\infty}^j, {\cal F}_{j+n}^\infty) =
\alpha(n)$ and $\beta({\cal F}_{-\infty}^j, {\cal F}_{j+n}^\infty) =
\beta(n)$ for every integer $j$. The (strictly stationary) sequence
$X$ is said to satisfy the Rosenblatt [29] ``strong mixing''
condition, or ``$\alpha$-mixing,'' if $\alpha(n) \to 0$ as $n \to \infty$; and it is said to satisfy the ``absolute regularity'' [33] condition if $\beta(n) \to 0$ as $n \to \infty$.  
For the first of those conditions, the term ``$\alpha$-mixing'' will be used here in order to avoid ambiguity from conflicting uses of the phrase ``strong mixing'' in the literature.

Obviously absolute regularity implies $\alpha$-mixing. For strictly stationary, finite-state sequences, one also has the following:

(i) First, $\alpha$-mixing does not imply absolute regularity (see e.g.\ [5, Vol.\ 1, Theorem 9.10(II)]).

(ii) The ``weak Bernoulli'' condition, defined and studied by
Friedman and Ornstein [17], is equivalent to absolute regularity,
and (as was shown in that paper) it implies the Bernoulli property.

(iii)  It is unknown (an open problem posed by Donald Ornstein in the 1970s) whether $\alpha$-mixing implies the Bernoulli property.

(iv) The Bernoulli property does not imply $\alpha$-mixing (and hence also does not imply absolute regularity); that was shown by Smorodinsky [32].

  The strictly stationary, 3-state, triplewise independent,
absolutely regular counterexample (to the CLT) developed in [3,4], alluded to above, is in two respects ``optimal'' under absolute regularity: 
First, with its random variables being bounded (in fact,
finite-state), its ``mixing rate'' $\beta(n) = O(1/n)$ (as $n \to\infty$) is essentially as rapid as possible (it cannot satisfy $\beta(n) = o(1/n)$ or even $\alpha(n) = o(1/n)$, by a CLT of Merlev\`ede and Peligrad [21]). 
Second, as was pointed out in [4, section 1] with a brief explanation, if a given strictly stationary
sequence of nondegenerate, square-integrable random variables
satisfies both $\alpha$-mixing (or absolute regularity) and
4-tuplewise independence, then it satisfies the CLT.  
(This fact is an elementary corollary, via a truncation argument, of the CLT under $\alpha$-mixing given in Theorem 1.4 below.)\ \ As a consequence:
\medskip

  {\sc Remark 1.3.} The random sequence $X$ in Theorem 1.1
cannot satisfy $\alpha$-mixing (or absolute regularity).
\medskip

  In the latter part of the book [5, Vol.\ 3], there is a
detailed presentation of a large collection of strictly stationary, absolutely regular (but in most cases not 
pairwise independent or finite-state) counterexamples to the CLT, including examples of Davydov [11], and also including a (slightly embellished) presentation of the well known example of Herrndorf [18] in which the random variables are uncorrelated.
\medskip

  {\it The double tail $\sigma$-field.}\ \ Let us digress
to take a quick look at some of the ways property (D) in Theorem 1.1 fits in with the various dependence conditions above, for strictly stationary sequences.

(i) First, it is well known and elementary (see e.g.\ [5,
Vol.\ 1, Proposition 5.17]) that absolute regularity implies a trivial double tail $\sigma$-field.

(ii) Next, $\alpha$-mixing does not imply a trivial double tail $\sigma$-field.  See for example the counterexample in [2] or [5, Vol.\ 2, Theorem 24.14] with real state space, or the finite-state counterexample constructed by Burton, Denker and Smorodinsky [8].
Each of those examples is ``bilaterally deterministic'' 
(that is, the entire random sequence is, modulo null-sets, measurable with respect to its double tail $\sigma$-field).

(iii) In the finite-state case, the Bernoulli property does not imply a trivial double tail $\sigma$-field.  Ornstein and Weiss [25] showed instead that within the class of strictly stationary, finite-state random sequences that are Bernoulli, the ones that are also bilaterally deterministic are in a certain sense ``ubiquitous.''

(iv) In the finite-state case, a trivial double tail $\sigma$-field does not imply $\alpha$-mixing (see e.g.\ [5, Vol.\ 1, Theorem 9.11(II)]), and hence also does not imply absolute regularity.
\medskip

  {\it A classic CLT under $\alpha$-mixing.}\ \  Now let us take a quick look at the
following classic theorem:
\bigskip

 {\sc Theorem 1.4.} {\it Suppose $X := (X_k, k \in {\bf Z})$
is a strictly stationary sequence of (real-valued) random variables
such that $EX_0 = 0$, $EX_0^2 < \infty$, $\sigma_n^2 := ES_n^2 \to
\infty$ as $n \to \infty$, and $\alpha(n) \to 0$ as $n \to \infty$.
Then the following two conditions (I), (II) are equivalent:

  (I) The family of random variables
$(S_n^2/\sigma_n^2,\ n \in {\bf N})$ is uniformly integrable.

  (II) $S_n/\sigma_n \Rightarrow N(0,1)$ as $n \to \infty$.}
\bigskip

  Even if the assumption of $\alpha$-mixing were omitted
altogether, (II) implies (I), by a well known, elementary argument
(see e.g.\ [13] or [22] or [5, Vol.\ 1, top half of p.\ 38]). The interest here in Theorem 1.4 is the fact that under all of the given assumptions, (I) implies (II).  That fact was shown by Cogburn [9, Theorem 13] (with $\alpha$-mixing replaced by a similar but technically weaker condition), a reference that did not seem to be well known for a long time; and its proof was also given by Denker [13] and by Mori and Yoshihara [22]. (A proof of Theorem 1.4 is spelled out in generous detail in [5, Vol.\ 1, Theorem 1.19].)

  By properties (A), (B), and (E) in Theorem 1.1, the
sequence $X$ in Theorem 1.1 satisfies $EX_0 = 0$, $EX_0^2 < \infty$,
$\sigma_n^2 := ES_n^2 = n$ for $n \in {\bf N}$, as well as the
uniform integrability of $(S_n^2/n,\ n \in {\bf N})$. Hence by
properties (A), (C), (D), and (F) in Theorem 1.1, together with
Remark 1.2, one has the following:
\medskip

  {\sc Remark 1.5.} In Theorem 1.4 (for the assertion that
(I) implies (II)), even if the random variables $X_k$ are
finite-state and 5-tuplewise independent, the assumption of
$\alpha$-mixing cannot be replaced by either (or both) of the
assumptions that \hfil\break (i) $X$ is Bernoulli (or even that Property (C) in Theorem 1.1 holds),

\noindent (ii) $X$ has a trivial double tail $\sigma$-field.
\medskip

  Suppose the hypothesis of Theorem 1.4 holds, along with
the extra assumptions that (i) $\sigma_n^2 = n \cdot h(n)$ where $h:
(0,\infty) \to (0,\infty)$ is slowly varying at $\infty$, and  (ii)
$ \limsup_{n \to \infty} \|S_n\|_2/ [(\pi/2)^{1/2}E|S_n|\thinspace
]\ \leq 1$. 
Then the CLT (specifically, conclusion (II) of Theorem 1.4) holds.  
That is a result of Dehling, Denker, and Philipp [12,
Theorem 4]. 
(For an exposition of their argument in generous detail,
see [5, Vol.\ 2, Theorem 17.11].)\ \ 
It seems to be an open question whether
that still holds if the assumption of $\alpha$-mixing is replaced by (say) properties (B), (C), and (D) in Theorem 1.1.  Because of the absolute-value signs, the quantities $E|S_n|$ seem to be hard to estimate effectively for the construction given below for Theorem 1.1.

  The rest of this paper is devoted to the proof of
Theorem 1.1.  Here is how that proof will be organized, in Sections 2 through 10:

\noindent \S 2.  Some basic notations

\noindent \S 3.  Some special probability measures on $\{-1,1\}^m$
with $m$ being powers of~6

\noindent \S 4.  A particular class of functions on certain infinite
sequences of vectors

\noindent \S 5. A special Markov chain (based on Section 4) and a
related random sequence

\noindent \S 6. Scaffolding: primarily a random field $(W_k^{(n)},\
n \in {\bf N},\ k \in {\bf Z})$ based on Section 5

\noindent \S 7. More scaffolding: including  the random sequence $X$
for Theorem 1.1

\noindent \S 8. More scaffolding: random sequences $(X_k^{(n)}, k
\in {\bf Z})$ and $(Y_k^{(n)}, k \in {\bf Z})$ for $n \in {\bf N}$

\noindent \S 9.  The proofs of most properties in Theorem 1.1

\noindent \S 10.  Proof of property (D) in Theorem 1.1
\bigskip
\vfill\eject

\centerline {\bf 2. Notations and conventions}  
\bigskip
This section gives some specific notations and conventions that will be used throughout this paper.

The cardinality of a set $S$ will be denoted $\card\, S$.

For a given probability space $(\Omega, \cF, P)$, the indicator
function (on $\Omega$) of a given event $A$ will be denoted $I(A)$,
and the $\sigma$--field $(\subset \cF)$ generated by a given
collection $(V_j,\ j\in J)$ of random variables (where $J$ is an index set) will be denoted $\sigma(V_j$, $j\in J)$.

A ``left-infinite sequence'' (of elements of some set) is a family of
elements $(x_k$, $k\le j)$ where $j$ is an integer (and $k$ is
restricted to integers).  Often $j=0$ in that context.
\medskip

{\sc Notations 2.1.} 
(A)~For typographical convenience, for a given nonnegative integer
$n$, whenever the integer $6^n$ appears in a subscript or
superscript, we shall use the notation
$$\hbox{\rm sxtp}(n) := 6^n.
\leqno{(2.1)}$$ (The letters ``sxtp'' are an abbreviation of the
word ``sextuple.'')

(B) For a given $n\in \N$ and a given vector $a:=
(a_1,a_2,\dots,a_n)\in {\R}^n$, the sum and product of the elements
will be denoted by
$$\hbox{\rm sum}(a) := \sum^n_{k=1} a_k\quad\hbox{and}\quad
\hbox{prod}(a) := \prod^n_{k=1}a_k. \leqno{(2.2)}$$

(C) Suppose $m\in \N$, and for each $u\in \{1,2,\dots,6\}$, $v_u:= (v_{u,1},v_{u,2}$, $\dots,v_{u,m})$ is a vector with 
$m$ coordinates. 
Then the notation $\langle v_1,v_2,v_3,v_4,v_5, \allowbreak
v_6\rangle$ will mean the vector $w:= (w_1,w_2,\dots,w_{6m})$ with $6m$ coordinates such that for each $u\in \{1,2,\dots,6\}$ and each $j\in \{1,2,\dots,m\}$, $w_{(u-1)m+j} = v_{u,j}$. That is, in the vector $\langle v_1,v_2,v_3,v_4,v_5,v_6\rangle$, with $6m$ coordinates, the first $m$ coordinates are (in order) the coordinates of $v_1$, the next $m$ coordinates are those of $v_2$, and so on, with the last $m$ coordinates being those of $v_6$.

(D) On the family of all nonempty finite subsets of $\Z$, the
following partial ordering will be used: For such sets $A$ and $B$,
the notation $A<B$ (or $B>A$) means that $\max A<\min B$
(equivalently, $a<b$ for all $a\in A$, $b\in B$).

(E) Suppose $S$ is a nonempty finite set of integers. If $(a_k$,
$k\in S)$ is a family of elements of some set $A$, then the notation
$a_S$ denotes the vector defined by
$$a_S := \left(a_{s(1)},a_{s(2)},\dots,a_{s(n)}\right)
\leqno{(2.3)}$$ where $n=$ card $S$ and $s(1)<s(2)< \dots<s(n)$ are
{\it in increasing order}\/ the elements of $S$. Similarly, if
$(X_k$, $k\in S)$ is a family of random variables on a probability
space $(\Omega,\cF,P)$, then the notation $X_S$ denotes the random
vector defined by
$$X_S:= \left(X_{s(1)},X_{s(2)},\dots,X_{s(n)}\right)
\leqno{(2.4)}$$ where the $s(j)$'s are as above. The vector $a_S$ in
(2.3) and the random vector $X_S$ in (2.4) will also be expressed
respectively as $(a_{s(i)}$, $1\le i\le n)$ and $(X_{s(i)}$, $1\le
i\le n)$.

(F) (i) A set $S\subset \Z$ is said to be ``doubly infinite'' if it
contains infinitely many negative integers and infinitely many
positive integers.

(ii) If $(a_k$, $k\in \Z)$ is a (``two-sided'') sequence of elements
of some set $A$, and $S$ is a doubly infinite subset of $\Z$, then
the notation $(a_j$, $j\in S)$ refers to the sequence $(a_{s(\ell)}$,
$\ell\in \Z)$ where
$$\dots < s(-2)<s(-1)<s(0) \le 0 < 1 \le s(1)
<s(2)<s(3)<\dots $$
and $S = \{\,\dots, s(-1),s(0),s(1),\dots\, \}$.

(G) If $(a_k$, $k\in \Z)$ is a (``two-sided'') sequence of elements
of some set $A$, and $S$ is a subset of $\Z$ which is {\it both
infinite and bounded above,} then a convention ``opposite'' to that
of sections (E) and (F) will be used: The notation $(a_j$, $j\in S)$
refers to the sequence $(a_{s(0)},a_{s(-1)}, a_{s(-2)},\dots\,)$
where $s(0)>s(-1)>s(-2)>\dots$ and $S= \{s(0),s(-1),s(-2),\dots\,\}$.

(H) Suppose $\cE = \{E_1,E_2,E_3,\dots\,\}$ is a nonempty family of
finitely many or countably many sets $E_i \subset \Z$. (The sets
$E_i$ themselves can be countable or finite, or even empty, and they
need not be disjoint.) Then define the notation
$$\hbox{union } \cE := E_1\cup E_2 \cup E_3 \cup \dots\,\,,$$
the union of all sets $\in \cE$.
\medskip

{\sc Definition 2.2.} 
(A) For each $n\in \{0,1,2,\dots\,\}$, we shall define a function
$\psi_n:\{0,1,2,3,4,5,6\}^{\N} \longrightarrow \{0,1,2,\dots\,\}$.
Suppose $a:= (a_0,a_1,a_2,\dots\,)$ is a sequence of elements of
$\{0,1,\dots,6\}$. If the set $\{k\ge 0:a_k=1\}$ is infinite, then define the nonnegative integers  
$\psi_0(a),\ \psi_1(a),\ \psi_2(a),\ \dots$ (uniquely) by the
conditions
$$0\le \psi_0(a) <\psi_1(a) <\psi_2(a)<\dots\quad\hbox{and}$$
$$\{k\ge 0:a_k=1\} = \{\psi_0(a),\psi_1(a),\psi_2(a),\dots\,\}.$$
If instead the set $\{k\ge 0: a_k =1\}$ is finite, then
$\psi_n(a):=0$ for all $n\ge 0$. (This last sentence will be an
irrelevant formality.)

(B) {\it Remark.\/} In that definition (both cases), if
$a_0=1$ then $\psi_0(a) = 0$.
\medskip

{\sc Definition 2.3.} 
(A) A (``two-sided'') sequence $w:= (w_k$, $k\in \Z)$ of elements of $\{0,1,2,3,4,5,6\}$ satisfies ``Condition~$\cS$'' if the following three statements hold:

(i) For each $i\in \{1,2,3,4,5,6\}$ the set $\{k\in {\Z}:w_k=i\}$ is doubly infinite. (That is not required for $i=0$.)

(ii) For each $i\in \{1,2,3,4,5\}$ and each $\ell\in \Z$ such that
$w_\ell =i$, either (a)~$w_{\ell+1} = i+1$ or (b)~for some $m\ge 2$,
one has that $w_{\ell+1} = w_{\ell+2} = \dots =w_{\ell+m-1} =0$ and
$w_{\ell+m} =i+1$.

(iii) For each $\ell\in\Z$ such that $w_\ell =6$, either
(a)~$w_{\ell+1} =1$ or (b)~for some $m\ge 2$, one has that
$w_{\ell+1} = w_{\ell+2} = \dots =w_{\ell+m-1} =0$ and $w_{\ell+m}
=1$.

That is, the sequence $w$ satisfies Condition $\cS$ if ``from time immemorial,'' the non-zero elements in $w$ cycle through
$1,2,\dots,6$ in order, over and over again, with perhaps some 0's in between.
\medskip

(B) {\it Remark.} Suppose $w:= (w_k$, $k\in \Z)$ is a sequence of elements of $\{0,1,\dots,6\}$ that satisfies 
Condition~$\cS$. 
Then for any two integers $J$ and $L$ such that $J\le L$, one has that
$$\leqalignno{
&\card \{k\in {\Z}: J\le k\le L\ \ \hbox{and}\ \ w_k=1\} &(2.5)\cr
& \indent
\le\, 1 + (1/6) \cdot \card \{k\in {\Z}: J\le k\le L\ \ 
\hbox{and}\ \  w_k\not= 0 \}.  \cr}$$ 
The point is that between any two ``consecutive 1's,'' there are (exactly) five non-zero elements ($2,3,4,5,6$, each once). Hence, if the left side of (2.5) equals $\ell$ for some $\ell \ge 2$, then the set in the right side of (2.5) has at least $\ell +5(\ell-1)$ elements, and thus (2.5) holds.
Equation (2.5) holds trivially when its left side is 0 or~1.
\medskip

(C) {\it Remark.} Suppose $S:= \{\, \dots,s(-1),s(0),s(1),\dots\,\}$
is a doubly infinite set of integers where $\dots <s(-1)<s(0)
<s(1)<\dots\,\,$. 
Suppose $w:= (w_k$, $k\in \Z)$ is a (``two-sided'') sequence of elements of $\{0,1,\dots,6\}$ such that
(i)~$w_k =0$ for all $k\in {\Z}-S$, and (ii)~the (``two-sided'')
sequence $(\, \dots,w_{s(-1)},w_{s(0)}, w_{s(1)},\dots\,)$ satisfies Condition~$\cS$. Then the entire sequence $w$ satisfies Condition~$\cS$.

(In words, if a given two-sided sequence satisfies 
Condition~$\cS$, and one sticks some zeros between
its entries, the resulting new sequence still satisfies
Condition~$\cS$.)  
\medskip

{\sc Notations 2.4.} 
(A) The $6\times 6$ identity matrix will be denoted~$I_6$.

(B) The transpose of a vector $v$ will be denoted~$v^t$.

(C) The elements of $\{0,1\}^6$ will be represented as ``row''
vectors.  Suppose $m$ is a positive integer, and for each $i\in
\{1,2,\dots,m\}$, $\alpha_i := (\alpha_{i1},\alpha_{i2}$,
$\dots,\alpha_{i6})$ is an element of $\{0,1\}^6$. Then for the
$6\times m$ matrix whose columns are the transposes of the
$\alpha_i$'s respectively, we shall use the following notation:
$$\left[\alpha^t_1\mid \alpha^t_2 \mid \cdots \mid \alpha^t_m\right]
:= \left[ \matrix{ \alpha_{11} & \alpha_{21} & \cdots &
\alpha_{m1}\cr \alpha_{12} &\alpha_{22} & \cdots &\alpha_{m2} \cr
\vdots & \vdots & \vdots\vdots\vdots & \vdots \cr \alpha_{16} &
\alpha_{26} & \cdots & \alpha_{m6}\cr} \right]\, . \leqno{(2.6)}$$

(D) A (``left-infinite'') sequence $(\beta_0, \beta_{-1},
\beta_{-2},\dots\,)$ of elements of $\{0,1\}^6$ is said to be
``back-standard'' if
$$\left[ \beta^t_{\ell-5}\mid \beta^t_{\ell-4} \mid \cdots \mid
\beta^t_{\ell}\right] = I_6 \leqno{(2.7)}$$
holds for infinitely
many integers $\ell \le 0$. A (``two-sided'') sequence $(\beta_k$,
$k\in \Z)$ of elements of $\{0,1\}^6$ is said to be ``two-sided
standard'' if the set $\{\ell\in\Z:$ (2.7) holds$\}$ is doubly
infinite (see Section~2.1(F)).

(E) In the proofs of lemmas, we shall often use the notation
$$\beta^\ell_j := \left(\beta_j,\beta_{j+1},\dots,\beta_\ell\right)
\leqno{(2.8)}$$
when $j\le \ell $ are integers and $\beta_k \in
\{0,1\}^6$ for each $k\in \{j,\dots,\ell\}$, and also the
(``reverse'') notation
$$\beta^\ell_{-\infty} := \left(\beta_\ell, \beta_{\ell-1},
\beta_{\ell-2},\dots\,\right) \leqno{(2.9)}$$
when $\ell \in \Z$ and
$\beta_k \in \{0,1\}^6$ for each $k\le \ell$.
\medskip

{\sc Notations 2.5.}
(A) Suppose that on a probability space $(\Omega,\cF,P)$, $\eta$ is
a random variable (or random vector, etc.) taking its values in a
measurable space $(A,\cA)$. The distribution (or ``law'') of $\eta$
(on $(A,\cA))$ will be denoted $\cL(\eta)$. For any event $F$ such
that $P(F)>0$, the conditional distribution of $\eta$ given $F$ will
be denoted $\cL(\eta\mid F)$.

(B) Suppose $m$ is a positive integer and $\lambda$ is a probability
measure on $\{-1,1\}^m$. Then $\lambda^{[6]}$ will denote the
six-fold ``product measure'' of $\lambda$. That is,  
$\lambda^{[6]}$ is the
distribution (on $\{-1,1\}^{6m})$ of a $\{-1,1\}^{6m}$--valued
random vector
$$\langle Z_1,Z_2,\dots,Z_6\rangle$$
(in Notation 2.1(C)) where $Z_1,Z_2,\dots,Z_6$ are six independent
$\{-1,1\}^m$--valued random vectors, each having
distribution~$\lambda$.

(C) Suppose $(\Omega,\cF,P)$ is a probability space. For two events
$A$ and $B$, the notation $A\, \dot = \, B$ will mean that
$P(A\triangle B) =0$, where $\triangle$ denotes the 
symmetric difference.  
If $A$ is an event and $\cB$ is a $\sigma$-field
$\subset \cF$, then
the notation $A\, \dot\in\, \cB$ will mean that there 
exists an event $B \in \cB$ such that $A \dot = B$.  
If $\cA$ and $\cB$ are $\sigma$--fields $\subset \cF$,
the notation $\cA\, \dot\subset \, \cB$ will mean that for every $A\in \cA$, one has that $A\, \dot\in\, \cB$.
\medskip

In arguments below, we shall sometimes simply show 
$\cA\dot\subset \cB$ with a quick
verification, in place of $\cA \subset \cB$, when that latter literal inclusion is not needed and requires a longer argument.
\medskip

{\sc Definition~2.6.} 
(A) Suppose that on a probability space $(\Omega,\cF,P)$, $Y:=
(Y_k$, $k\in \Z)$ is a sequence of random variables taking their
values in a measurable space $(A,\cA)$ and $Z:= (Z_k$, $k\in \Z)$ is
a sequence of random variables taking their values in a measurable
space $(B,\cB)$. The ordered pair $(Y,Z)$ is said to satisfy
``Condition $\cM$'' if the following two conditions are satisfied:

(i) The sequence $Z$ is strictly stationary.

(ii) There exists a measurable function $f:B^{\N}\to A$ (that is,
$f^{-1}(H)\in \cB^{\N}$ for every $H\in \cA$, for the ``infinite
product'' $\sigma$--field $\cB^{\N}:= \cB \times \cB \times \cB\times
\dots\,$) such that for each $k\in \Z$, $Y_k =
f(Z_k,Z_{k-1},Z_{k-2},\dots\,)$ a.s.
\medskip

(B) {\it Remark.} Obviously, if $(Y,Z)$ satisfies Condition~$\cM$,
then (i)~the sequence $Y$ is strictly stationary, and
(ii)~$\sigma(Y)\, \dot\subset \, \sigma(Z)$ (see section~2.5(C)).
\medskip

{\sc Notations 2.7.} Suppose $(\Omega,\cF,P)$ is a probability space.

(A) An ordered triplet $(\cA,\cB,\cC)$ of $\sigma$--fields $\subset
\cF$ is a ``Markov triplet'' if one has that for all $A\in \cA$ and
all $C\in \cC$, $P(A\cap C|\cB) = P(A|\cB) \cdot P(C|\cB)$ a.s.

A ``restricted'' version of (A) will also be needed for some random
sequences that are not Markov chains but have some ``limited'' Markov
properties:

(B) Suppose $B\in \cF$, and suppose $\cA$ and $\cC$ are
$\sigma$--fields $\subset \cF$. The ordered triplet $(\cA,B, \cC)$
is a ``restricted Markov triplet'' if (i)~$P(B)>0$, and (ii)~for all
$A\in \cA$ and all $C\in \cC$, $P(A\cap C|B) = P(A|B) \cdot P(C|B)$.

(C) {\it Remark.} If $(\cA,B,\cC)$ is a restricted Markov triplet,
$A\in \cA$, $P(A\cap B)>0$, and $C\in \cC$, then (by a trivial
calculation) $P(C|A\cap B) = P(C|B)$.

(D) {\it Remark.} Suppose $\cG$ and $\cH$ are independent
$\sigma$--fields $\subset \cF$. Suppose $G\in \cG$, and $\cG_1$ and
$\cG_2$ are $\sigma$--fields $\subset \cG$, and $(\cG_1,G,\cG_2)$ is
a restricted Markov triplet.  Suppose $H\in \cH$, and $\cH_1$ and
$\cH_2$ are $\sigma$--fields $\subset \cH$, and $(\cH_1,H,\cH_2)$ is
a restricted Markov triplet. Then $(\cG_1\vee \cH_1$, $G\cap H$,
$\cG_2\vee \cH_2)$ is a restricted Markov triplet.
\medskip

{\it Proof of (D).} Of course $P(G\cap H) = P(G) \cdot P(H)>0$. By a
simple calculation, if $G_1 \in \cG_1$, $G_2\in \cG_2$, $H_1\in
\cH_1$, and $H_2\in \cH_2$, then
$$\eqalign{
&P(G_1\cap G_2\cap H_1 \cap H_2\mid G\cap H)= P(G_1\cap G_2\mid G)
\cdot P(H_1 \cap H_2\mid H) \cr
&=P(G_1|G) \cdot P(G_2|G) \cdot
P(H_1|H) \cdot P(H_2|H) \cr 
&= P(G_1 \mid G\cap H) 
\cdot P(G_2 \mid G\cap H)
\cdot P(H_1 \mid G\cap H)
\cdot P(H_2 \mid G\cap H). \cr}
$$
Thus under the probability measure  $Q(\, . \, ) := P(\, . \, |G\cap
H)$ on $(\Omega,\cF)$, the four $\sigma$--fields $\cG_1$, $\cG_2$,
$\cH_1$, and $\cH_2$ are independent, and hence (under $Q$) the
$\sigma$--fields $\cG_1\vee \cH_1$ and $\cG_2\vee \cH_2$ are
independent. Thus (D) holds.
\medskip

{\sc Remark 2.8} (a standard trivial but useful fact). Suppose
$(\Omega,\cF,P)$ is a probability space, $C = \bigcup_iC_i$ where
$C_1,C_2,C_3,\dots$ is a finite or countable sequence of (pairwise)
disjoint events, $P(C_i)>0$ for at least one $i$, $D$ is an event,
$p\in [0,1]$, and $P(D|C_i) =p$ (resp.\ $\ge p$ resp.\ 
$\le p$) for every $i$ such
that $P(C_i)>0$. Then $P(D|C) =p$ (resp.\ $\ge p$ 
resp.\ $\le p$).

\bigskip

\centerline {\bf 3. Some key probability measures} 
\bigskip
This section is devoted to the definitions and key properties of certain discrete probability measures that will play a pervasive role in the construction for Theorem~1.1. 

In connection with the design of error-correcting codes, the book by McWilliams and Sloane [20] includes an extensive treatment of constructions, for a given positive integer $k$, of ``big''$k$-tuplewise independent random vectors from ``small'' ones.  
The material here in section~3 fits into that general framework (in a somewhat hidden way).  
However, the details will need to be spelled out here, in order to facilitate the proofs of certain properties
(in particular, the bounds on sixth moments of partial sums) in Theorem~1.1.

The main tool in this section is a well known, elementary ``parity'' trick, built into Definitions 3.1 and 3.3 below.  The same general type of ``parity'' trick was used (with different distributions) by Pruss [28] and by the
author and Pruss [6] in the constructions in those papers, and was also used (for different purposes) in the book [20] alluded to above.

In this section, for convenience, for a given $n\in\N$ and a given
element $x\in\{-1,1\}^{{\rm sxtp}(n)}$ (see (2.1)), $x$ will be
represented as $(x_0,x_1,\dots,x_{{\rm sxtp}(n)-1})$ (instead of
$x_1,x_2,\dots,x_{{\rm sxtp}(n)})$).
\medskip

{\sc Definition~3.1.} 
Referring to (2.2),  
 define the set
$$\Upsilon :=\{x:= (x_0,x_1,\dots,x_5)\in \{-1,1\}^6: \hbox{prod} (x)
=-1\}.$$
That is, $\Upsilon$ is the set of all 6-tuples of $-1$'s and
 $+1$'s with an odd number of $-1$'s. For a given $x:=
(x_0,x_1,\dots,x_5)\in \Upsilon$, sum$(x)$ (see
(2.2)) 
is $-4$ (resp.\ 0 resp.\ 4) if exactly 5 (resp.\ 3
resp.\ 1) of the $x_i$'s are $-1$.
\medskip

Let $\nu^{\key}$ denote the uniform probability measure on
$\Upsilon$. That is, \break
\noindent $\nu^{\key}(\{x\}) = 1/32$ for
each $x\in \Upsilon$.
\medskip

{\sc Remark 3.2.} 
If $V:= (V_0,V_1,\dots,V_5)$ is a
$\{-1,1\}^6$--valued random vector with distribution $\nu^{\key}$ (thus $V\in \Upsilon$ a.s.), then by trivial arguments, the following statements hold:

(A) The distribution of the random vector $-V$ is $\nu^{\key}$.

(B) For any permutation $\sigma$ of $\{0,1,\dots,5\}$, the
distribution of the random vector
$(V_{\sigma(0)},V_{\sigma(1)},\dots,V_{\sigma(5)})$ is $\nu^{\key}$.

(C) For each $i\in \{0,1,\dots,5\}$, $P(V_i=-1) =P(V_i=1) =1/2$.

(D) For every set $S\subset \{0,1,\dots,5\}$ such that card~$S=5$,
the random variables $V_i$, $i\in S$ are independent.

(E) prod$(V) =-1$ a.s. (see (2.2)). 

(F) $P($sum$(V)=0) =5/8$ and $P($sum$(V)=-4) =P($sum$(V)=4) = 3/16$
(see (2.2)). 

(G) For each $i\in \{0,1,\dots,5\}$, $P(V_i =1|$sum$(V) =4) =5/6$ and
$P(V_i=-1|$sum$(V) =4) = 1/6$, and hence $E(V_i|$sum$(V) =4) =2/3$.
\medskip

{\sc Definition 3.3.} 
Refer to (2.1) 
and (2.2).
For each $n\in \N$, we shall define four probability
measures $\nu^{(n)}_{\rm ord}$, $\nu^{(n)}_{\rm cen}$, $\nu^{(n)}_{\rm
fri}$, and $\nu^{(n)}_{\rm pos}$ on the set \break
\noindent $\{-1,1\}^{\sxtp(n)}$, such that the following holds:
$$\eqalign{&\hbox{\rm
If}\,\,W:=\left(W_0,W_1,\dots,W_{\sxtp(n)-1}\right)\,\, \hbox{\rm is
a $\{-1,1\}^{\sxtp(n)}$--valued} \cr
&\hbox{random vector with
distribution $\nu^{(n)}_{\rm pos}$, then sum$(W) =4^n$ a.s.}\cr} \leqno(3.1)$$
(The subscripts ``ord'', ``cen,'' ``fri,'' and ``pos'' are
respectively abbreviations of the words ``ordinary'' ``center,''
``fringe,'' and ``positive.'' The relevance of those subscripts will
become clear below.) For $n=0$, we shall also define two probability
measures $\nu^{(0)}_{\rm fri}$ and $\nu^{(0)}_{\rm pos}$ on the set
$\{-1,1\}$, such that (3.1) holds for $n=0$. The definition is
recursive and is as follows:

Start with $n=0$. On the set $\{-1,1\}$, define the probability
measures $\nu^{(0)}_{{\rm fri}}$ and $\nu^{(0)}_{{\rm pos}}$ by
$\nu^{(0)}_{{\rm fri}}(\{-1\}) = \nu^{(0)}_{{\rm fri}}(\{1\}) =1/2$
and $\nu^{(0)}_{{\rm pos}}(\{1\})=1$. Trivially (3.1) holds for
$n=0$.

Now suppose $n\ge 0$, and suppose the probability measure
$\nu^{(n)}_{{\rm pos}}$ on $\{-1,1\}^{\sxtp(n)}$ has already been
defined such that (3.1) holds. Let $W^{(i)}:=$ $(W^{(i)}_0,
W^{(i)}_1,\dots,W^{(i)}_{\sxtp(n)-1})$, $i\in \{0,1,\dots,5\}$ be
six independent \qquad\break \noindent $\{-1,1\}^{\sxtp(n)}$--valued
random vectors, each having distribution $\nu^{(n)}_{\rm{pos}}$. Let
$V:=(V_0,V_1,\dots,V_5)$ be a $\{-1,1\}^6$-valued random vector
which is independent of the family $(W^{(i)}$, $0\le i\le 5)$ and
has distribution $\nu^{\key}$ (see Definition~3.1). 
Let $Z:= (Z_0,Z_1,\dots,Z_{\sxtp(n+1)-1})$ be the
$\{-1,1\}^{\sxtp(n+1)}$-valued random vector defined as follows:
$$\eqalign{\forall\,\, i\in \{0,1,\dots,5\},&\quad \forall\,\, j\in
\{0,1,\dots,6^n-1\}, \cr Z_{i\cdot \sxtp(n)+j}&:=V_i\cdot
W^{(i)}_j.\cr} \leqno{(3.2)}$$ 
Then (see (2.2), (3.1), 
and (3.2)) 
with probability~1,
$$\hbox{sum}(Z) = \sum^5_{i=0} \sum^{\sxtp(n)-1}_{j=0} V_iW_j^{(i)} =
(\hbox{sum}(V))\cdot 4^n, \leqno{(3.3)}$$ 
and hence (see Remark~3.2(F)) 
sum$(Z)$ takes the
value 0 resp.\ $-4^{n+1}$ resp.\ $4^{n+1}$ with probability 5/8
resp.\ 3/16 resp.\ 3/16.

Let $\nu^{(n+1)}_{\rm ord}$ denote the distribution on
$\{-1,1\}^{\sxtp(n+1)}$ of the random vector $Z$. Let
$\nu^{(n+1)}_{\rm cen}$ resp.\ $\nu^{(n+1)}_{\rm{fri}}$ 
resp.\
$\nu^{(n+1)}_{\rm{pos}}$ denote the conditional distribution on
$\{-1,1\}^{\sxtp(n+1)}$ of the random vector $Z$ given the event
$\{$sum$(Z)=0\}$ resp.\ $\{|$sum$(Z)| =4^{n+1}\}$ resp.\
$\{$sum$(Z)=4^{n+1}\}$. Thus (3.1) 
holds with $n$
replaced by $n+1$. This completes the recursive definition.
\medskip

{\sc Remark 3.4.} 
For each $n\in \N$, $\nu^{(n)}_{\rm ord}
=(5/8)\nu^{(n)}_{\rm{cen}} +(3/8)\nu^{(n)}_{\rm{fri}}$.
\medskip

{\sc Proof.} If $Z$ is a $\{-1,1\}^{\sxtp(n)}$--valued random vector
with distribution $\nu^{(n)}_{\rm ord}$, and $B\subset
\{-1,1\}^{\sxtp(n)}$, then (see the comments after (3.3),
but with $n+1$ replaced by $n$)
$$\eqalign{
\nu^{(n)}_{\rm ord}(B) &= P(Z\in B) \cr &= P(Z\in B|
\hbox{sum}(Z)=0)\cdot P(\hbox{sum}(Z)=0) \cr &\qquad +P(Z\in B|\,
|\hbox{sum}(Z)|=4^n)\cdot P(|\hbox{sum}(Z)|=4^n) \cr &=
\nu^{(n)}_{\rm{cen}}(B) \cdot (5/8) + \nu^{(n)}_{\rm{fri}}(B) \cdot
(3/8). \cr}$$

\medskip

The various distributions in Definition~3.3 
have pervasive symmetries. 
We will need later on, and will verify in the next two sections, only a couple of mild aspects or manifestations of those symmetries.
\medskip

{\sc Remark 3.5.} 
(A)~For a given $n\ge 0$, if the random
vectors $W^{(i)}$, $i\in \{0,1,\dots,5\}$, $V$, and $Z$ are as in
Definition~3.3, 
then the distribution of the random collection $(V,W^{(0)},
W^{(1)},\dots,W^{(5)})$ is (by independence) a product measure on
$\{-1,1\}^6 \times (\{-1,1\}^{\sxtp(n)})^6$, and (by Remark~(3.2)(A))
is the same as that of
$(-V,W^{(0)},W^{(1)},\dots,W^{(5)})$, and hence by (3.2)
the distribution on $\{-1,1\}^{\sxtp(n+1)}$ of the
random vector $-Z$ is the same as that of $Z$.

(B)~By an elementary argument, it
follows that for any given $n\in \N$, the distributions
$\nu^{(n)}_{\rm ord}$, $\nu^{(n)}_{\rm{cen}}$, and
$\nu^{(n)}_{\rm{fri}}$ on $\{-1,1\}^{\sxtp(n)}$ satisfy the symmetry
conditions $\nu^{(n)}_{\rm ord}(\{-x\})=\nu^{(n)}_{\rm ord}(\{x\})$,
$\nu^{(n)}_{\rm{cen}}(\{-x\}) = \nu^{(n)}_{\rm{cen}}(\{x\})$, and
$\nu^{(n)}_{\rm{fri}}(\{-x\}) = \nu^{(n)}_{\rm{fri}}(\{x\})$. For
$n=0$, this holds trivially for $\nu^{(0)}_{\rm{fri}}$.

(C) Here for convenient reference are a few other features and
elementary consequences of Definition 3.3, for a given $n\ge 0$ and
a given $x\in \{-1,1\}^{\sxtp (n)}$: (i)~If sum$(x) = 4^n$ then
$\nu^{(n)}_{\Fri}(\{-x\}) = \nu^{(n)}_{\Fri}(\{x\}) = (1/2)
\nu^{(n)}_{\pos}(\{x\})$. (ii)~If $\nu^{(n)}_{\pos}(\{x\})>0$ then
sum$(x) = 4^n$. (iii)~If $\nu^{(n)}_{\Fri}(\{x\})>0$ then sum$(x) =
-4^n$ or $4^n$. (iv)~If $(n\ge 1$ and) $\nu^{(n)}_{\cen}(\{x\})>0$
then sum$(x) =0$. (v)~If ($n\ge 1$ and) $\nu^{(n)}_{\ord}(\{x\})>0$
then sum$(x) = -4^n$, $0$, or $4^n$.
\medskip

{\sc Lemma 3.6.} 
{\it Suppose $n\ge 0$. Suppose $Y:= (Y_0,Y_1,\dots,Y_{\sxtp(n)-1})$
is a $\{-1,1\}^{\sxtp(n)}$--valued random vector whose distribution
is  $\nu ^{(n)}_{\pos}$. Then $EY_k = (2/3)^n$ for every $k\in
\{0,1,\dots,6^n-1\}$.}
\medskip

{\sc Proof.} Lemma~3.6 
holds trivially for $n=0$. Now for induction, suppose it holds for a
given $n\ge 0$. Let the random vectors $W^{(i)}$, $i\in
\{0,1,\dots,5\}$, $V$, and $Z$ be as in Definition~3.3.
By (3.3),  
the events $\{$sum$(Z) =4^{n+1}\}$ and
$\{$sum$(V) =4\}$ are identical (modulo a null set).
By (3.2), 
Remark~3.2(G), 
and our induction
hypothesis, for each $i\in \{0,1,\dots,5\}$ and each  \break
\noindent $j\in \{0,1,\dots,6^n-1\}$,
$$\eqalign{
E\Bigl(&Z_{i\cdot \sxtp (n)+j}\, \Big|\, \hbox{sum}(Z)=4^{n+1}\Bigr) =
E\left(V_i\cdot W^{(i)}_j\, \Big|\, \hbox{sum}(V)=4\right) \cr
&= EW^{(i)}_j
\cdot E(V_i|\hbox{sum}(V)=4) = (2/3)^n \cdot (2/3) = (2/3)^{n+1}.
\cr}$$
That is, $E(Z_k| \hbox{sum}(Z) =4^{n+1}) = (2/3)^{n+1}$ for
every $k\in \{0,1,\dots,6^{n+1}-1\}$. By Definition~3.3 
itself, Lemma~3.6 
holds for $n+1$. That completes the
induction step and the proof.
\medskip

{\sc Lemma 3.7.} 
{\it Suppose $n\ge 0$. Suppose
$Y:=(Y_0,Y_1,\dots,Y_{\sxtp(n+1)-1})$ is a
$\{-1,1\}^{\sxtp(n+1)}$--valued random vector with distribution
$(\nu^{(n)}_{\rm fri})^{[6]}$ (see section~2.5(B)). Suppose
$Z:=(Z_0,Z_1,\dots,Z_{\sxtp(n+1)-1})$ is a
$\{-1,1\}^{\sxtp(n+1)}$--valued random vector with the distribution
$\nu^{(n+1)}_{\ord}$. Then the following three statements hold:

(A) For every set $S\subset \{0,1,\dots,6^{n+1}-1\}$ such that {\rm
card}~$S=5$, the random vectors $Y_S$ and $Z_S$ (see (2.4)) have the
same distribution on $\{-1,1\}^5$.

(B) For every nonempty set $Q\subset \{0,1,\dots,6^{n+1}-1\}$, one
has that
$$E\biggl(\, \sum_{k\in Q}Y_k\biggr)^6 \ge 
E\biggl(\, \sum_{k\in Q}Z_k\biggr)^6.
\leqno{(3.4)}$$

(C) Also, }
$$E(\hbox{\rm sum}(Y))^6 = E(\hbox{\rm sum}(Z))^6+ 720\cdot 4^{6n}.
\leqno{(3.5)}$$

{\sc Proof.} Suppose $n\ge 0$. Without loss of generality, to prove
Lemma~3.7, 
we shall construct particular convenient
random vectors $Y$ and $Z$ with the required distributions, and then
prove statements (A), (B), and (C) for those two random vectors.

Let the random vectors $W^{(i)}$, $i\in \{0,1,\dots,5\}$, $V$, and
$Z$ be as in Definition~3.3. 
In that context, let $U_i$,
$i\in \{0,1,\dots,5\}$ be independent, identically distributed
$\{-1,1\}$--valued random variables with $P(U_i=-1) =P(U_i=1) = 1/2$,
with the family $(U_i$, $0\le i\le 5)$ being independent of the
family $(V;\ W^{(i)},\ 0\le i\le 5)$. Define the
$\{-1,1\}^{\sxtp(n+1)}$--valued random vector
$Y:=(Y_0,Y_1,\dots,Y_{\sxtp(n+1)-1})$ as follows:
$$\leqalignno{\forall\,\, i\in \{0,1,\dots,5\}, &\quad \forall\,\, j\in
\{0,1,\dots,6^n-1\}, &(3.6)\cr
Y_{i\cdot \sxtp(n)+j} &:= U_i\cdot
W^{(i)}_j.  \cr}$$
By Remark~3.5(B)(C), 
Definition~3.3, 
and a simple argument, for each \break
\noindent $i\in \{0,1,\dots,5\}$, the distribution on
$\{-1,1\}^{\sxtp(n)}$ of the
random vector $(U_i\cdot W_j^{(i)}$, $0\le j\le 6^n-1)$ (see the
sentence after (2.4)) 
is $\nu^{(n)}_{\rm{fri}}$. From this
and the definition of $\nu^{(n+1)}$ in Definition~3.3, the random vectors $Y$ and $Z$ constructed here have the
distributions specified in the statement of Lemma~3.7.
 
Just for this proof, define for each $i\in \{0,1,\dots,5\}$ the set
$$K_i = K^{(n)}_i := \{k\in {\Z}: i\cdot 6^n \le k\le (i+1)\cdot 6^n-1\}. \leqno{(3.7)}$$ 
These sets $K_0,K_1,\dots,K_5$ form a
partition of the set $\{0,1,\dots,6^{n+1}-1\}$. As a consequence of
Remark~3.2(C)(D) 
and the above conditions, for any set $\Lambda\subset
\{0,1,\dots,5\}$ with card~$\Lambda=5$, the random family $(V_i$,
$i\in \Lambda$; $W^{(i)}$, $i\in \Lambda)$ has the same
distribution---a ten-fold product measure on $\{-1,1\}^5
\times(\{-1,1\}^{\sxtp(n)})^5$ --- as the random family $(U_i$,
$i\in \Lambda$; $W^{(i)}$, $i\in \Lambda)$. Hence by (3.2) and
(3.6), for every set $\Lambda \subset \{0,1,\dots,5\}$ with
card~$\Lambda=5$ (see section 2.5(A) and the sentence after (2.4)),
$${\cal L}\biggl(Y_k, \, k\in \bigcup_{i\in \Lambda}K_i\biggr) = {\cal
L}\biggl(Z_k,\,  k\in
\bigcup_{i\in \Lambda}K_i\biggr)
\leqno(3.8) $$ 
(an equality of distributions on $\{-1,1\}^{5\cdot \sxtp(n)})$.

Now if $S\subset \{0,1,\dots,6^{n+1}-1\}$ is such that card~$S=5$,
then $S\subset \bigcup_{i\in \Lambda}K_i$ for some set $\Lambda\subset
\{0,1,\dots,5\}$ with card~$\Lambda\le 5$. Hence statement (A) in
Lemma~3.7 
follows from (3.8). 

{\it Proof of statements (B) and (C).} For any nonempty set $Q\subset
\{0,1,\dots$, $6^{n+1}-1\}$,
$$\leqalignno{\qquad E\biggl(&\sum_{k\in Q}Y_k\biggr)^6 - E\biggl(\sum_{k\in
Q}Z_k\biggr)^6 &(3.9)\cr & = \sum_{k(0)\in Q} \sum_{k(1)\in Q}
\cdots \sum_{k(5)\in Q}\left[ E\left(\prod^5_{i=0} Y_{k(i)}\right) -
E\left(\prod^5_{i=0} Z_{k(i)}\right)\right]. \cr}$$

For any set $\Lambda \subset \{0,1,\dots,5\}$ with card~$\Lambda=5$,
and any choice of (not necessarily distinct) elements $k(0)$,
$k(1),\dots$, $k(5)$ of $\bigcup_{i\in \Lambda}K_i$, the term in the
brackets in (3.9) 
 equals $0$ by (3.8). 
On the
other hand, if $k(i)\in K_i$ for every $i\in \{0,1,\dots,5\}$, then
representing $k(i) = i\cdot 6^n +j(i)$ where $j(i) \in
\{0,1,\dots,6^n-1\}$, one has
$$E\left(\prod^5_{i=0} Y_{k(i)}\right) = \left[\prod^5_{i=0}
EU_i\right]\cdot\left[\prod^5_{i=0} EW^{(i)}_{j(i)}\right]=0$$
by (3.6) 
and the trivial fact $EU_i=0$, and one also has
$$E\left(\prod^5_{i=0}Z_{k(i)}\right) =
\left[E\left(\prod^5_{i=0}V_i\right)\right] \cdot \left[
\prod^5_{i=0} EW_{j(i)}^{(i)}\right] = -1 \cdot [(2/3)^n]^6
=-(2/3)^{6n}$$
by (3.2), 
Remark~3.2(E), and 
Lemma~3.6. 
Hence by (3.9) 
for any nonempty set
$Q\subset \{0,1,\dots,6^{n+1}-1\}$,
$$[\hbox{Left side of (3.9)}] = (2/3)^{6n} \cdot 6! \cdot
\prod^5_{i=0} \hbox{card}(Q\cap K_i). \leqno{(3.10)}$$ 
Statement
(B) in Lemma~3.7 
 follows. Since card~$K_i=6^n$ for each
$i$, statement (C) follows from (3.10) 
with $Q =\{0,1,\dots,6^{n+1}-1\}$ itself. That completes the proof.
\bigskip

\centerline {\bf 4. Some particular functions}
\bigskip 

A key ``building block''
in the construction of the random sequence $X$ in Theorem~1.1 will be a particular strictly stationary, finite-state, irreducible, aperiodic Markov chain.
It will be given a particular, explicit representation as a
``causal moving function'' (\`a la eq.\ (1.3)) 
of an i.i.d.\ finite-state sequence.
This representation will be spelled out in detail in
Section~4 here and Section~5 together (the Markov chain
itself will be identified in Lemma 5.3), in order to facilitate transparent proofs (in Sections 7 and 10 respectively) of property (C) and (especially)
property (D) in Theorem~1.1.
The particular form of the representation that will be
used here is an old one that goes back several decades;
its origin is hard to trace.
Its spirit goes back at least to a paper of Rosenblatt [30]
in which it is shown that every strictly stationary, countable-state, irreducible, aperiodic Markov chain can be represented as a ``causal moving function'' of an i.i.d.\ sequence.
\medskip

{\sc Definition 4.1.} For each $n\in \N$, we shall define a function
$h_n: (\{0,1\}^6)^n \to \{1,2,3,4,5,6\}$. The definition will be
recursive and is as follows:

First, define the function $h_1: \{0,1\}^6 \to \{1,2,\dots,6\}$ as
follows: For $\alpha:= (\alpha_1,\alpha_2,\dots,\alpha_6)
\in\{0,1\}^6$,
$$h_1(\alpha) := \cases{
6 &if $\alpha_1=0$ \cr
1 &if $\alpha_1 =1$. \cr}
\leqno(4.1.1)$$

Now suppose $n\ge 2$, and the function $h_{n-1} :(\{0,1\}^6)^{n-1}
\to \{1,2,\dots,6\}$ has already been defined. Define the function
$h_n:(\{0,1\}^6)^n \to \{1,2,\dots,6\}$ as follows:

Suppose that for each $i\in \{1,2,\dots,n\}$, $\beta_i :=
(\beta_{i,1},\beta_{i,2},\dots,\beta_{i,6})\in \{0,1\}^6$. If
$h_{n-1}(\beta_1,\beta_2,\dots,\beta_{n-1}) = j\in \{1,2,3,4,5\}$,
then
$$h_n(\beta_1,\beta_2,\dots,\beta_n) := \cases{
j &if $\beta_{n,j+1}=0$ \cr j+1 & if $\beta_{n,j+1} =1$. \cr}
\leqno{(4.1.2)}$$ If instead
$h_{n-1}(\beta_1,\beta_2,\dots,\beta_{n-1})=6$, then
$$h_n(\beta_1,\beta_2,\dots,\beta_n) := \cases{
6 & if $\beta_{n,1}=0$ \cr
1 & if $\beta_{n,1} =1$. \cr}
\leqno{(4.1.3)}$$
That completes the recursive definition of the functions~$h_n$.
\medskip

{\sc Lemma 4.2.} {\it Suppose $\ell\in \{-5,-6,-7,\dots\,\}$, and
$\beta_\ell,\beta_{\ell+1},\dots,\beta_{-1}$ are each an element of
$\{0,1\}^6$, and that
$$\left[\beta^t_{-5}\mid \beta^t_{-4} \mid \beta^t_{-3} \mid
\beta^t_{-2} \mid \beta^t_{-1}\right] = \left[ \matrix{ 1 & 0 & 0 &
0 & 0 \cr 0 & 1 & 0 & 0 & 0 \cr 0 & 0 & 1 & 0 & 0 \cr 0 & 0 & 0 & 1
& 0 \cr 0 & 0 & 0 & 0 & 1 \cr 0 & 0 & 0 & 0 & 0 \cr}
\right]\leqno{(4.2.1)}$$ (the first five columns of the identity
matrix $I_6$ --- see Section~2.4(A)(B)(C)). Then }
$$h_{-\ell}(\beta_\ell,\beta_{\ell+1},\dots,\beta_{-1}) =5.
\leqno{(4.2.2)}$$
\medskip

{\sc Proof.} Consider first the case $\ell =-5$. By (4.2.1),
(4.1.1), and then four applications of (4.1.2), one obtains
$h_1(\beta_{-5})=1$, $h_2(\beta_{-5},\beta_{-4}) =2$,
$h_3(\beta_{-5},\beta_{-4},\beta_{-3}) =3$, $h_4(\beta^{-2}_{-5})
=4$ (we start using the notations in Section~2.4(E)), and finally
$h_5(\beta^{-1}_{-5}) =5$, which is (4.2.2).

Now consider the case where $\ell \le -6$. We shall give the argument
for the case where
$$h_{-5-\ell}(\beta_{\ell}, \beta_{\ell+1},\dots,\beta_{-6}) =3.
\leqno{(4.2.3)}$$ The argument is similar for the other possible
values $(1,2,4,5$, and $6)$ for the right hand side of (4.2.3).

Now $\beta_{-5,4} =0$ by (4.2.1), and hence
$h_{-4-\ell}(\beta^{-5}_\ell) =3$ by (4.2.3) and (4.1.2). Next,
$\beta_{-4,4} =0$ by (4.2.1), and hence $h_{-3-\ell}(\beta^{-4}_\ell
)=3$ now follows from (4.1.2). Similarly $\beta_{-3,4}=0$ and hence
$h_{-2-\ell}(\beta^{-3}_\ell) =3$. Next, $\beta_{-2,4} =1$ by
(4.2.1), and hence $h_{-1-\ell}(\beta^{-2}_\ell) =4$ by (4.1.2).
Finally, $\beta_{-1,5} =1$, and hence (4.2.2) now follows from
(4.1.2). That completes the proof of Lemma 4.2.  
\medskip

In Definitions 4.3 and 4.5 below, two closely related functions on \break
$(\{0,1\}^6)^{\N}$ will be defined. Because of the way those
functions will be used, ``left-infinite'' sequences will be used in
their definitions.
\medskip

{\sc Definition 4.3.} Define the function 
$g_{\rm basic}
(\{0,1\}^6)
^{\N} 
\to \{1,2,3,4,5,6\}$ as follows:

Suppose $\beta_0,\beta_{-1},\beta_{-2},\dots$ each $\in \{0,1\}^6$.

If the sequence $(\beta_0,\beta_{-1},\beta_{-2},\dots\,)$ is not back-standard (see Section \break 
2.4(A)(B)(C)(D)), then define $g_{\rm basic}
(\beta_0,\beta_1,\beta_{-2},\dots\,) := 6$.

Now suppose instead that the sequence
$(\beta_0,\beta_{-1},\beta_{-2},\dots\,)$ is back- \break
standard (again
see Section~2.4). Let $L$ denote the greatest integer $\in
\{0,-1,\allowbreak -2,\dots\,\}$ such that 
$[\beta^t_{L-5} \mid
\beta^t_{L-4}\mid \dots \mid \beta^t_{L}] = I_6$.

(i) If $L=0$ (that is, if $[\beta^t_{-5}\mid \beta^t_{-4}\mid \dots
|\beta^t_0] =I_6)$, then define $g_{\basic}(\beta_0$, $\beta_{-1}$,
$\beta_{-2}$, $\dots\,):=6$.

(ii) If instead $L\le -1$, then referring to Definition~4.1, define
$$g_{\basic}(\beta_0,\beta_{-1},\beta_{-2},\dots\,) :=
h_{-L}(\beta_{L+1},\beta_{L+2},\dots,\beta_0). \leqno{(4.3.1)}$$
\medskip

{\sc Lemma 4.4.} {\it Suppose
$(\beta_0,\beta_{-1},\beta_{-2},\dots\,)$ is a (``left-infinite'') back- \break
standard sequence of elements of $\{0,1\}^6$ (see
Section~2.4(D)).

(A) If $g_{\basic} (\beta_{-1},\beta_{-2},\beta_{-3},\dots\,) = j\in
\{1,2,3,4,5\}$, then
$$g_{\basic}(\beta_0,\beta_{-1},\beta_{-2},\dots\,) = \cases{
j & if $\beta_{0,j+1}=0$ \cr
j+1 & if $\beta_{0,j+1}=1$. \cr}
\leqno{(4.4.1)}$$

(B) If instead $g_{\basic}(\beta_{-1},\beta_{-2},\beta_{-3},\dots\,)
=6$, then}
$$
g_{\basic} (\beta_0,\beta_{-1},\beta_{-2},\dots\,) = \cases{
6 & if $\beta_{0,1}=0$ \cr
1 & if $\beta_{0,1}=1$. \cr}
\leqno{(4.4.2)}$$
\medskip

{\sc Proof.} Let $L$ denote the greatest {\it negative}\/ integer
(the integer 0 is excluded here) such that (see 
Section 2.4(A)(B)(C))
$$\left[ \beta^t_{L-5} \mid \beta^t_{L-4}\mid \dots \mid
\beta^t_L\right] = I_6. \leqno{(4.4.3)}$$

The proofs of statements (A) and (B) in Lemma~4.4 will be handled
together, and will be divided into three cases.  The notations in
Section 2.4(E) (both eqs.\ (2.8) and (2.9)) will be used.

{\it Case 1:} $[\beta^t_{-5}\mid \beta^t_{-4}\mid \cdots \mid
\beta^t_{0}] = I_6$. Then $\beta_{\ell,6} =0$ for $\ell \in
\{-5,-4$, $\dots,-1\}$, and hence by (4.4.3) (which implies
$\beta_{L,6}=1$) and its entire sentence, $L\le -6$ must hold. By
Definition~4.3 and Lemma~4.2, $g_{\basic} (\beta^{-1}_{-\infty}) =
h_{-L-1}(\beta^{-1}_{L+1})=5$. Hence here in Case~1, the hypothesis
of statement (A) (in Lemma~4.4) holds with $j=5$ there, and
statement (B) there is vacuous. Since (by the definition of $I_6$)
$\beta_{0,6} =1$, the right side of (4.4.1) equals 6. Also, by
Definition 4.3(i), the left side of (4.4.1) equals~6. Hence (4.4.1)
holds, and statement (A) is verified. That completes the argument
for Case~1.

{\it Case 2:} $[\beta^t_{-5}\mid \beta^t_{-4}\mid \dots \mid
\beta^t_0] \not= I_6$  and $L=-1$ (redundant --- see (4.4.3)).  Then
by Definition~4.3(i) (see (4.1.3)),
$g_{\basic}(\beta^{-1}_{-\infty})=6$, and hence the hypothesis of
statement (B) (in Lemma~4.4) is satisfied, and statement (A) there
is vacuous.  By Definition~4.3(ii) and Definition~4.1,
$g_{\basic}(\beta^0_{-\infty}) = h_1(\beta_0) = 6$ resp.\ 1 if
$\beta_{0,1}=0$ resp.\ 1. Thus (4.4.2) holds, and statement (B) is
verified. That completes the argument for Case~2.

{\it Case 3:} $[\beta^t_{-5}\mid \beta^t_{-4}\mid \cdots\mid
\beta^t_0] \not= I_6$ and $L\le -2$.  By Definition  4.3(ii),
$$g_{\basic}(\beta^{-1}_{-\infty}) = h_{-L-1}(\beta^{-1}_{L+1})
\leqno{(4.4.4)}$$
and
$$g_{\basic}(\beta^0_{-\infty}) = h_{-L}(\beta^0_{L+1}).
\leqno{(4.4.5)}$$

If the hypothesis of statement(A) (in Lemma~4.4) holds for some
$j\in\{1$, $2$, $\dots$, $5\}$, then (for that $j$)
$h_{-L-1}(\beta^{-1}_{L+1}) =j$ by (4.4.4), and then
$h_{-L}(\beta^0_{L+1}) =j$ resp.\ $j+1$ if $\beta_{0,j+1} =0$ resp.\ 1 by equation (4.1.2) in Definition~4.1, and
then (4.4.1) --- the conclusion of statement (A) --- holds by
(4.4.5).

If instead the hypothesis of statement (B) holds, then
$h_{-L-1}(\beta^{-1}_{L+1}) =6$ by (4.4.4), $h_{-L}(\beta^0_{L+1})
=6$ resp.\ 1 if $\beta_{0,1} =0$ resp.\ 1 by (4.1.3) in Definition
4.1, and then (4.4.2) --- the conclusion of statement (B) --- holds
by (4.4.5). That completes the argument for Case~3, and the proof of
Lemma~4.4. \medskip

{\sc Definition 4.5.} Define the function $g_{\spaced}
(\{0,1\}^6)^{\N} \to \{0,1,2,3,4$, $5,6\}$ as follows: For any given
(``left-infinite'') sequence
$(\beta_0,\beta_{-1},\beta_{-2},\dots\,)$ of elements of
$\{0,1\}^6$, define
$$\leqalignno{\qquad
&g_{\spaced}(\beta_0,\beta_{-1},\beta_{-2},\dots\,)&(4.5.1)\cr
 &\,:= \cases{
g_{\basic}(\beta_0,\beta_{-1},\beta_{-2},\dots\,) &if
$g_{\basic}(\beta_0,\beta_{-1},\beta_{-2},\dots\,)$ \cr & $\quad
\not= g_{\basic}(\beta_{-1},\beta_{-2},\beta_{-3},\dots\,)$ \cr 0 &
otherwise.\cr} \cr}$$
\medskip

(The subscript ``spaced'' in (4.5.1) is motivated by conclusion (A)
in the next lemma---think of elements $\{1,\dots,6\}$ ``spaced
apart'' with 0's in between.)
\medskip

{\sc Lemma 4.6.} {\it Refer to Sections 2.4(D) and 2.3(A). Suppose
$(\beta_k$, $k\in \Z)$ is a two-sided standard sequence of elements
of $\{0,1\}^6$. For each $k\in \Z$, define the numbers
$$u_k := g_{\basic} (\beta_k,\beta_{k-1},\beta_{k-2},\dots\,)
\leqno{(4.6.1)}$$
and (see (4.5.1))
$$w_k := g_{\spaced} (\beta_k,\beta_{k-1},\beta_{k-2},\dots\,) =
\cases{ u_k &if $u_k \not= u_{k-1}$ \cr
0 &if $u_k =u_{k-1}.$ \cr}
\leqno{(4.6.2)}$$
Then (A) the sequence $(w_k$, $k\in \Z)$ satisfies condition $\cS$;
and (B) one has that
$$\{k\in {\Z} :w_k=1\} = \{k\in {\Z} : u_{k-1} =6\ \ 
{\rm and}\ \ 
u_k =1\},
\leqno{(4.6.3)}$$
and for each $i\in \{2,3,4,5,6\}$, }
$$\{k\in{\Z} : w_k =i\} = \{k\in {\Z} :u_{k-1} =i-1\ \  \hbox{and}\ \ 
u_k =i\}.
\leqno{(4.6.4)}$$
\medskip

{\sc Proof.} By (4.6.1) and Lemma~4.4(A)(B), for each 
$k \in \Z$,
$$u_k -u_{k-1} \equiv 0 \,\, \hbox{or 1 mod 6.}
\leqno{(4.6.5)}$$
Hence by (4.6.2),
$$\{k\in {\Z} :w_k =6\} = 
\{k\in {\Z}: u_{k-1} =5 \,\,\hbox{and}\,\,
u_k =6\}. \leqno{(4.6.6)}$$ Also, for each $\ell\in \Z$ such that
$[\beta_{\ell-5}\mid \beta_{\ell-4} \mid \cdots \mid \beta_{\ell}]
=I_6$, one has that $u_\ell =6$ by (4.6.1) and Definition~4.3(i) and
$u_{\ell-1} =5$ by (4.6.1), Definition 4.3(ii), 
and Lemma~4.2. Hence (by the hypothesis
of Lemma~4.6) the set in (4.6.6) is doubly infinite (see section
2.1(F)). Hence also the (``larger'') set $\{k\in{\Z} :u_{k-1} \not=
u_k\}$ is doubly infinite.

Now suppose $k\in\Z$ is such that $w_k =6$. Then $u_k =6$ by
(4.6.6). Now by (4.6.5), either $u_{k+1}=1$ or there exists $m\ge 2$
such that $(u_{k+1},u_{k+2}$, $\dots,u_{k+m}) =(6,6,\dots,6,1)$.
Hence by (4.6.2), either $w_{k+1}=1$ or there exists $m\ge 2$ such
that $(w_{k+1},w_{k+2},\dots,w_{k+m}) = (0,0,\dots,0,1)$.

By a similar argument, if $i\in \{1,2,3,4,5\}$ and $k\in \Z$ is such
that $w_k =i$, then either $w_{k+1} =i+1$ or there exists $m\ge 2$
such that $(w_{k+1},\dots,w_{k+m}) = (0,0,\dots,0,i+1)$.

Since the set $\{k\in {\Z}:w_k =6\}$ in (4.6.6) is doubly infinite (as
was noted above), it now follows (by trivial induction) that for
each $i\in \{1,2,3,4,5\}$, the set $\{k\in {\Z}:w_k=i\}$ is also
doubly infinite. From the preceding two paragraphs, we now have that
the sequence $(w_k$, $k\in \Z)$ satisfies Condition $\cS$ (again see
Definition~2.3(A)). Equation (4.6.3) and (for $2\le i\le 6$)
equation (4.6.4) now follow from (4.6.1), (4.6.2), and (4.6.5).
Lemma~4.6 is proved.
\medskip

{\sc Lemma 4.7.} {\it For each $n\in\N$, there exist functions
$\rho_n:(\{0,1\}^6)^n \to \{1,2,3,4,5,6\}$ and $\varphi_n:
(\{0,1\}^6)^n \to \{0,1,2,3,4,5,6\}$ such that following holds:

For any (``left infinite'') back-standard sequence
$(\beta_n,\beta_{n-1},\beta_{n-2},\dots\,)$ of elements of
$\{0,1\}^6$ (see Section 2.4(D)) such that
$$\left[ \beta^t_{-5}\mid \beta^t_{-4} \mid \cdots \mid
\beta^t_0\right] = I_6,
\leqno{(4.7.1)}$$
one has that
$$g_{\basic}(\beta_n,\beta_{n-1},\beta_{n-2},\dots\,) =
\rho_n(\beta_1,\beta_2,\dots,\beta_n) \leqno{(4.7.2)}$$ and}
$$g_{\spaced} (\beta_n,\beta_{n-1},\beta_{n-2},\dots\,) =
\varphi_n(\beta_1,\beta_2,\dots,\beta_n). \leqno{(4.7.3)}$$
\medskip

{\sc Proof.} Again the notations of Section 2.4(E) will be used.

In the arguments below, keep in mind that the equality
$$\left[\beta^t_{\ell-5}\mid \beta^t_{\ell-4}\mid \cdots \mid
\beta^t_{\ell}\right] = I_6 \leqno{(4.7.4)}$$
{\it cannot}\/ hold
for $\ell \in \{1,2,\dots,5\}$, for that would contradict (4.7.1).

The first task is to define the functions $\rho_n$, $n\in\N$.

For each $n\in \{1,2,3,4,5\}$, define the function
$\rho_n:(\{0,1\}^6)^n \to \{1$, $2$, $\dots, 6\}$ by $\rho_n=h_n$
from Definition~4.1.  If $1 \leq n \leq 5$ and (4.7.1) holds, then (4.7.2) holds by
Definition~4.3(ii).

Now suppose instead that $n\ge 6$. Define the function
$\rho_n(\{0,1\}^6)^n \to \{1,2,\dots,6\}$ as follows: Suppose
$\beta_1,\beta_2,\dots,\beta_n$ each $\in \{0,1\}^6$. Let $S$ denote
the set of all integers $\ell\in \{6,7,\dots,n\}$ such that (4.7.4)
holds. If that set $S$ is empty, then define $\rho_n(\beta^n_1) :=
h_n(\beta^n_1)$ (from Definition~4.1). If that set $S$ is nonempty
but does not contain $n$, then define $\rho_n(\beta^n_1) :=
h_{n-L}(\beta^n_{L+1})$ where $L$ is the greatest element of $S$. If
$n\in S$, then define $\rho_n(\beta^n_1) :=6$. That completes the
definition of $\rho_n$. Now using Definition~4.3(i)(ii), one can
verify, case by case, that if (4.7.1) holds then (4.7.2) holds.

Our next task is to define the functions $\varphi_n$, $n\in\N$.

Define the function $\varphi_1: \{0,1\}^6 \to \{0,1,\dots,6\}$ as
follows: For $\alpha := (\alpha_1,\alpha_2,\dots,\alpha_6) \in
\{0,1\}^6$,
$$ \varphi_1(\alpha) :=\alpha_1 = \cases{
1 & if $\alpha_1=1$ \cr
0 &if $\alpha_1=0$ \cr}
\leqno{(4.7.5)}$$
Now for each $n\ge 2$, define the function $\varphi_n
:(\{0,1\}^6)^n\to \{0,1,\dots,6\}$ as follows: If
$\beta_1,\beta_2,\dots,\beta_n$ each $\in \{0,1\}^6$, define
$$\varphi_n(\beta^n_1) := \cases{
\rho_n(\beta^n_1) &if $\rho_n(\beta^n_1) \not=
\rho_{n-1}(\beta^{n-1}_1)$ \cr
0 &otherwise. \cr}
\leqno{(4.7.6)}$$
That completes the definition of $\varphi_n$ for $n\in\N$.

For $n\ge 2$, if (4.7.1) holds, then (4.7.3) holds by (4.7.6),
(4.7.2), and Definition~4.5. For $n=1$, note that if (4.7.1) holds,
then $g_{\basic}(\beta^0_{-\infty}) =6$ by Definition 4.3(i),
$g_{\basic}(\beta^1_{-\infty} )=6$ resp.\ 1 if $\beta_{1,1} =0$
resp.\ 1 by Lemma 4.4(B), $g_{\spaced}(\beta^1_{-\infty}) =0$ resp.\
1 if $\beta_{1,1}=0$ resp.\ 1 by Definition 4.5, and hence (4.7.3)
holds by (4.7.5). That completes the proof of Lemma~4.7.
\medskip

{\sc Lemma 4.8.} {\it Suppose $(\beta_k$, $k\in\Z)$ is a two-sided
standard sequence of elements of $\{0,1\}^6$ (see Section 2.4(D)).
Suppose $m$ is a positive integer, and for each $u\in
\{1,2,\dots,m\}$, $z_u := (0,0,0,0,0,0)$. Then
$$\forall\,\, u\in \{1,\dots,m\},\quad
g_{\spaced}(z_u,z_{u-1},\dots,z_1,\beta_0,\beta_{-1},\beta_{-2},\dots\,)
=0; \leqno{(4.8.1)}$$
and}
$$\leqalignno{&&(4.8.2) \cr
\forall\,\, n\in\N, \quad\quad
&g_{\spaced}(\beta_n,\beta_{n-1},\dots,\beta_1,z_m,z_{m-1},\dots,z_1,\beta_0,\beta_{-1},
\beta_{-2},\dots\,) \cr &\quad =
g_{\spaced}(\beta_n,\beta_{n-1},\dots,\beta_1,\beta_0,\beta_{-1},\beta_{-2},\dots\,).
\cr}$$
\medskip

{\sc Proof.} It suffices to prove this lemma for the case $m=1$, for
then the lemma as stated follows by induction on~$m$.

Accordingly, let $z := z_1 = (0,0,0,0,0,0)$. Then
$$g_{\basic}(z,\beta_0,\beta_{-1},\beta_{-2},\dots\,)
=g_{\basic}(\beta_0,\beta_{-1},\beta_{-2},\dots\,) \leqno{(4.8.3)}$$
by Lemma 4.4 ((A) or (B), whichever applies). Hence
$g_{\spaced}(z,\beta_0,\beta_{-1}$, $\beta_{-2},\dots\,)=0$ by
Definition 4.5, giving (4.8.1) (in the case $u=m=1$).
Our remaining task is to prove (4.8.2) (for $m=1$).

If $g_{\basic}(\beta_0,\beta_{-1},\beta_{-2},\dots\,) =j \in
\{1,2,\dots,5\}$, then by (4.8.3) and Lemma 4.4(A),
$$g_{\basic} (\beta_1,z,\beta_0,\beta_{-1},\beta_{-2},\dots\,) =
g_{\basic} (\beta_1,\beta_0,\beta_{-1},\beta_{-2},\dots\,),
\leqno(4.8.4)$$
with the common value being $j$ resp.\ $j+1$ if
$\beta_{1,j+1}=0$ resp.\ 1.  If instead
$g_{\basic}(\beta_0,\beta_{-1},\beta_{-2},\dots\,) =6$, then
(4.8.4) holds similarly by (4.8.3) and Lemma~4.4(B).

Starting with (4.8.4) and applying induction on $n$, using Lemma
4.4(A)(B) one obtains that
$$\leqalignno{\quad\forall\,\, n\in \N,\quad&
g_{\basic}(\beta_n,\beta_{n-1},\dots,\beta_1,z,\beta_0,\beta_{-1},\beta_{-2},\dots\,)
&(4.8.5)\cr &\qquad =
g_{\basic}(\beta_n,\beta_{n-1},\dots,\beta_1,\beta_0,
\beta_{-1},\beta_{-2},\dots\,). \cr}$$ By (4.8.3), (4.8.5), and
Definition~4.5, one has that for all $n\in\N$, (4.8.5) holds with
the subscript ``basic'' replaced on both sides by ``spaced.'' That
is, (4.8.2) holds for $m=1$. That completes the proof.
\bigskip

\centerline {\bf 5. A Markov chain and a related process} 
\bigskip
This section will
build on section~4, and will give a study of two particular strictly
stationary sequences --- one a Markov chain, and the other closely
related to it --- that together will play a key role as a ``
building block'' in the construction of the random sequence $X$ for
Theorem~1.1.

Throughout this section, the setting is a probability space
$(\Omega,\cF,P)$, rich enough to accommodate all random variables
defined in Construction 5.1 below.
\medskip

{\sc Construction 5.1.} (A) On the probability space
$(\Omega,\cF,P)$, let $(\xi_{k,i}$, $k\in \Z$, $i\in
\{1,2,3,4,5,6\})$ be an array of independent, identically
distributed $\{0,1\}$--valued random variables such that for each
$(k,i)$,
$$P(\xi_{k,i} =0) = 5/8 \quad \hbox{and}\quad (\xi_{k,i}=1) =3/8.
\leqno{(5.1.1)}$$
For each $k\in \Z$, define the random vector
$$\xi_k :=
(\xi_{k,1},\xi_{k,2},\xi_{k,3},\xi_{k,4},\xi_{k,5},\xi_{k,6})
\leqno{(5.1.2)}$$
(with the above $\xi_{k,i}$'s). Then $\xi:= (\xi_k$, $k\in \Z)$ is a
sequence of independent, identically distributed $\{0,1\}^6$--valued
random vectors.
\medskip

(For technical convenience, we assume that the random variables
$\xi_{k,i}$ are all defined in $\{0,1\}$ at {\it every}
$\omega\in\Omega$, not just on a set of probability~1.)
\medskip

(B) Refer to Definitions 4.3 and 4.5.  Define the sequence $U:=
(U_k$, $k\in \Z)$ of $\{1,2,3,4,5,6\}$--valued random variables, and
the sequence of $W:= (W_k$, $k\in\Z)$ of $\{0,1,2,3,4,5,6\}$--valued
random variables as follows:  For each $k\in\Z$,
$$U_k:= g_{\basic}
\left(\xi_k,\xi_{k-1},\xi_{k-2},\dots\,\right) \quad \hbox{and}
\leqno{(5.1.3)}$$
$$\leqalignno{ W_k&:=
g_{\spaced}\left(\xi_k,\xi_{k-1},\xi_{k-2},\dots\,\right) &(5.1.4) \cr
&\ = U_k \cdot I(U_k \not= U_{k-1})  \cr}
$$
where the second equality in (5.1.4) comes from (5.1.3) and
Definition~4.5.
\medskip

(C) For convenient reference, here is a list (with some redundancy)
of some basic properties of these random variables.  First,
$$\forall\,\,k\in{\Z},\,\,\forall\,\,\omega\in\Omega,\quad U_k(\omega)
\in \{1,2,3,4,5,6\} \,\,\hbox{and}\,\,W_k(\omega) \in
\{0,1,2,3,4,5,6\}
\leqno{(5.1.5)}$$
by (5.1.3), (5.1.4), and Definitions 4.3 and 4.5. Next,
$$\forall\,\, k\in {\Z},\,\, \forall\,\, \omega\in \Omega,\,\,
U_k(\omega) -  U_{k-1}(\omega) \equiv 0 \,\,\hbox{or 1 mod 6.}
\leqno{(5.1.6)}$$
For any $\omega\in\Omega$ such that the (``left-infinite'') sequence
$(\xi_j(\omega)$, $\xi_{j-1}(\omega)$, $\xi_{j-2}(\omega),\dots\,)$ of
elements of
$\{0,1\}^6$ is back-standard (see Section 2.4(D)) for some (hence
every) $j\in \Z$, (5.1.6) holds by (5.1.3), (5.1.5), and Lemma 4.4.
For all other $\omega\in\Omega$, (5.1.6) holds trivially (with
$U_k(\omega) =6$ for all $k\in \Z$) by (5.1.3) and the third
paragraph of Definition 4.3.  Also, referring to (5.1.5), one has
that for each $k\in\Z$, by (5.1.3), (5.1.4), (5.1.6), and
Definition~4.5,
$$\leqalignno{ \qquad&\{W_k=0\} = \{U_k =U_{k-1}\}; &(5.1.7)\cr
&\{W_k =1\} = \{ U_k =1\} \cap \{U_{k-1} =6\}; \,\hbox{and}\cr
&\forall\,\,i \in \{2,3,4,5,6\},\,\, \{W_k =i\} = \{U_k =i\} \cap
\{U_{k-1} =i-1\}.  \cr}$$
\medskip

{\sc Lemma 5.2.} {\it In the context of Construction 5.1, the
following statements hold (see Section 2.4(A)(D)):

(A) One has that $P([ \xi^t_{-5}\mid \xi^t_{-4}\mid\cdots \mid
\xi^t_0] =I_6) = (5/8)^{30} \cdot (3/8)^6$.

(B) The sequence $\xi$ is two-sided standard a.s.

(C) Defining the random variable
$$\tau:= \min \left\{ n\ge 6: \left[\xi^t_{n-5}\mid \xi^t_{n-4} \mid
\cdots \mid \xi^t_n\right] = I_6\right\} \leqno{(5.2.1)}$$ (note
that that set is a.s.\ nonempty, in fact infinite, by (B)), one has
that $E\tau\le 6 \cdot (8/5)^{30} \cdot (8/3)^6$.

(D) Suppose $\cA$ is a $\sigma$-field ($\subset \cF$, in the
underlying probability space $(\Omega,\cF,P))$ such that $\cA$ is
independent of the sequence $\xi$. Suppose $\kappa :=
(\dots,\kappa(-1),\kappa(0),\kappa(1),\dots\,)$ is a random,
$\cA$--measurable, strictly increasing sequence of integers. Then
(i)~$\widetilde\xi := (\widetilde\xi_j, j \in {\Z})  
:= (\xi_{\kappa(j)}$, $j\in \Z)$ is a sequence of
independent, identically distributed $\{0,1\}^6$--valued random
variables with the same marginal distribution as that of the random
variables $\xi_k$, $k\in\Z$; (ii)~this sequence $\widetilde \xi$ is
independent of the $\sigma$-field $\cA$; (iii)~this sequence
$\widetilde \xi$ is two-sided standard a.s.; and (iv)~defining the
random variable $\widetilde \tau$ to be the analog of the right hand
side of (5.2.1), with $\xi^t_{n-\ell}$ replaced by
$\widetilde\xi^t_{n-\ell}$ for each $\ell\in\{0,1,\dots,5\}$, one
has that $E\widetilde \tau \le 6 \cdot (8/5)^{30}\cdot (8/3)^6$. }
\medskip

{\sc Proof.} Statement (A) holds trivially by (5.1.1) and its entire
sentence (since the $6\times 6$ identity matrix $I_6$ has 30 0's and
6 1's).

Next, the $\{0,1\}$--valued random variables
$$Z_k:= I\left(\left[ \xi^t_{k-5}\mid \xi^t_{k-4}\mid\cdots \mid
\xi^t_k\right] =I_6\right),
\leqno{(5.2.2)}$$
$k\in \{\dots,-12,-6,0,6,12,\dots\,\}$ are independent and
identically distributed with $P(Z_k=1)=(5/8)^{30}\cdot (3/8)^6>0$.
Hence by two applications of the strong law of large numbers,
statement (B) holds.

Next, defining the random variable $T:= \min \{m\ge 1: Z_{6m } =1\}$
(see (5.2.2)), one has that $\tau \le 6T$ a.s.\ by (5.2.1), and hence
$E\tau \le 6\cdot ET$. Of course $T$ is finite a.s.\ and has a
geometric distribution: For each $n\in\N$, $P(T =n) = (1-p)^{n-1}p$
where $p:= (5/8)^{30}\cdot (3/8)^6$. Hence as a standard fact (see
e.g.\ [27, p.~212, Example 3, Problem 1]), $ET = 1/p =
(8/5)^{30}\cdot (8/3)^6$. Statement (C) follows.
\medskip

{\sc Proof of statement (D).} Suppose first that $m$ is a positive integer.

Next, suppose $A\in \cA$ is such that $P(A)>0$.

Next, suppose $j(-m),j(-m+1),\dots,j(m)$ are $2m+1$ integers such
that the event $F:= \bigcap^m_{u=-m}\{\kappa(u) = j(u)\}$ satisfies
$P(F\cap A)>0$. Then (see Section 2.5(A)), since $F\cap A\in \cA$,
$$\eqalign{ &\cL\left(\widetilde \xi_{-m},\widetilde\xi_{-m+1},\dots,
\widetilde\xi_{m}\mid A\cap F\right) \cr
&=
\cL\left(\xi_{\kappa(-m)},\xi_{\kappa(-m+1)},\dots,\xi_{\kappa(m)}\mid
A\cap F\right) \cr
&=\cL\left( \xi_{j(-m)},\xi_{j(-m+1)},\dots,\xi_{j(m)}\mid A\cap
F\right) \cr
&=\cL\left(\xi_{j(-m)},\xi_{j(-m+1)},\dots,\xi_{j(m)}\right) \cr
&=\lambda\times\dots\times \lambda, \cr}$$
the $(2m+1)$-order product measure (here and below), where $\lambda$
is the marginal distribution (on $\{0,1\}^6$) of the $\xi_\kappa$'s.

It follows from a simple application of Remark 2.8 that
$\cL(\widetilde\xi_{-m}$, $\widetilde\xi_{-m+1}$,
$\dots,\widetilde\xi_m\mid A) = \lambda\times \dots\times\lambda$.
By the same argument with $A$ replaced by $\Omega$,
$$\cL(\widetilde\xi_{-m},\widetilde\xi_{-m+1},\dots,\widetilde\xi_m)
= \lambda\times\dots\times \lambda.$$
Since $A\in\cA$ (with $P(A)>0$) was arbitrary, it now also follows
that the random vector
$(\widetilde\xi_{-m},\widetilde\xi_{-m+1},\dots,\widetilde\xi_m)$ is
independent of the $\sigma$--field~$\cA$.

Since $m\in\N$ was arbitrary, conclusions (i) and (ii) in statement
(D) now follow from standard measure--theoretic arguments. Finally,
conclusions (iii) and (iv) in statement (D) now follow from an
application of statements (B) and (C) to the sequence
$\widetilde\xi$. This completes the proof of statement (D), and of
Lemma~5.2.
\medskip

{\sc Lemma 5.3.} {\it The random sequence $U$ in Construction 5.1(B)
(equation (5.1.3)) is a strictly stationary, irreducible, aperiodic
Markov chain with state space $\{1,2,3,4,5,6\}$, with invariant
marginal distribution given by
$$P(U_0=j) =1/6 \quad \forall\,\, j\in \{1,2,3,4,5,6\},
\leqno{(5.3.1)}$$ and with one-step transition probability matrix
$(p_{ij}$, $1\le i,j \le 6)$ (where $p_{ij} = P(U_1 =j\mid U_0 =i))$
given by
$$\leqalignno{ p_{ii} &:= 5/8 \quad \hbox{for $i\in
\{1,2,3,4,5,6\}$,} &(5.3.2) \cr p_{i,i+1} &:= 3/8 \quad \hbox{for
$i\in \{1,2,3,4,5\}$,} \cr p_{61} &:= 3/8, \quad \hbox{and} \cr
p_{ij} &:= 0 \quad \hbox{for all other ordered pairs $(i,j)$.} \cr}$$
} 
\medskip

{\sc Proof.} Since the sequence $\xi := (\xi_k$, $k\in \Z)$ of
$\{0,1\}^6$--valued random variables is i.i.d., the strict
stationarity of the sequence $U$ is a standard consequence of
(5.1.3). Also, the random variables $U_k$, $k\in \Z$ take their
values in the set $\{1,2,\dots,6\}$ by e.g.\ (5.1.5).

Our next task is to show the following: If $i$ and $j$ are each an
element of $\{1,2,\dots,6\}$, $K$ is an integer, $A\in
\sigma(U_K,U_{K-1},U_{K-2},\dots\,)$ is an event, and $P(A\cap \{U_K
=i\})>0$, then
$$\leqalignno{ &P(U_{K+1} =j\mid A\cap \{U_K=i\}) &(5.3.3) \cr
&\qquad =P\left(U_{K+1} =j \mid U_K =i\right) = p_{ij} \cr}$$ where
$p_{ij}$ is as defined in (5.3.2). Once that is established, it will
follow that the sequence $U$ is a Markov chain with one-step
transition probability matrix (5.3.2). From that matrix, it will
then be easy to see that the Markov chain $U$ is irreducible (by the
second and third lines in (5.3.2)) and aperiodic (by the first line
in (5.3.2)). Also, it is easy to show that the (unique) invariant
distribution on the state space $\{1,2,\dots,6\}$ for the one-step
transition probability matrix (5.3.2) is the uniform distribution on
$\{1,2,\dots,6\}$ (as in (5.3.1)). Thus, once (5.3.3) is verified,
the proof of Lemma~5.3 will be complete.

We shall verify (5.3.3) for the case where $K=0$, $i\in
\{1,2,3,4,5\}$, and $j=i+1$. By similar arguments and strict
stationarity, one can verify (5.3.3) in the other cases.

Suppose $A\in \sigma(U_0,U_{-1},U_{-2},\dots\,)$ and $P(A\cap \{U_0
=i\})>0$. For any $\omega\in\Omega$ such that $U_0(\omega) =i$, one
has by (5.1.3) and Lemma 4.4(A) that
$$U_1(\omega) = \cases{
i &if $\xi_{1,i+1}(\omega) =0$ \cr
i+1 &if $\xi_{1,i+1}(\omega) =1$. \cr}
\leqno(5.3.4)$$
Also, the event $A\cap \{U_0=i\}$ is a member of the $\sigma$-field
$\sigma(\xi_0,\xi_{-1},\xi_{-2},\dots\,)$ (see (5.1.3) again) and is
therefore independent of the $(\{0,1\}^6$--valued) random variable
$\xi_1$. Hence by (5.3.4) and (5.1.1),
$$\leqalignno{\qquad P(U_1 =i+1\mid A\cap \{U_0 =i\}) &=P(\xi_{1,i+1} =1\mid
A\cap \{U_0=i\}) &(5.3.5) \cr
&= P(\xi_{1,i+1}=1) = 3/8. \cr}$$
Applying the same argument with $A$ replaced by $\Omega$, one has
that $P(U_1 =i+1\mid U_0 =i) =3/8$. Since $p_{i,i+1} := 3/8$ by
(5.3.2), one now has by (5.3.5) that (5.3.3) holds (for the case
$K=0$, $i\in \{1,2,\dots,5\}$, and $j=i+1$). That completes the proof
of Lemma~5.3.
\medskip

{\sc Lemma 5.4.} {\it The random sequence $W$ in Construction 5.1(B)
(equation (5.1.4)) is strictly stationary and has the following
properties:

(A) The sequence $W$ satisfies Condition $\cS$ a.s.\ (see
Definition~2.3(A)).

(B) The (marginal) distribution of $W_0$ (on $\{0,1,2,3,4,5,6\}$) is
given by $P(W_0 =0) =5/8$ and $P(W_0 =i) =1/16$ for all $i\in
\{1,2,3,4,5,6\}$.

(C) For each $j\in \Z$, the $\{0,1\}$--valued random variable $I(W_j
\not= 0)$ is independent of the $\sigma$--field $\sigma(U_k,W_k$,
$k\le j-1)$.

(D) The $\{0,1\}$--valued random variables $I(W_k \not= 0)$, $k\in
\Z$ are independent and identically distributed, each taking the
value 0 resp.\ 1 with probability $5/8$ resp.\ $3/8$. In fact for each
$j\in\Z$, the random sequence $(I(W_k \not=0)$, $k\in
\{j,j+1,j+2,\dots\,\})$ is independent of the $\sigma$--field
$\sigma(U_k,W_k$, $k\le j-1)$.}
\medskip

The redundancies here (statement (C) is a special case of (D), and
$\sigma(U_k,W_k$,  $k\le j-1) =\sigma(U_k$, $k\le j-1)$ by (5.1.4))
are for convenient reference.
\medskip

{\sc Proof.}  Since the sequence $\xi:= (\xi_k$, $k\in\Z)$ of
$\{0,1\}^6$--valued random variables is i.i.d., the strict
stationarity of the sequence $W$ is a standard consequence of
(5.1.4). Also, statement (D) (in Lemma 5.4) follows from statements
(B) and (C), strict stationarity, and an elementary induction
argument.  Our remaining task is to prove statements (A), (B), and
(C).

Statement (A) holds by Lemma 5.2(B), equation (5.1.4), and Lemma
4.6(A).

To prove statement (B), first note that (see (5.1.5)) $P(W_0=0)=1
-\sum^6_{i=1}P(W_0=i)$. Hence to prove (B) it suffices to prove that
$P(W_0=i) =1/16$ for each $i\in \{1,2,\dots,6\}$. But that holds by
(5.1.7) and Lemma 5.3; for example, for $i=6$ one thereby has
$$\eqalign{P(W_0=6) = P(U_0 =6,\, U_{-1} =5) &= P(U_{-1} =5) \cdot
P(U_0 =6\mid U_{-1} =5) \cr
&=(1/6) \cdot (3/8) = 1/16. \cr}$$
Thus statement (B) holds.
\medskip

{\sc Proof of statement (C).} The argument is the same for any $j\in
\Z$. We shall give it for $j=1$. By (5.1.7), $\sigma(U_k,W_k$, $k\le
0) = \sigma(U_k$, $k\le 0)$. Our task is to prove that the
$\{0,1\}$--valued random variable $I(W_1\not= 0)$ is independent of
the $\sigma$-field $\sigma(U_k$, $k\le 0)$.

Let $m\in\N$ be arbitrary but fixed.  By a standard
measure--theoretic argument, it suffices to show that the random
variable $I(W_1\not= 0)$ is independent of the random vector
$(U_0,U_{-1},U_{-2},\dots,U_{-m})$.

Referring to (5.1.5), suppose $i_0,i_1,\dots,i_m$ are each an
element of $\{1$, $2$, $\dots,6\}$, and that the event $F:=
\bigcap^m_{k=0}\{U_{-k} =i_k\}$ satisfies $P(F)>0$. It suffices to
show that
$$P(I(W_1\not= 0) = 0\mid F) = P(I(W_1\not= 0) =0). \leqno{(5.4.1)}$$
(For then by taking complements, one also obtains $P(I(W_1 \not=0)
=1\mid F) = P(I(W_1\not= 0)=1)$.) The right side of (5.4.1) is of
course simply $P(W_1=0)$, which is 5/8 by statement (B). By (5.1.7)
and Lemma 5.3, the left side of (5.4.1) is simply
$$\eqalign{
P(W_1 =0\mid F) &= P(U_1 =U_0\mid F) =P(U_1 = i_0\mid F) \cr
&= P(U_1 =i_0\mid U_0 =i_0) =5/8. \cr}$$
Thus (5.4.1) holds. That completes the proof of statement (C), and of
Lemma~5.4.
\medskip

{\sc Lemma 5.5.} {\it Refer to Construction 5.1(B), Lemma 5.4(A)(B),
and Definition 2.3(A). Define the random variable
$$T:= \min \{k\in{\N}: W_k =1\}. \leqno{(5.5.1)}$$
Then $E(T\mid W_0 =1)=16$ and $\var(T\mid W_0 =1) = 80/3$. }
\medskip

{\it Proof.} For convenience, for each $k\ge 1$, define the
$\{0,1\}$--valued random variable $V_k:= I(W_k \not=0)$. By Lemma
5.4(D) (both sentences of it), conditional on the event $\{W_0=1\}$,
the random variables $V_1,V_2,V_3,\dots$ are independent and
identically distributed, each taking the value 0 resp.\ 1 with
probability 5/8 resp.\ 3/8. By Lemma 5.4(A) (see Definition 2.3(A))
and (5.5.1), for $P$--a.e.\ $\omega\in\{W_0 =1\}$,
$$T(\omega) = \min \left\{ n\in{\N}: \sum^n_{k=1} V_k (\omega)
=6\right\}.$$
Hence conditional on the event $\{W_0=1\}$, the random variable $T$
has the ``negative binomial'' distribution with parameters $r=6$ and
$p=3/8$. Hence by standard elementary calculations (see e.g.\ [27,
pp.\ 213--215, Example 4, Problem 2]), $E(T\mid W_0 =1) = r/p =16$
and $\var(T\mid W_0 =1) = r(1-p)/p^2 = 80/3$. Thus Lemma 5.5 holds.
\medskip

{\sc Lemma 5.6.} {\it Refer to Section 2.7(B) and the sequence $W$ in
Construction 5.1(B). For any given $J\in\Z$, the ordered triplet
$$(\sigma(W_k, k\le J),\,\, \{W_J=1\},\,\, \sigma(W_k,k\ge J))
\leqno{(5.6.1)}$$ is a restricted Markov triplet.}
\medskip

{\sc Proof.} By strict stationarity (see Lemma 5.4), it suffices to
give the argument for $J=0$.

Refer to Section 2.7(A). We shall use some standard elementary
properties of Markov triplets (see e.g.\ [5, Vol.\ 1, Appendix, Section A701]).
By Lemma 5.3,
$$(\sigma(U_k,k\le -2),\,\, \sigma(U_{-1},U_0),\,\, \sigma(U_k,k\ge
1))$$
is a Markov triplet.  Hence so is
$$(\sigma(U_k,k\le 0), \,\, \sigma(U_{-1},U_0),\,\, \sigma(U_k,k\ge
-1)).$$
Hence (trivially) by (5.1.7), so is
$$(\sigma(W_k,k\le 0), \,\, \sigma(U_{-1},U_0),\,\, \sigma(W_k,k\ge
0)).
\leqno{(5.6.2)}$$
Now by (5.1.7), $\{W_0=1\}=\{U_0 =1\} \cap \{U_{-1} =6\}$, which is
an atom of the $\sigma$--field $\sigma(U_{-1},U_0)$. Hence from the
Markov triplet (5.6.2), one has that (5.6.1) is (for $J=0$) a
restricted Markov triplet. That completes the proof.
\bigskip

\centerline {\bf 6. Scaffolding (part 1)}
\bigskip 
Throughout the rest of this paper, the
context will be a particular probability space $(\Omega,\cF,P)$, rich
enough to accommodate all random variables defined henceforth. With
or without explicit mention, many of those ``random variables'' will
be random vectors or random sequences, taking their values in spaces
such as $(\{0,1\}^6)^{\N}$ or $\{-1,1\}^{\sxtp(n)}$.

Sections 6 (here), 7, and 8 will provide ``scaffolding'' that will be
used for the definition (in Section~7) of the random sequence $X$ for
Theorem~1.1 and for the proofs of the various properties of $X$
stated in that theorem.
\medskip

{\sc Construction 6.1.} (A) On our probability space
$(\Omega,\cF,P)$, let $\xi^{(n)}_{ki}$, $n\in\N$, $k\in\Z$, $i\in
\{1,2,3,4,5,6\}$ be an array of independent, identically distributed
$\{0,1\}$--valued random variables such that for each $(n,k,i)$,
$$P\left(\xi^{(n)}_{ki} =0\right) = 5/8 \quad\hbox{and}\quad
P\left(\xi^{(n)}_{ki}=1\right) = 3/8.
\leqno{(6.1.1)}$$
From those random variables, for convenient reference, for each
$n\in\N$ and each $k\in\Z$, define the $\{0,1\}^6$--valued random
variable (random vector)
$$\xi^{(n)}_k :=
\left(\xi^{(n)}_{k1},\xi^{(n)}_{k2},\dots,\xi^{(n)}_{k6}\right).
\leqno{(6.1.2)}$$ Of course those random variables $\xi^{(n)}_k$,
$n\in\N$, $k\in\Z$ are independent and have the same distribution (a
product measure) on $\{0,1\}^6$. From those random variables, again
for convenient reference, for each $n\in\N$, define the random
sequence $\xi^{(n)} := (\xi^{(n)}_k$, $k\in\Z)$.

(B) For each $n\in\N$, the random sequence $\xi^{(n)}$ here has
(trivially) the same distribution (on $(\{0,1\}^6)^{\Z}$) as the random
sequence $\xi$ in Construction~5.1. Henceforth, {\it theorems
involving the random sequence $\xi$ in Section 5 will be applied freely to each of the
random sequences $\xi^{(n)}$, $n\in\N$ here.}

(C) Referring to (6.1.2), for each $n\in\N$ and each $k\in\Z$,
define the $(\{0,1\}^6)^n$--valued random variable
$$\overline{\xi}^{(n)}_k := \left( \xi^{(1)}_k,
\xi^{(2)}_k,\dots,\xi^{(n)}_k\right). \leqno{(6.1.3)}$$
Of course, for any particular $n\in\N$, those random variables
$\overline{\xi}^{(n)}_k$, $k\in\Z$ are independent and identically
distributed. For convenient reference, for each $n\in\N$, define the
resulting random sequence $\overline{\xi}^{(n)}:=
(\overline{\xi}^{(n)}_k$, $k\in\Z)$.

(D) Again referring to (6.1.2), for each $k\in\Z$, define the
$(\{0,1\}^6)^{\N}$--valued random variable
$$\overline{\xi}^{(\infty)}_k := \left(
\xi^{(1)}_k,\xi^{(2)}_k,\xi^{(3)}_k,\dots\,\right). \leqno{(6.1.4)}$$
Of course those random variables $\overline{\xi}^{(\infty)}_k$ are
independent and identically distributed. Define the resulting random
sequence $\overline{\xi}^{(\infty)}:= (\overline{\xi}^{(\infty)}_k$,
$k\in\Z)$.
\medskip

{\sc Construction 6.2.} For each $n\in\N$, we shall define a sequence
$W^{(n)}:= (W^{(n)}_k$, $k\in\Z)$ of $\{0,1,2,3,4,5,6\}$--valued
random variables such that (see Section 2.6 and Construction 6.1(C))
$$\hbox{the ordered pair $(W^{(n)}$, $\overline{\xi}^{(n)})$
satisfies Condition $\cM$.} \leqno{(6.2.1)}$$
Also, for each $n\in\N$, we shall define an event
$\Omega^{(n)}_{\good}$ such that
$$P\left(\Omega^{(n)}_{\good}\right) =1 \leqno{(6.2.2)}$$
and (see Section 2.3)
$$\leqalignno{&\hbox{for each $\omega\in\Omega^{(n)}_{\good}$, the sequence
$W^{(n)}(\omega)$ } &(6.2.3)\cr &\hbox{of elements of
$\{0,1,\dots,6\}$ satisfies condition $\cS$.} \cr}$$ 
Also, for each $n \in \N$,  we shall
define nonnegative integer--valued random variables $\Psi(n,k,j)$,
$k\in\Z$, $j\in\{0,1,2,\dots\,\}$, such that for each $j\in
\{0,1,2,\dots\,\}$, defining the random sequence $\Psi(n,j) :=
(\Psi(n,k,j)$, $k\in\Z)$, one has that (see Section 2.6 again)
$$\hbox{the ordered pair $(\Psi(n,j)$, $W^{(n)})$ satisfies Condition
$\cM$.}  \leqno{(6.2.4)}$$
The definition will be recursive --- starting with just
$W^{(1)}$ and $\Omega^{(1)}_{\rm good}$, and then (for the
recursion step)  
set up with $\Psi(n,j)$ defined
together with  $W^{(n+1)}$ and $\Omega^{(n+1)}_{\good}$.
\medskip

{\it Initial step.}\ \ First, define the sequence $W^{(1)}:= (W^{(1)}_k$, $k\in\Z)$ of
$\{0,1,\dots,6\}$--valued random variables as follows: For each
$k\in\Z$ (see (6.1.2) and Definition 4.5),
$$W^{(1)}_k := g_{\spaced}
\left(\xi^{(1)}_k,\xi^{(1)}_{k-1},\xi^{(1)}_{k-2},\dots\,\right).
\leqno{(6.2.5)}$$ 
Obviously (see Sections 6.1(B) and 2.6) equation
(6.2.1) holds for $n=1$.

Also, let $\Omega^{(1)}_{\good}$ denote the set of all
$\omega\in\Omega$ such that the sequence $\xi^{(1)}(\omega) :=
(\xi^{(1)}_k(\omega)$, $k\in\Z)$ of elements of $\{0,1\}^6$ (see the very end of Construction 6.1(A)) is two-sided standard (see Notations 2.4(D)). Equation
(6.2.2) holds for $n=1$ by Lemma 5.2(B), and equation (6.2.3) holds
for $n=1$ by Lemma~4.6(A).
\medskip

{\it Recursion step.} Now suppose $n\in\N$, and the event
$\Omega^{(n)}_{\good}$ and the sequence $W^{(n)}:= (W^{(n)}_k$,
$k\in\Z)$ of $\{0,1,\dots,6\}$--valued random variables are already
defined and satisfy (6.2.1), (6.2.2), and (6.2.3).

(A) For each $k\in\Z$ and each $j\in\{0,1,2,\dots\,\}$, define the
nonnegative integer--valued random variable $\Psi(n,k,j)$ by (see
Definition 2.2)
$$\Psi(n,k,j) :=
\psi_j\left(W^{(n)}_k,W^{(n)}_{k-1},W^{(n)}_{k-2},\dots\,\right).
\leqno{(6.2.6)}$$ Now by (6.2.1) and Definition 2.6, the sequence
$W^{(n)}$ is strictly stationary. Hence by (6.2.6), 
equation (6.2.4) holds for every $j\ge 0$.

(B) Next, define the sequence $W^{(n+1)}:= (W^{(n+1)}_k$, $k\in\Z)$
of $\{0,1, \allowbreak \dots,6\}$--valued random variables as follows: For each
$k\in\Z$ and each $\omega\in\Omega$, referring to Definition 4.5 and
(6.2.6),
$$W^{(n+1)}_k(\omega):= \cases{
0 &if $W^{(n)}_k(\omega)\not= 1$ \cr
g_{\spaced}\biggl(\xi^{(n+1)}_k(\omega),\xi^{(n+1)}_{k-\Psi(n,k,1)(\omega)}(\omega),
& \cr
\,\,\xi^{(n+1)}_{k-\Psi(n,k,2)(\omega)}(\omega),\xi^{(n+1)}_{k-\Psi(n,k,3)(\omega)}(\omega),\dots\,\biggr)
&if $W^{(n)}_k(\omega) =1$. \cr} \leqno{(6.2.7)}$$

{\it Verification of (6.2.1) with $n$ replaced by $n+1$.} For each
$j\ge 0$, each $k\in\Z$, and each $\omega\in\Omega$, by (6.2.6),
$$\xi^{(n+1)}_{k-\Psi(n,k,j)(\omega)}(\omega) = \sum^\infty_{u=0}
\xi^{(n+1)}_{k-u}(\omega)\cdot I\left(\Psi(n,k,j) =u\right)(\omega).
\leqno{(6.2.8)}$$
Now by (6.2.4), (6.2.1), and Definition 2.6 (and a trivial argument),
for each $j\ge 0$, the ordered pair $(\Psi(n,j),
\overline{\xi}^{(n)})$ satisfies Condition~$\cM$. Hence by (6.2.8) and
Construction 6.1(C), for each $j\ge 0$, the ordered pair
$$\left(\left(\xi^{(n+1)}_{k-\Psi(n,k,j)},\, k\in\Z\right),
\overline{\xi}^{(n+1)}\right)$$
satisfies Condition~$\cM$. Hence by (6.2.1) for the given $n$, (6.2.7),
and Construction 6.1(C), equation (6.2.1) holds with $n$ replaced by
$n+1$.

(C) Our final task in the recursion step is to define the event
$\Omega^{(n+1)}_{\good}$ and then to verify (6.2.2) and (6.2.3) with
$n$ replaced by $n+1$. To facilitate this, we shall define a couple
of other random sequences.

Define the random strictly increasing sequence $\kappa^{(n)}:=
(\kappa(n,j)$, $j\in\Z)$ of integers as follows: Referring to (6.2.2)
and (6.2.3) (for the given $n$) and to Section 2.3, for each
$\omega\in\Omega^{(n)}_{\good}$, define the integers
$\kappa(n,j)(\omega)$, $j\in\Z$ (uniquely) by
$$\leqalignno{\dots &<\kappa(n,-2)(\omega)<\kappa(n,-1)(\omega)<\kappa(n,0)(\omega)
\le 0 \cr
 &<1 \le
\kappa(n,1)(\omega)<\kappa(n,2)(\omega)<\kappa(n,3)(\omega)<\dots
&(6.2.9) \cr}$$
and
$$\left\{k\in{\Z}: W^{(n)}_k(\omega) =1\right\} =
\{\dots,\kappa(n,-1)(\omega),\kappa(n,0)(\omega),\kappa(n,1)(\omega),\dots\,\}.
\leqno{(6.2.10)}$$
(The sequence $\kappa^{(n)}$ can be left undefined on the null set
$(\Omega^{(n)}_{\good})^c$.) By (6.2.1) and (6.2.9)--(6.2.10),
$\sigma(\kappa^{(n)})$ $\dot\subset$ $\sigma(\overline{\xi}^{(n)})$, and
hence by Construction 6.1(A)(C), the sequence $\kappa^{(n)}$ is
independent of the sequence $\xi^{(n+1)}$. Hence by Lemma 5.2(D) (and Section 6.1(B)),  
$$\leqalignno{\qquad \qquad&\hbox{the sequence $(\xi^{(n+1)}_{\kappa(n,j)}$, $j\in\Z)$ is i.i.d.\
with the}&(6.2.11)  \cr
&\hbox{ same marginal distribution as the $\xi_k$'s in Construction
5.1;} \cr}$$
$$\hbox{the sequence $(\xi^{(n+1)}_{\kappa(n,j)}$, $j\in\Z)$
 is independent of the $\sigma$-field $\sigma(\overline{\xi}^{(n)})$;
 and} \leqno{(6.2.12)}$$
$$\hbox{the sequence $(\xi^{(n+1)}_{\kappa(n,j)}$, $j\in\Z)$ is
two-sided standard a.s.\ (see Section 2.4(D)).} \leqno{(6.2.13)}$$
\medskip
\noindent Let $\Omega^{(n+1)}_{\good}$ denote the set of all
$\omega\in\Omega^{(n)}_{\good}$ such that the sequence
\break
$(\xi^{(n+1)}_{\kappa(n,j)(\omega)}(\omega)$, $j\in\Z)$ of elements
of $\{0,1\}^6$ is two-sided standard. Then by (6.2.13), equation
(6.2.2) holds with $n$ replaced by $n+1$.

{\it Verification of (6.2.3) with $n$ replaced by $n+1$.} Refer to (6.2.6), 
(6.2.9)--(6.2.10), and Definition 2.2. For each
$\omega\in\Omega^{(n+1)}_{\good}$, each $\ell\in\Z$, and each $j\ge
0$,
$$\Psi(n,\kappa(n,\ell)(\omega),j)(\omega) = \kappa(n,\ell)(\omega)
- \kappa(n,\ell-j)(\omega)$$
and hence
$$\kappa(n,\ell)(\omega) - \Psi(n,\kappa(n,\ell)(\omega),j)(\omega) =
\kappa(n,\ell-j)(\omega).$$
Hence for each $\omega\in\Omega^{(n+1)}_{\good}$ and each
$\ell\in\Z$, by (6.2.7) and (6.2.10),
$$W^{(n+1)}_{\kappa(n,\ell)(\omega)}(\omega) =
g_{\spaced}\left(\xi^{(n+1)}_{\kappa(n,\ell)(\omega)}(\omega),\xi^{(n+1)}_{\kappa(n,\ell-1)(\omega)}
(\omega),\xi^{(n+1)}_{\kappa(n,\ell-2)(\omega)}(\omega),\dots\,\right).
\leqno{(6.2.14)}$$ By (6.2.14), the definition of
$\Omega^{(n+1)}_{\good}$, and Lemma 4.6(A), one has that for each
$\omega\in\Omega^{(n+1)}_{\good}$, the sequence
$(W^{(n+1)}_{\kappa(n,\ell)(\omega)}(\omega)$, $\ell\in\Z)$ of
elements of $\{0,1,\dots,6\}$ satisfies Condition~$\cS$. Also, by
(6.2.7) and (6.2.10), for each $\omega\in\Omega^{(n+1)}_{\good}$,
one has that
$$\leqalignno{\qquad 
&\qquad W^{(n+1)}_k(\omega) =0 &(6.2.15)\cr
&\quad \indent \hbox{for all $k\in {\Z} -
\{\dots,\kappa(n,-1)(\omega),\kappa(n,0)(\omega),\kappa(n,1)(\omega),\dots\,\}$.}\cr}$$
Hence by Remark 2.3(C), equation (6.2.3) holds with $n$ replaced by
$n+1$. This completes the recursion step and Construction~6.2.

\medskip

{\sc Section 6.3.} This builds on the material in Constructions 6.1 and 6.2.

(A) Referring to the events $\Omega^{(n)}_{\good}$, $n\in\N$ in
Construction 6.2, define the event
$$\Omega_0 := \bigcap^\infty_{n=1} \Omega^{(n)}_{\good}.
\leqno{(6.3.1)}$$
Then by (6.2.2),
$$P(\Omega_0) =1.
\leqno{(6.3.2)}$$
Also, by (6.2.3), (6.2.10), and the definitions of the events
$\Omega^{(n)}_{\good}$, $n\in\N$ in Construction 6.2, one has the
following:

(i) For each $\omega\in\Omega_0$ and each $n\in\N$, the sequence
$(W^{(n)}_k(\omega)$, $k\in\Z)$ of elements of $\{0,1,\dots,6\}$
satisfies Condition $\cS$.
(Again recall Definition 2.3(A).)

(ii) For each $\omega\in\Omega_0$, the sequence
$(\xi^{(1)}_k(\omega)$, $k\in\Z)$ of elements of $\{0,1\}^6$ is
two-sided standard.

(iii) For each $\omega\in\Omega_0$ and each $n\in\N$, the sequence
$(\xi^{(n+1)}_j(\omega)$, $j\in\{k\in{\Z}:W^{(n)}_k(\omega) =1\})$ (see
Section 2.2(F)(ii)) is two-sided standard.

(B) {\it Remark.} From (6.2.7), one has that for a given $n\in\N$,
$k\in\Z$, and $\omega\in\Omega$, if $W^{(n)}_k(\omega) \not= 1$ then
$W^{(n+1)}_k(\omega) =0$. Thus for a given $k\in\Z$ and
$\omega\in\Omega$, the sequence
$(W^{(1)}_k(\omega),W^{(2)}_k(\omega),W^{(3)}_k(\omega),\dots\,)$ of
elements of $\{0,1,\dots,6\}$ will have one of the following five
forms:
\item{(i)} $(0,0,0,\dots\,)$,

\item{(ii)} $(i,0,0,0,\dots\,)$ where $i\in \{2,3,4,5,6\}$,

\item{(iii)} $(1,1,\dots,1,0,0,0,\dots\,)$,

\item{(iv)} $(1,1,\dots,1,i,0,0,0,\dots\,)$ where $i\in
\{2,3,4,5,6\}$,
\item{(v)} $(1,1,1,\dots\,)$.

(In (iii) and (iv), the number of 1's can be any positive integer.)
\medskip

(C) {\it Remark.} Recall the first sentence of Remark (B) above.  By
Remark 2.3(B) and remark (i) in Section (A) above, one has the
following: For any $n\in\N$, any $\omega\in\Omega_0$ (see (6.3.1)),
and any pair of integers $J$ and $L$ such that $J\le L$, one has that
$$\leqalignno{\qquad\qquad&\card \{k\in{\Z}:\ J\le k\le L 
\ \ \hbox{and}\ \ 
W_k^{(n+1)}(\omega) =1\} &(6.3.3) \cr &\le 1 +(1/6) \cdot
\card\{k\in{\Z}:\ J\le k\le L \ \ \hbox{and}\ \ W_k^{(n+1)}(\omega)
\not= 0\}\cr &\le 1+(1/6) \cdot \card \{k\in{\Z}:\ J\le k\le
L\ \ \hbox{and}\ \ W^{(n)}_k(\omega) =1\}. \cr}$$

(D) {\it Remark.} For each $n\in\N$, each $k\in\Z$ and each
$\omega\in\Omega_0$,
$$\leqalignno{
\qquad\qquad\Psi(n,k,0)(\omega) &= \min\left\{i\ge
0:W^{(n)}_{k-i}(\omega) =1\right\} &(6.3.4) \cr &\le \min\left\{i\ge
0: W^{(n+1)}_{k-i}(\omega) =1\right\} = \Psi(n+1,k,0)(\omega) \cr}$$
by (6.2.6) and Definition 2.2, since (again see Remark (B) above) the
second set in (6.3.4) is a subset of the first; equivalently
$$k-\Psi(n,k,0)(\omega) \ge k-\Psi(n+1,k,0)(\omega). \leqno{(6.3.5)}$$

(E) Referring to (6.2.5) and (6.2.7), for each $n\in\N$ and each
$k\in\Z$, define the $\{0,1,\dots,6\}^n$--valued random variable
$$\overline{W}^{(n)}_k :=
\left(W^{(1)}_k,W^{(2)}_k,\dots,W^{(n)}_k\right). \leqno{(6.3.6)}$$
Accordingly, for each $n\in\N$, define the random sequence
$\overline{W}^{(n)}:= (\overline{W}^{(n)}_k$, $k\in\Z)$.  By (6.2.1) and Definition 2.6(A) (and Construction 6.1(C)), 
for each $n\in\N$,
$$\hbox{The ordered pair $(\overline{W}^{(n)},\overline{\xi}^{(n)})$
satisfies Condition $\cM$.} \leqno{(6.3.7)}$$

(F) {\it Remark.} Of course the random sequence $W^{(1)}$ in (6.2.5)
has the same distribution (on $\{0,1,\dots,6\}^{\Z})$ as the random
sequence $W$ in Construction 5.1(B), by (6.2.5), (5.1.4), and
Construction 6.1(B).

(G) {\it Remark.} For each $n\in\N$, the random sequence
$(W^{(n+1)}_{\kappa(n,j)}$, $j\in\Z)$ (recall (6.3.1), (6.3.2),
(6.2.9), (6.2.10), and (6.2.14)) has the same distribution (on
$\{0,1,\dots,6\}^{\Z})$ as the random sequence $W:= (W_j$, $j\in\Z)$ in
Construction 5.1(B). (This holds by (6.2.11), (6.2.14), (5.1.4), and
a standard measure--theoretic argument.)

(H) {\it Remark.} For each $n\in\N$, the random sequence
$(W^{(n+1)}_{\kappa(n,j)}$, $j\in\Z)$ (recall (G) above) is
independent of the sequence $\overline{\xi}^{(n)}$, by (6.2.14) and
(6.2.12).
\medskip

{\sc Lemma 6.4.} {\it (A) For each $i\in \{1,2,\dots,6\}$, $P(W^{(1)}_0
=i) =1/16$.

(B) If $n\in\N$ and $A\in \sigma(\overline{W}^{(n)})$ (see section
6.3(E)) are such that $P(A\cap \{W^{(n)}_0 =1\})>0$, then for each
$i\in \{1,2,\dots,6\}$, $P(W^{(n+1)}_0 = i\mid A\cap \{W^{(n)}_0
=1\}) = 1/16$.

(C) For each $n\in\N$ and each $i \in \{1,2,\dots,6\}$, $P(W^{(n)}_0
=i) = 16^{-n}$.

(D) Suppose $n\in\N$. Suppose $S$ is a nonempty finite set $\subset
\Z$. Suppose $A\in \sigma(\overline{W}^{(n)})$. Suppose also that the
event
$$F:= A\bigcap \left[ \bigcap_{k\in S} \left\{W^{(n)}_k =1\right\}\right]
\leqno{(6.4.1)}$$
satisfies $P(F)>0$. Then conditional on $F$, the $\{0,1\}$--valued
random variables $I(W^{(n+1)}_k \not=0)$, $k\in S$ are independent and
identically distributed, each taking the value 0 resp.\ 1 with
probability 5/8 resp.\ 3/8.}
\medskip

{\sc Proof.}  Statement (A) holds by Lemma~5.4(B) and 
Remark 6.3(F).

To prove statement (B), note first that by (6.2.9)--(6.2.10),
$\kappa(n,0)(\omega) =0$ for each $\omega\in\Omega_0\cap \{W^{(n)}_0
=1\}$ (recall (6.3.1)--(6.3.2)). By (6.3.6), (6.3.7), and the
hypothesis of statement (B), $A\cap \{W^{(n)}_0 =1\} \in
\sigma(\overline{W}^{(n)})$ $\dot\subset$
$\sigma(\overline{\xi}^{(n)})$. Hence for a given $i\in
\{1,2,\dots,6\}$, by Remark 6.3(H), Remark 6.3(G), and Lemma 5.4(B),
$$\eqalign{P\left(W^{(n+1)}_0 =i\,\Bigl|\, A\cap \{W^{(n)}_0 =1\}\right) &=
P\left(W^{(n+1)}_{\kappa(n,0)} =i\,\Bigl|\, A\cap \{W^{(n)}_0 =1\}\right)
\cr
&=P\left(W^{(n+1)}_{\kappa(n,0)}=i\right) = P(W_0 =i) = 1/16 \cr}$$
(where $W_0$ is as in Construction 5.1(B)). Thus statement (B) holds.

Statement (C) holds for $n=1$ by statement (A). Also, for a given $n\in\N$ for which Statement (C) holds, 
and a given 
$i\in \{1,2,\dots,6\}$, one has that $\{W^{(n+1)}_0 =i\}
\subset \{W^{(n)}_0 =1\}$ by Remark 6.3(B), and hence by statement
(B) (in Lemma 6.4),
$$\eqalign{P(W^{(n+1)}_0 =i) &= P(W^{(n)}_0 =1) \cdot P(W^{(n+1)}_0
=i\mid W^{(n)}_0 =1) \cr
&= P(W^{(n)}_0 =1) \cdot 1/16 = (1/16)^{n+1}. \cr}$$
Now statement (C) holds for all $n \in \N$ by induction. 
\medskip

{\it Proof of (D).} Let $m:= \card$ $\cS$. Let $\lambda$ denote the
``Bernoulli'' probability measure on $\{0,1\}$ given by
$\lambda(\{0\}) = 5/8$ and $\lambda(\{1\}) =3/8$. In the argument
below, the notation $\lambda\times \dots\times \lambda$ will mean the
$m$--fold product measure on $\{0,1\}^m$.

Let the elements of $S$ be denoted by $s(i)$, $1\le i\le m$ where
$s(1)<s(2)<\dots<s(m)$. Refer to (6.2.9)--(6.2.10)
and to (6.3.1)--(6.3.2). For each
$m$--tuple $j := (j(1),j(2),\dots,j(m))$ of integers such that
$j(1)<j(2)<\dots<j(m)$, define the event
$$F(j):= F\bigcap \Omega_0 \bigcap 
\left[ \bigcap^m_{u=1}\{\kappa(n,j(u)) =
s(u)\}\right].$$
Those events $F(j)$ are (pairwise) disjoint (and some of them will be
empty), and by (6.4.1) and (6.2.9)--(6.2.10) their union is 
$F \cap \Omega_0 \doteq F$.
Also, for each such $j$, one has that $F(j) \dot\in
\sigma(\overline{W}^{(n)})$ $\dot\subset$
$\sigma(\overline{\xi}^{(n)})$ by (6.3.6), (6.3.7), and the
hypothesis of statement (D) (since by (6.2.9)--(6.2.10) the random
variables $\kappa(n,i)$, $i\in\Z$ are $\sigma(W^{(n)})$--measurable modulo the null-set $\Omega_0^c$).

Hence by Remark 6.3(H), Remark 6.3(G), and Lemma 5.4(D) (see Section 2.5(A)),
for each such $m$-tuple $j$, 
$$\eqalign{&\cL \left(I(W^{(n+1)}_{s(1)} \not=0),I(W^{(n+1)}_{s(2)}
\not=0),\dots,I(W^{(n+1)}_{s(m)} \not=0)\,\Bigl|\, F(j)\right) \cr
&=\cL\left(I(W^{(n+1)}_{\kappa(n,j(1))}
\not=0),I(W^{(n+1)}_{\kappa(n,j(2))} \not=0),
\dots,I(W^{(n+1)}_{\kappa(n,j(m))} \not=0)\,\Bigl|\, F(j)\right) \cr
&=\cL\left(I(W^{(n+1)}_{\kappa(n,j(1))}\not=
0),I(W^{(n+1)}_{\kappa(n,j(2))}
\not=0),\dots,I(W^{(n+1)}_{\kappa(n,j(m))}\not=0)\right) \cr
&=
\lambda\times\dots\times \lambda. \cr}$$ 
Hence by Remark 2.8,
$$\cL\left( I(W^{(n+1)}_{s(1)} \not=0),I(W^{(n+1)}_{s(2)}
\not=0),\dots, I(W^{(n+1)}_{s(m)}\not= 0)\,\Bigl|\, F\right) =
\lambda\times\dots\times\lambda. $$ Thus statement (D) holds. Lemma
6.4 is proved.
\medskip

{\sc Lemma 6.5.} {\it Suppose $n\in\N$. Then for each integer $J$,
the ordered triplet
$$\left(\sigma(\overline{W}^{(n)}_k, k\le J),\,\, \{W^{(n)}_J =1\},
\,\, \sigma(\overline{W}^{(n)}_k, k\ge J+1)\right)$$ is a restricted
Markov triplet (see section 2.7(B)).}
\medskip

The slight ``asymmetry'' in this statement is just a matter of
convenience. This statement is slightly stronger than what will be
needed in its application later on.
\medskip

{\sc Proof.}  For $n=1$, this holds by Lemma 5.6 and Remark~6.3(F).

Now for induction, suppose $N$ is a positive integer and Lemma 6.5
holds for $n=N$. Our task is to show that it holds for $n=N+1$.

Of course by (6.3.7) (with $n=N+1$) and Definition 2.6(A), the
sequence $\overline{W}^{(N+1)}$ is strictly stationary. Hence it
suffices to carry out the argument in the case where ($n=N+1$ and)
$J=0$.

Refer to (6.2.9)--(6.2.10) for $n=N$. For convenient notation,
define the sequence $Z^*:= (Z^*_j$, $j\in\Z)$ of
$\{0,1,\dots,6\}$--valued random variables as follows: For each
$j\in\Z$
$$Z^*_j := W^{(N+1)}_{\kappa(N,j)}. \leqno{(6.5.1)}$$

By the induction hypothesis,
$$\left(\sigma(\overline{W}^{(N)}_k,k\le 0),\,\,  \{W^{(N)}_0 =1\}, \,\, \sigma(\overline{W}^{(N)}_k,k\ge 1)\right)$$
is a restricted Markov triplet. By (6.5.1), Remark 6.3(G), and Lemma 5.6,
$$\left(\sigma(Z^*_j,j\le 0),\,\, \{Z^*_0=1\},\,\, \sigma(Z^*_j, j\ge
0)\right)$$
is a restricted Markov triplet. Hence by 
(6.3.7), Remark 6.3(H), and Remark~2.7(D),
$$\leqalignno{
\qquad\quad\bigl(\sigma(\overline{W}^{(N)}_k,k\le 0) \vee \sigma(Z^*_j,j\le 0),\,\,
 &\{W^{(N)}_0=1\} \cap \{Z^*_0 =1\}, &(6.5.2) \cr
&\sigma(\overline{W}^{(N)}_k, k\ge
1) \vee \sigma(Z^*_j,j\ge 1)\bigr) \cr}$$
is a restricted Markov triplet.

Next, by (6.2.9)--(6.2.10) and a simple argument,
$$\sigma(\kappa(N,j), j\le 0) \dot\subset 
\sigma(W^{(N)}_k,k\le
0)\quad\hbox{and} \leqno{(6.5.3)}$$
$$\sigma(\kappa(N,j),j\ge 1) \dot\subset 
\sigma(W^{(N)}_k, k\ge 1).
\leqno{(6.5.4)}$$

Of course by (6.2.9)--(6.2.10), (6.5.1), and (6.2.7), for a given
$k\le 0$ and a given $\omega\in\Omega_0$ (see (6.3.1)--(6.3.2)),
$$W^{(N+1)}_k (\omega) = \cases{
0 &if $k\not\in \{\kappa(N,0)(\omega),
\kappa(N,-1)(\omega),\kappa(N,-2)(\omega),\dots\,\}$ \cr
Z^*_j(\omega) &if $k=\kappa(N,j)(\omega)$ for some $j\le 0$. \cr}
\leqno{(6.5.5)}$$
Hence by (6.5.3) (and (6.3.6)),
$$\leqalignno{
\sigma(W^{(N+1)}_k,k\le 0)\,\,&\dot\subset\,\, \sigma(\kappa(N,j),
j\le 0) \vee \sigma(Z^*_j,j\le 0) &(6.5.6)\cr
&\dot\subset
\sigma(\overline{W}^{(N)}_k, k\le 0) \vee \sigma(Z^*_j,j\le 0).
\cr}$$
Hence by (6.3.6) again,
$$\sigma(\overline{W}^{(N+1)}_k,k\le 0) \,\,\dot\subset\,\,
\sigma(\overline{W}^{(N)}_k,k\le 0) \vee \sigma(Z^*_j, j\le 0).
\leqno{(6.5.7)}$$

Next, as an analog of (6.5.5), one has that for $k\ge 1$ and
$\omega\in\Omega_0$, $W^{(N+1)}_k(\omega) =0$ resp.\ $Z^*_j(\omega)$
if $k\not\in \{\kappa(N,1)(\omega)$, $\kappa(N,2)(\omega)$,
$\kappa(N,3)(\omega),\dots\,\}$ resp.\ $k = \kappa(N,j)(\omega)$
for some $j\ge 1$. Then by (6.5.4), one obtains an analog of (6.5.6)
with the inequalities $k\le 0$ and $j\le 0$ replaced by $k\ge 1$ and
$j\ge 1$. Thereby one obtains the following analog of (6.3.7):
$$\sigma(\overline{W}^{(N+1)}_k, k\ge 1) \,\,\dot\subset \,\,
\sigma(\overline{W}^{(N)}_k,k\ge 1) \vee \sigma(Z^*_j,j\ge 1).
\leqno{(6.5.8)}$$

Now by Remark 6.3(B) and (6.2.9)--(6.2.10) and a simple argument,
$$\{W^{(N+1)}_0 =1\} \subset \{W_0^{(N)} =1\} \doteq \{\kappa(N,0) =0\}.$$
As a trivial consequence, by (6.5.1) for $j=0$,
$$\leqalignno{\qquad\qquad\{W^{(N+1)}_0 =1\} &\doteq 
\{W^{(N+1)}_0 =1\} \cap
\{W^{(N)}_0 =1\} \cap \{\kappa(N,0) =0\}&(6.5.9)\cr
& = \{W^{(N+1)}_{\kappa(N,0)}
=1\} \cap \{W^{(N)}_0 =1\} \cap \{\kappa(N,0)=0\}  \cr
&= \{Z^*_0 =1\} \cap \{W^{(N)}_0 =1\} \cap \{\kappa(N,0) =0\} \cr
&\doteq \{Z^*_0=1\} \cap \{W^{(N)}_0 =1\}. \cr}$$

Now by (6.5.7), (6.5.8), (6.5.9), and the entire sentence containing
(6.5.2), one has that
$$\left(\sigma(\overline{W}^{(N+1)}_k,k\le 0), \,\, \{W^{(N+1)}_0
=1\}, \,\,  \sigma(\overline{W}^{(N+1)}_k,k\ge 1)\right)$$ is a
restricted Markov triplet. That completes the induction step and the
proof of Lemma~6.5.
\medskip

{\sc Lemma 6.6.} {\it Suppose $n\in\N$. Referring to section 6.3(A)
(equations (6.3.1) and (6.3.2) and statement (i)), Definition 2.3(A),
and Lemma 6.4(C), define the positive integer--valued random variable
$$T^{(n)} := \min\{k\in{\N}: W^{(n)}_k =1\}. \leqno{(6.6.1)}$$
Then }
$$E\left(T^{(n)}\,\Bigl|\, W^{(n)}_0 =1\right) = 16^n, \leqno{(6.6.2)}$$
$$\var\left(T^{(n)}\,\Bigl|\, W^{(n)}_0 =1\right) \le 16^{2n}, \quad
\hbox{\it and (hence)} \leqno{(6.6.3)}$$
$$E\left((T^{(n)})^2 \,\Bigl|\, W^{(n)}_0 =1\right) \le 2\cdot 16^{2n}.
\leqno{(6.6.4)}$$
\medskip

{\sc Proof.} For $n=1$, this holds by Lemma 5.5 and Remark 6.3(F).
Now for induction suppose $N\in\N$ and Lemma 6.6 holds for $n=N$. Our
task is to show that it holds for $n=N+1$.

By (6.2.9)--(6.2.10),
$$T^{(N)} = \kappa(N,1)\ \ {\rm a.s.} \leqno{(6.6.5)}$$
and also
$$\left\{W^{(N)}_0 =1\right\} \doteq \{\kappa(N,0) =0\}. \leqno{(6.6.6)}$$

As in the proof of Lemma 6.5, define the sequence $Z^* := (Z^*_j$,
$j\in\Z)$ of $\{0,1,\dots,6\}$--valued random variables by (6.5.1).
Then (6.5.9) holds. We shall refer to both (6.5.1) and (6.5.9) freely
below.

By (6.5.1), Remark 6.3(G), and Lemma 5.4(A), the sequence $Z^*$
satisfies Condition~$\cS$ a.s.\ (see Definition 2.3(A) again).
Accordingly, define the positive integer--valued random variable $M$
as follows:
$$M:= \min \{j\in{\N}: Z^*_j =1\}. \leqno{(6.6.7)}$$
By (6.5.1), Remark 6.3(G), and Lemma 5.5,
$$E(M\mid Z^*_0 =1) =16 \quad \hbox{and}\quad \var(M\mid Z^*_0 =1) =
80/3. \leqno{(6.6.8)}$$

Next, by (6.5.1), Remark 6.3(H), and (6.2.1), the sequences $W^{(N)}$
and $Z^*$ are independent of each other. It follows from (6.5.9),
(6.6.7), and a simple calculation that $E(M\mid Z^*_0 =1) = E(M\mid
W^{(N+1)}_0 =1)$ and $\var (M\mid Z^*_0 =1) = \var(M\mid W^{(N+1)}_0
=1)$. Hence by (6.6.8),
$$E(M\mid W^{(N+1)}_0 =1) =16 \quad \hbox{and} \quad \var(M\mid
W^{(N+1)}_0 =1) = 80/3. \leqno{(6.6.9)}$$

Now for every $\omega\in\Omega_0$ (see (6.3.1)--(6.3.2) again), one
has the following: For any $j\in\{1,2,\dots,M(\omega) -1\}$ (if
$M(\omega) \ge 2$), $W^{(N+1)}_{\kappa(N,j)(\omega)}(\omega) =
Z^*_j(\omega) \not= 1$ by (6.5.1) and (6.6.7). For any $k\in {\Z} -
\{\kappa(N,j)(\omega): j\in \Z\}$, $W^{(N+1)}_k(\omega) =0$ by
(6.2.9)--(6.2.10) and Remark 6.3(B). Hence 
(for $\omega \in \Omega_0$), for all $k\in
\{1,2,\dots,\kappa(N,M(\omega))(\omega)-1\}$, $W^{(N+1)}_k (\omega)
\not=1$. Also (for $\omega \in \Omega_0$), 
$W^{(N+1)}_{\kappa(N,M(\omega))(\omega)}(\omega) =
Z^*_{M(\omega)}(\omega) =1$ by (6.5.1) and (6.6.7).

By the preceding two sentences (see (6.3.1)--(6.3.2) again) and
(6.6.1),
$$T^{(N+1)} = \kappa(N,M)\ \ \hbox{a.s.} \leqno{(6.6.10)}$$

Now let us look at the random variables 
$\kappa(N,j) - \kappa(N,j-1)$, $j\in\N$. By Lemma 6.5 and the strict stationarity
of the sequence $W^{(N)}$ (recall (6.2.1) and 
Remark 2.6(B)), one has
the following: For any positive integers $J$ and $\ell$, and any
event $A\in \sigma(W^{(N)}_k$, $k\le J-1)$ such that $P(A\cap
\{W^{(N)}_J =1\})>0$, one has that (see Remark (2.7)(C)) 
$$\eqalign{&P\left(W^{(N)}_{J+i} \not= 1 \,\,\forall\,\,
i\in\{1,\dots,\ell-1\}\quad\hbox{and}\quad W^{(N)}_{J+\ell}=1\,\Bigl|\,
A\cap \{W^{(N)}_J =1\}\right) \cr
&=P\left(W^{(N)}_{J+i} \not= 1\,\,\forall\,\,
i\in\{1,\dots,\ell-1\}\quad\hbox{and}\quad W^{(N)}_{J+\ell} =1\,\Bigl|\,
W^{(N)}_J =1\right) \cr
&= P\left(W^{(N)}_i \not= 1 \,\,\forall\,\, i\in \{1,\dots,\ell-1\}
\quad\hbox{and}\quad W^{(N)}_\ell =1 \,\Bigl|\, W^{(N)}_0 =1\right). \cr}$$
(Of course for $\ell =1$, omit the phrases $W^{(N)}_{J+i} \not= 1$
resp.\ $W^{(N)}_i \not= 1$ $\forall\,\, i\in \{1,\dots,\ell-1\}$.)
Hence (see (6.2.9)--(6.2.10) again), by a standard induction
argument, one can show that conditional on the event $\{W^{(N)}_0
=1\}$, the random variables $\kappa(N,1)$ 
(or $\kappa(N,1) - \kappa(N,0)$ --- see (6.6.6)), 
$\kappa(N,2) - \kappa(N,1)$,
$\kappa(N,3) - \kappa(N,2)$, $\kappa(N,4) - \kappa(N,3),\dots$ are
independent and identically distributed, with (see (6.6.6), (6.6.5),
and (6.2.9)--(6.2.10))
$$\eqalign{
\cL(\kappa(N,1) - \kappa(N,0) \mid W^{(N)}_0 =1) &=
\cL(\kappa(N,1) \mid W^{(N)}_0 =1)\cr &=
\cL(T^{(N)} \mid W^{(N)}_0 =1).\cr}$$ Hence by (6.5.9),
the sentence after (6.6.8) (recall that $\sigma(\kappa(N,i)) \dot\subset
\sigma(W^{(N)})$ for $i\in\Z$ by (6.2.9)--(6.2.10)), and a standard
trivial calculation, conditional on the event $\{W^{(N+1)}_0 =1\}$,
the random variables $\kappa(N,j) -\kappa(N,j-1)$, $j\in\N$ are
independent and identically distributed, with
$$\leqalignno{\qquad\quad&\cL(\kappa(N,1) - \kappa(N,0) \mid W^{(N+1)}_0 =1) &(6.6.11) \cr
&=
\cL(\kappa(N,1) - \kappa(N,0) \mid W^{(N)}_0 =1) = \cL(T^{(N)}\mid
W^{(N)}_0 =1). \cr}$$
Hence by the induction hypothesis of (6.6.2) and (6.6.3) for $n=N$,
$$E(\kappa(N,1) - \kappa(N,0) \mid W^{(N+1)}_0 =1) =16^N\quad
\hbox{and} \leqno{(6.6.12)}$$
$$\var(\kappa(N,1) - \kappa(N,0) \mid W^{(N+1)}_0 =1) \le 16^{2N}.
\leqno{(6.6.13)}$$

Now recall again from (6.2.9)--(6.2.10) that 
$\sigma(\kappa(N,i),
i\in{\Z})\dot\subset \sigma(W^{(N)})$. 
By that fact, (6.6.7), and the
sentence after (6.6.8), and equation (6.5.9), together with a
standard simple argument, conditional on the event $\{W^{(N+1)}_0
=1\}$, the random variable $M$ is independent of the random sequence
$(\kappa(N,j) - \kappa(N,j-1)$, $j\in\N)$. Hence by (6.6.10),
(6.6.9), (6.6.12), (6.6.13), the entire sentence containing (6.6.11),
and a well known elementary calculation (see e.g.\ [15, p.\ 301,
Exercise 1]), one has that (recall (6.6.6) and (6.5.9))
$$\eqalign{&E(T^{(N+1)}\mid W^{(N+1)}_0 =1) =
E(\kappa(N,M) \mid W^{(N+1)}_0 =1) \cr
&=E\biggl(\, \sum^M_{j=1} [\kappa(N,j) - \kappa(N,j-1) ]
\,\biggl|\, W^{(N+1)}_0
=1\biggl) \cr
&= E(M \mid W^{(N+1)}_0 =1) \cdot E(\, \kappa(N,1) -
\kappa(N,0) \mid W^{(N+1)}_0 =1) \cr
&= 16 \cdot 16^N = 16^{N+1} \cr}$$
and
$$\eqalign{
&\var(T^{(N+1)} \mid W^{(N+1)}_0 =1) =
\var(\kappa(N,M)\mid W^{(N+1)}_0 =1) \cr
&=\var\biggl(\, \sum^M_{j=1} [\kappa(N,j) -\kappa(N,j-1)]\,\biggl|\,
W^{(N+1)}_0 =1]\biggl) \cr &= 
E(M\mid W^{(N+1)}_0 =1)
\cdot \var (\kappa(N,1) - \kappa(N,0)\mid W^{(N+1)}_0 =1)
\cr &\qquad\quad +\var(M\mid W^{(N+1)}_0 =1) \cdot
\bigl[E( \kappa(N,1) - \kappa(N,0) \mid W^{(N+1)}_0
=1)\bigl]^2 \cr &\le 16\cdot 16^{2N} +(80/3) \cdot (16^N)^2
<16^{2(N+1)}. \cr}$$ Thus (6.6.2) and (6.6.3) (and hence also
(6.6.4)) hold for $n=N+1$. That completes the induction step and the
proof of Lemma 6.6.
\medskip

{\sc Lemma 6.7.} {\it Suppose $n\in\N$. Let $p_n$ denote the
probability that there exist at least two distinct integers $i,j
\in\{1,2,\dots,6\cdot 16^n\}$ such that $W^{(n)}_i = W^{(n)}_j =1$.
Then $p_n \ge 1/2$.}
\medskip

{\sc Proof.}  Suppose $n\in\N$. Refer to (6.2.9)--(6.2.10). 
It
suffices to prove that $P(\kappa(n,2) \le 6\cdot 16^n) \ge 1/2$, or
$$P(\kappa(n,2) >6 \cdot 16^n) \le 1/2.
\leqno{(6.7.1)}$$

Recall (from (6.2.1) and Remark 2.6(B)) that the sequence $W^{(n)}$
is strictly stationary. Recall from (6.2.9)--(6.2.10) that
$\kappa(n,1) = T^{(n)}$\ a.s., where $T^{(n)}$ is as in Lemma 6.6, and that
$\{\kappa(n,0) =0\} \doteq \{W^{(n)}_0 =1\}$. From (6.2.9)--(6.2.10) and
a trivial argument, followed by Lemma 6.4(C), one has that for each
positive integer $j$,
$$\leqalignno{
\qquad \quad P(\kappa(n,1) &=j) = \sum^\infty_{\ell=0} P(\kappa(n,0)
= -\ell \,\,\hbox{and}\,\, \kappa(n,1) =j) &(6.7.2) \cr &=
\sum^\infty_{\ell=0} P(\kappa(n,0) =0,\, \kappa(n,1) = j+\ell) \cr
&=\sum^\infty_{\ell=0} P(\kappa(n,1) = j+\ell\mid \kappa(n,0) =0)
\cdot P(\kappa(n,0) =0) \cr &= P(W^{(n)}_0 =1) \cdot
P(\kappa(n,1) \ge j\mid \kappa(n,0) =0) \cr &= 16^{-n} \cdot
P(T^{(n)}\ge j \mid W^{(n)}_0 =1). \cr}$$

Now for (say) any positive integer--valued random variable $Z$, one
has by a simple argument that $\sum^\infty_{j=1} j\cdot P(Z \ge j)
\le EZ^2$. Hence by (6.7.2) and Lemma 6.6,
$$\leqalignno{
\qquad \quad E\kappa(n,1) &= \sum^\infty_{j=1} j\cdot 16^{-n} \cdot
P(T^{(n)} \ge j\mid W^{(n)}_0 =1) &(6.7.3)\cr
 \le \,
&16^{-n} \cdot E\left((T^{(n)})^2
\, \Bigl| \, W^{(n)}_0 =1\right) \le
16^{-n} \cdot 2\cdot 16^{2n} = 2\cdot 16^n.  \cr}$$

Now by strict stationarity of the sequence $W^{(n)}$, and a standard
elementary argument using Remark 2.8 (akin to certain arguments in
the proof of Lemma 6.6), the random variable $\kappa(n,2) -
\kappa(n,1)$ is independent of $\kappa(n,1)$ and $\cL(\kappa(n,2) -
\kappa(n,1)) = \cL(\kappa(n,1) \mid W^{(n)}_0 =1)$. Recall again from
the paragraph after (6.7.1) that $\kappa(n,1) = T^{(n)}$, from Lemma
6.6. One  now has by Lemma 6.6 that
$$E(\kappa(n,2) - \kappa(n,1)) = E(T^{(n)} \mid W^{(n)}_0 =1) =
16^n.$$
Hence by (6.7.3), $E\kappa(n,2) \le 3\cdot 16^n$. Now (6.7.1) holds
by Markov's inequality. Lemma 6.7 is proved.
\bigskip

\centerline {\bf 7. Scaffolding (part~2)} 
\bigskip

Among other things, this section will
include (in Construction 7.5 below) the construction of the random
sequence $X$ itself for Theorem 1.1. This section will build on
Section~6. Recall from that section the given probability space
$(\Omega,\cF,P)$. \medskip

{\sc Construction 7.1.} Refer to Construction 6.2, including equation
(6.2.6). Define the $\{0,1,2,3,4,5,6\}$--valued random variables
$\delta^{(n)}_k$, $n\in\N$, $k\in\Z$ as follows:
$$\delta^{(1)}_k := W^{(1)}_k
\leqno{(7.1.1)}$$
and for $n\ge 2$,
$$\delta^{(n)}_k := W^{(n)}_{k- \Psi(n-1,k,0)}.
\leqno{(7.1.2)}$$
By (7.1.1), (7.1.2), and (6.2.6), one has that 
for all $k\in\Z$,
$$\sigma(\delta_k^{(1)}) \subset \sigma(W^{(1)}), \quad 
{\rm and} \quad
\forall n \geq 2,\ \ \sigma(\delta^{(n)}_k) \subset
\sigma(W^{(n-1)},W^{(n)}).
\leqno{(7.1.3)}$$

Define the $(\{0\} \cup {\N} \cup \{\infty\})$--valued random variables
$N_k$, $k\in\Z$ as follows, for each $\omega\in\Omega$,
$$N_k(\omega) := \cases{
0 &if $\delta^{(1)}_k(\omega) =0$ (that is, $W^{(1)}_k(\omega) =0$)
\cr m\in \N &if $\delta^{(u)}_k(\omega) \not= 0$ $\forall$
$u\in\{1,\dots,m\}$ and $\delta^{(m+1)}_k(\omega) =0$ \cr \infty &if
$\delta^{(u)}_k(\omega) \not= 0$ for all $u\in\N.$ \cr}
\leqno{(7.1.4)}$$

By (7.1.4) and (7.1.3) (and (6.3.6)), for each $m\in\N$ and each
$k\in\Z$,
$$\{N_k\ge m\} = \bigcap^m_{u=1} \{\delta^{(u)}_k \not= 0\} \subset
\sigma(\overline{W}^{(m)}); \leqno{(7.1.5)}$$
and for each integer $m\ge 0$,
$$\{N_k =m\} \subset \sigma(\overline{W}^{(m+1)}).
\leqno{(7.1.6)}$$

Also, define the integer-valued random variables $J(m,k)$, $m\in\N$,
$k\in\Z$ as follows:
$$J(m,k) := \sum^m_{u=1} 6^{u-1}\left(\delta^{(u)}_k -1\right).
\leqno{(7.1.7)}$$
By (7.1.3), for each $m\in\N$ and each $k\in\Z$,
$$\sigma(J(m,k)) \subset \sigma(\overline{W}^{(m)}).
\leqno{(7.1.8)}$$

{\sc Remarks 7.2.} Recall from Construction 7.1 that the random
variables $\delta^{(n)}_k$, $n\in \N$, $k\in\Z$ take their values in
the set $\{0,1,\dots,6\}$.

(A) Suppose $m\in\N$, $\omega\in\Omega$, $k\in\Z$, and
$W^{(m)}_k(\omega) =1$; then (i)~$W^{(n)}_k(\omega) =1$ for all
$n\in\{1,\dots,m\}$ by Remark 6.3(B), (ii)~for each
$n\in\{1,\dots,m\}$, $\Psi(n,k,0)(\omega) =0$ by (6.2.6) and 
Remark 2.2(B), 
hence (iii)~$\delta^{(n)}_k(\omega) =1$ for all
$n\in\{1,\dots,m\}$ by (7.1.1) and (7.1.2), and hence
(iv)~$N_k(\omega) \ge m$ by (7.1.5).

(B) For each $m\in\N$, each $k\in\Z$, and each $\omega\in\Omega$ such
that $N_k(\omega) \ge m$, one has that (i)~$\delta^{(u)}_k(\omega)
\in \{1,2,\dots,6\}$ for all $u\in\{1,\dots,m\}$ by (7.1.5), hence
(ii)~$\delta^{(u)}_k(\omega) -1 \in \{0,1,\dots,5\}$ for all $u\in
\{1,\dots,m\}$, and hence (iii)~$J(m,k)(\omega) \in \{0,1,2,\dots,6^m
-1\}$ by (7.1.7) and a simple argument.
\medskip

{\sc Construction 7.3.} (A) Refer to Definition 3.3. On the given
probability space $(\Omega,\cF,P)$, let $\zeta^{(n,\ord)}_k$,
$\zeta^{(n,\cen)}_k$, $\zeta^{(n,\fri)}_k$, $n\in\N$, $k\in\Z$ be an
array of independent random variables, with this array being
independent of the entire collection of random variables
$\xi^{(n)}_k$ $W^{(n)}_k$, $\Psi(n,k,j)$, $\delta^{(n)}_k$, $N_k$,
$J(n,k)$, $n\in\N$, $k\in\Z$, $j\in \{0,1,2,\dots\,\}$ in Section 6
and Construction 7.1 (the redundancy here is for emphasis), such
that for each $n\in\N$ and each $k\in\Z$, (i)~all three random
vectors $\zeta^{(n,\ord)}_k$, $\zeta^{(n,\cen)}_k$,
$\zeta^{(n,\fri)}_k$ take their values in the set
$\{-1,1\}^{\sxtp(n)}$ (see (2.1)), and (ii)~the distribution of
$\zeta^{(n,\ord)}_k$ resp.\ $\zeta^{(n,\cen)}_k$ resp.\
$\zeta^{(n,\fri)}_k$ is $\nu^{(n)}_{\ord}$ resp.\ $\nu^{(n)}_{\cen}$
resp.\ $\nu^{(n)}_{\fri}$.

(B) For a given $n\in\N$, the random vector $\zeta^{(n,\ord)}_k$ will
be represented by
$$\zeta^{(n,\ord)}_k :=
\left(\zeta^{(n,\ord)}_{k,0},\zeta^{(n,\ord)}_{k,1},\dots,\zeta^{(n,\ord)}_{k,\sxtp(n)-1}\right),$$
that is, with the $6^n$ indices running through $0,1,\dots,6^n-1$
(instead of $1,2,\dots,6^n$); and exactly the same convention will be
used for the random vectors $\zeta^{(n,\cen)}_k$ and
$\zeta^{(n,\fri)}_k$. (This fits the convention in Section 3 where,
for a given $n\in\N$, the elements $x\in\{-1,1\}^{\sxtp(n)}$ were
represented as $x:= (x_0,x_1,\dots,x_{\sxtp(n)-1})$.

(C) For each $n\in\N$ and each $k\in\Z$, define the
$(\{-1,1\}^{\sxtp(n)})^3$--valued random vector
$$\zeta^{(n)}_k :=
\left(\zeta^{(n,\ord)}_k,\zeta^{(n,\cen)}_k,\zeta^{(n,\fri)}_k\right).
\leqno{(7.3.1)}$$
Of course these random vectors are independent and identically
distributed. For each $n\in\N$, define the resulting random sequence
$\zeta^{(n)} := (\zeta^{(n)}_k$, $k\in\Z)$.

(D) For each $n\in\N$ and each $k\in\Z$, define the $(\{-1,1\}^6)^3
\times (\{-1,1\}^{36})^3 \allowbreak \times \dots \times
(\{-1,1\}^{\sxtp(n)})^3$--valued random vector
$$\overline{\zeta}^{(n)}_k := \left(\zeta^{(1)}_k,
\zeta^{(2)}_k,\dots,\zeta^{(n)}_k\right). \leqno{(7.3.2)}$$
For any given fixed value of $n\in\N$, these random vectors
$\overline{\zeta}^{(n)}_k$, $k\in\Z$ are independent and identically
distributed. For each $n\in\N$, define the resulting random sequence
$\overline{\zeta}^{(n)} := (\overline{\zeta}^{(n)}_k$, $k\in\Z)$.

(E) For each $k\in\Z$, define the random item (sequence)
$$\overline{\zeta}^{(\infty)}_k := \left(\zeta^{(1)}_k,
\zeta^{(2)}_k, \zeta^{(3)}_k,\dots\,\right). \leqno{(7.3.3)}$$
These random items $\overline{\zeta}^{(\infty)}_k$, $k\in \Z$ are
independent and identically distributed. Define the resulting random
sequence $\overline{\zeta}^{(\infty)} :=
(\overline{\zeta}^{(\infty)}_k$, $k\in\Z)$.

(F) Refer to (B) above.  Purely as a formality, for each 
$n \in \N$, each
$k \in \Z$, and each $u \in {\Z}-\{0, 1, \dots, 6^n-1\}$,
define the constant random variables
$\zeta_{k,u}^{(n, {\rm ord})}
:= \zeta_{k,u}^{(n, {\rm cen})}
:= \zeta_{k,u}^{(n, {\rm fri})} 
:= 1$.  (This formality will ultimately turn out to be frivolous.)

\medskip

{\sc Construction 7.4.} (A) Let $X^{(0)} := (X^{(0)}_k$, $k\in\Z)$ be
a sequence of independent, identically distributed $\{-1,1\}$--valued
random variables such that (for each $k\in\Z$),
$$P\left(X^{(0)}_k = -1\right) = P\left(X^{(0)}_k =1\right) = 1/2,
\leqno{(7.4.1)}$$
with this sequence $X^{(0)}$ being independent of the entire array of
random variables $\xi^{(n)}_k$, $\zeta^{(n)}_k$, $n\in\N$, $k\in\Z$
(and hence independent of the entire collection of random variables
in Section 6 and Constructions 7.1 and 7.3).

(B) Refer to (A) above and to (6.1.4) and (7.3.3). For each $k\in\Z$,
define the random ordered triplet
$$\eta_k := \left(X^{(0)}_k, \overline{\xi}^{(\infty)}_k,
\overline{\zeta}^{(\infty)}_k\right). \leqno{(7.4.2)}$$
Note that for any given $k\in\Z$, the three components of $\eta_k$
are independent of each other. These random ordered triplets
$\eta_k$, $k\in\Z$ are independent and identically distributed.
Define the resulting random sequence $\eta:= (\eta_k$, $k\in\Z)$.
\medskip

{\sc Construction 7.5.} Now on our given probability space
$(\Omega,\cF,P)$, define the sequence $X:= (X_k$, $k\in\Z)$ of
$\{-1,1\}$--valued random variables for Theorem 1.1 as follows: For
each $k\in\Z$ and each $\omega\in\Omega$ (see Remark 7.2(B)),
$$X_k(\omega):= \cases{
X^{(0)}_k(\omega) &if $N_k(\omega) =0$ \cr
\zeta^{(\ell,\cen)}_{k-\Psi(\ell,k,0)(\omega),J(\ell,k)(\omega)}(\omega)
&if $N_k(\omega) =\ell\in\N$ \cr 1 &if $N_k(\omega) =\infty$. \cr}
\leqno{(7.5.1)}$$
\medskip

{\sc Lemma 7.6.} {\it Refer to Constructions 7.4 and 7.5 and Definition
2.6(A). The ordered pair $(X,\eta)$ satisfies Condition~$\cM$.}
\medskip

{\sc Proof.} In what follows, keep in mind the three sentences after
(7.4.2).

For each $n\in\N$, the ordered pair $(W^{(n)}, \eta)$ satisfies
Condition~$\cM$ by (6.1.3), (6.1.4), (6.2.1), and (7.4.2). Hence for
each $n\in\N$ and each $j\ge 0$, the ordered pair $(\Psi(n,j), \eta)$
satisfies Condition~$\cM$ by (6.2.4) (see (6.2.6) and the phrase
right before (6.2.4)).

Next, for each $n\ge 2$ and each $k\in\Z$,
$$W^{(n)}_{k-\Psi(n-1,k,0)} = \sum^\infty_{u=0} W^{(n)}_{k-u}
I(\Psi(n-1,k,0) =u).$$ Hence by (7.1.1), (7.1.2), and both sentences
in the preceding paragraph above, for each $n\in\N$, the ordered
pair $((\delta^{(n)}_k$, $k\in\Z),\eta)$ satisfies Condition~$\cM$.

Next, by (7.1.4), for each $k\in\Z$,
$$\eqalign{N_k =\,\, &0 \cdot I(\delta^{(1)}_k =0) \cr
&+ \sum_{m\in\N} m\cdot\left[ \prod^m_{u=1} I\left(\delta^{(u)}_k
\not= 0\right)\right] \cdot I\left(\delta^{(n+1)}_k =0\right) \cr
&+\infty \cdot \prod_{u\in\N} I\left(\delta^{(u)}_k \not= 0\right).
\cr}$$ Hence by the last sentence of the preceding paragraph, the
ordered pair $((N_k$, $k\in\Z),\eta)$ satisfies Condition $\cM$.
Similarly, from (7.1.7), for each $m\in\N$, the ordered pair
$((J(m,k)$, $k\in\Z),\eta)$ satisfies condition~$\cM$.

Next, for each $\ell\in\N$ and each $k\in\Z$, by Remark 7.2(B) (see
also (2.1)),
$$\leqalignno{\qquad\qquad &I(N_k =\ell) \cdot
\zeta^{(\ell,\cen)}_{k-\Psi(\ell,k,0),J(\ell,k)} &(7.6.1) \cr
&=\sum^\infty_{u=0}\sum^{\sxtp(\ell)-1}_{v=0}\!\!\!\! I(N_k =\ell)
\cdot \zeta^{(\ell,\cen)}_{k-u,v} \cdot I(\Psi(\ell,k,0) =u) \cdot
I(J(\ell,k) =v). \cr}$$ 
(If $\omega\in\Omega$ is such that
$N_k(\omega) \not= \ell$, then trivially (7.6.1) holds 
for that $\omega$ with
both sides being 0; for the formal definition of the
left side of (7.6.1) for such $\omega$, recall Construction 7.3(F)
to cover the possible case $J(\ell,k)(\omega) < 0$.)\ \  
By (7.6.1) and the
observations made so far, together with equations (7.3.1), (7.3.3),
and (7.4.2), for each $\ell\in\N$, the ordered pair
$$\left(\left( I(N_k =\ell) \cdot
\zeta^{(\ell,\cen)}_{k-\Psi(\ell,k,0),J(\ell,k)},
k\in\Z\right),\eta\right)$$ satisfies Condition $\cM$.

Finally, for each $k\in\Z$, by (7.5.1),
$$\eqalign{X_k &= X^{(0)}_k \cdot I(N_k =0) + 1 \cdot I(N_k =\infty)
\cr &\qquad +\sum_{\ell\in\N} I(N_k =\ell) \cdot
\zeta^{(\ell,\cen)}_{k-\Psi(\ell,k,0), J(\ell,k)}. \cr}$$ Hence by
observations in the preceding two paragraphs together with (7.4.2),
the ordered pair $(X,\eta)$ satisfies Condition $\cM$. Lemma 7.6 is
proved.
\medskip

{\sc Remark 7.7.} Of course by Lemma 7.6 (and Definition 2.6(A)), the
random sequence $X$ is strictly stationary. Also, by (7.4.2),
(7.4.1), (7.3.3), (7.3.1), and (6.1.4), the random variables
$\eta_k$, $k\in\Z$ in (7.4.2) can be regarded as taking their values
in the set $S:= \{-1,1\}\times\{0,1\}^{\N} \times\{-1,1\}^{\N}$.
Trivially that set is bimeasurably isomorphic to the set
$\{0,1\}^{\N}$, and that set in turn is well known to be bimeasurably
isomorphic to the open unit interval $(0,1)$ (with its Borel
$\sigma$-field) and hence also to the real line $\R$ (with its Borel
$\sigma$-field $\cR$). Applying a particular bimeasurable isomorphism
$\Theta : S\to\R$ to each of the random variables $\eta_k$ in
(7.4.2), one has that the random sequence $\eta$ can thereby be
``coded'' as a sequence of independent, identically distributed {\it
real-}valued random variables. Thus by Lemma 7.6, property (C) in
Theorem 1.1 holds.

The other properties in Theorem 1.1 will be verified in sections 9
and 10, after some further preparation in Section~8. The following
technical lemma will be needed in Section~9.
\medskip

{\sc Lemma 7.8.} {\it For each $n\in\N$ and each $k\in\Z$, $P(N_k\ge
n) =(3/8)^n$.}
\medskip

{\sc Proof.} For each $k\in\Z$, by (7.1.4) (or (7.1.5)), (7.1.1), and
Lemma 6.4(A) (and stationarity),
$$P(N_k \ge 1) = P(\delta^{(1)}_k \not= 0) =
P(W^{(1)}_k \not= 0) = 3/8.$$

Now for induction, suppose that $n\in\N$, $k\in\Z$, and $P(N_k \ge n)
= (3/8)^n$. If it is shown that $P(N_k \ge n+1\mid N_k\ge n) = 3/8$,
then (since $\{N_k \ge n+1\} \subset \{N_k \ge n\}$) it will follow
that $P(N_k \ge n+1) = (3/8)^{n+1}$. Then Lemma 7.8 will hold by
induction.

Suppose $j\in \{0,1,2,\dots\,\}$ and that $P(\{N_k \ge n\}\cap
\{\Psi(n,k,0) =j\})>0$. If $\omega\in\Omega_0$ (see
(6.3.1)--(6.3.2)) is such that $\Psi(n,k,0)(\omega) =j$, then
$W^{(n)}_{k-j}(\omega) =1$ by (6.2.6) and Definition 2.2 (and
statement 6.3(i) and Definition 2.3(A)). Also, $\{N_k \ge n\}\in
\sigma(\overline{W}^{(n)})$ by (7.1.5), and $\{\Psi(n,k,0) =j\}\in
\sigma(\overline{W}^{(n)})$ by (6.2.6). Hence by (7.1.4) (see the
equality in (7.1.5) with $m=n$ and with $m=n+1$), (6.3.2), (7.1.2),
and Lemma 6.4(D),
$$\eqalign{
&P\left(N_k \ge n+1\mid \{N_k\ge n\}\cap \{\Psi(n,k,0,) =j\}\right)
\cr &=P\left(\delta^{(n+1)}_k \not= 0 
\, \Bigl| \, \{N_k \ge n\} \cap
\{\Psi(n,k,0) =j\}\cap \{W^{(n)}_{k-j} =1\}\right) \cr
&=P\left(W^{(n+1)}_{k-j} \not= 0
\, \Bigl| \, \{N_k \ge n\} \cap
\{\Psi(n,k,0) =j\}\cap \{W^{(n)}_{k-j} =1\}\right) \cr &=
3/8. \cr}$$

Hence by Remark 2.8, $P(N_k \ge n+1\mid N_k \ge n) =3/8$. 
That completes the induction argument and the proof.
\bigskip

\centerline {\bf 8. Scaffolding (part~3)}
\bigskip 

In this section, some more foundations
will be laid for the proofs, in sections 9 and 10, of the properties
in Theorem 1.1 not verified in Remark 7.7.
\medskip

{\sc Definition 8.1.} Refer to (7.1.4), Section 6.3(A), and Definition
2.3(A).  For each $m\in\N$ and each $\omega\in\Omega_0$, let
$\cE_m(\omega)$ denote the family of all sets $E\subset \Z$ such that
the following holds:

There exist integers $j$ and $\ell$ such that
$$\leqalignno{j& < \ell, \quad 
W^{(m)}_j(\omega) = W^{(m)}_\ell(\omega)
=1,\quad \hbox{and}   &(8.1.1) \cr
& W^{(m)}_k (\omega) \not= 1\quad\hbox{for
all}\quad k\in\{j+1,j+2,\dots,\ell-1\},
\cr }$$
and
$$E=\{k\in\{j,j+1,\dots,\ell-1\}: N_k(\omega) \ge m\}.
\leqno{(8.1.2)}$$

(Note that the ``inner set'' in (8.1.2) contains $j$ but not~$\ell$.)
\medskip

{\sc Remark 8.2.} (A) Refer to Definition 8.1. Suppose $m\in\N$,
$\omega\in\Omega_0$, $E\in \cE_m(\omega)$, and $j$ and $\ell$ are
integers such that (8.1.1) and (8.1.2) hold. Since $W^{(m)}_j(\omega)
=1$ (see (8.1.1)), one has that $N_j(\omega) \ge m$ by Remark
7.2(A)(iv), and hence $j\in E$ (see (8.1.2)). In fact (for the $j$ in
(8.1.1))
$$j= \min E. \leqno{(8.2.1)}$$
Thus trivially the set $E$ is nonempty and also (see (8.1.2)) finite.

(B) For a given $m\in\N$, $\omega\in\Omega_0$, and $E\in
\cE_m(\omega)$, the integers $j$ and $\ell$ in (8.1.1) are unique, by
(8.2.1) and (8.1.1) itself.

(C) If $m\in\N$, $\omega\in \Omega_0$, $E\in \cE_m(\omega)$, and also
$\widetilde E \in \cE_m(\omega)$, then (see section 2.1(D)) either
$E<\widetilde E$, $E=\widetilde E$, or $E>\widetilde E$.

(D) Refer to Definition 8.1 and Section 2.1(H). Suppose $m\in\N$ and
$\omega\in \Omega_0$. Then (recall Statement 6.3(A)(i))
$$\left\{ k\in{\Z}: N_k (\omega) \ge m\right\} = \hbox{union}\,\,
\cE_m(\omega). \leqno{(8.2.2)}$$
From this and (C) above and a simple argument (recall Statement 6.3(A)(i)), (i)~one has a representation of the 
form $\cE_m(\omega)
:= \{\dots,E_{-1},E_0,E_1,\dots\,\}$ (the $E_i$'s depend on 
$m$ and
$\omega$) where (see section 2.1(D)) $\dots< E_{-1}<E_0<E_1<\dots\,$,
and (ii)~those sets $E_i$, $i\in\Z$ form a partition of the set
$\{k\in{\Z}: N_k (\omega) \ge m\}$.

(E) Suppose $m\in\N$, $\omega\in\Omega_0$, and $E\in \cE_m(\omega)$.
(i)~If $W^{(m+1)}_{\min E}(\omega) =0$, then $N_k(\omega) =m$ for all
$k\in E$. (ii)~If instead $W^{(m+1)}_{\min E}(\omega) \not= 0$, then
$N_k(\omega) \ge m+1$ for all $k\in E$.
(Note that by (i) and (ii) together, the respective converses of (i) and (ii) each hold.)
\medskip

{\it Proof of (E).} Let the integers $j$ and $\ell$ be as in
(8.1.1)--(8.1.2). For each $k\in \{j,j+1,\dots,\ell-1\}$,
$\Psi(m,k,0)(\omega) = k-j$ by (8.1.1)--(8.1.2), (6.2.6), and
Definition 2.2, hence $j= k-\Psi(m,k,0)(\omega)$, and hence
$\delta^{(m+1)}_k(\omega) = W^{(m+1)}_j (\omega)$ by (7.1.2). Of
course for each $k\in E$, one has that $\delta^{(n)}_k(\omega) \not=
0$ for all $n\in \{1,\dots,m\}$ by (8.1.2) and (7.1.4). If
$W^{(m+1)}_j(\omega) =0$, then for all $k\in E$,
$\delta^{(m+1)}_k(\omega) =0$ and hence (see (7.1.4)) $N_k(\omega)
=m$. If instead $W^{(m+1)}_j(\omega) \not= 0$, then for all $k\in
E$, $\delta^{(m+1)}_k(\omega) \not= 0$ and hence (see (7.1.4))
$N_k(\omega) \ge m+1$. Thus (recall (8.2.1) once more) statement
(E)(i)(ii) holds.
\medskip

{\sc Definition 8.3.} (A) Refer to (6.3.1)--(6.3.2) and (7.1.4). For
each $\omega\in \Omega_0$, let $\cD_0(\omega)$ denote the family of
all ``singleton'' sets $\{k\}$ (with $k\in \Z$) such that
$N_k(\omega) =0$.

(B) For each $m\in\N$ and each $\omega\in \Omega_0$, let
$\cD_m(\omega)$ denote the family of all sets $E\in \cE_m(\omega)$
such that (see Remark 8.2(E)) $N_k(\omega) =m$ for all $k\in E$.
\medskip

{\sc Remark 8.4.} (A) Refer to Definition 8.3(B), Remark 8.2(D)(E),
and section 2.1(H). For each $m\in\N$ and each $\omega\in \Omega_0$,
$$\{k\in{\Z}: N_k(\omega) =m\} = \hbox{union}\,\, \cD_m(\omega),
\leqno{(8.4.1)}$$
and in fact the members of $\cD_m(\omega)$ form a partition of the
set $\{k\in{\Z}: N_k(\omega) =m\}$. That sentence also holds trivially
for $m=0$ and $\omega\in \Omega_0$; see Definition 8.3(A).

(B) Refer to Remark 8.2(A), Definition 8.3(A)(B), and equations
(8.2.2) and (8.4.1). (i)~For a given $\omega\in \Omega_0$, none of
the families $\cD_m(\omega)$ $(m\ge 0)$ or $\cE_m(\omega)$ $(m\ge
1)$ contains the empty set as a member. (ii)~If $\omega\in \Omega_0$
and $0\le m<n$, then for any $D\in \cD_m(\omega)$ and any $E\in
\cE_n(\omega)$ (in particular, any $E\in \cD_n(\omega))$, one has
that $D\cap E = \emptyset$.

(C) By (A) and (B) above and Remark 8.2(C)(D), the following holds:
If $m\in\N$ and $\omega\in \Omega_0$, then the family
$$\left[ \cD_0(\omega) \cup \cD_1(\omega) \cup \dots \cup
\cD_{m-1}(\omega)\right] \cup \cE_m(\omega)$$ gives a partition of
the set $\Z$ itself into countably many nonempty finite sets.
\medskip

{\sc Lemma 8.5.} {\it Suppose $m\in\N$, $\omega\in \Omega_0$, and
$E\in \cE_m(\omega)$. Then the following statements hold:

(A) One has that $\card\, E = 6^m$.

(B) Representing the set $E$ by $E= \{i(1),i(2),\dots,i(6^m)\}$
where $i(1)<i(2)<i(3)<\dots<i(6^m)$,
one has (see (7.1.7)) that $J(m, i(v))(\omega) = v-1$ for each $v\in
\{1,2,\dots,6^m\}$.

(C) If also $m\ge 2$, then there exist six sets
$E_1,E_2,\dots,E_6\in \cE_{m-1}(\omega)$ such that $E_1<E_2<\dots
<E_6$ (see Section 2.1(D)) and $E= E_1\cup E_2\cup \dots \cup E_6$.}
\medskip

{\sc Proof.} We shall first prove statements (A) and (B) for the case
$m=1$.

Suppose $\omega\in \Omega_0$ and $E\in \cE_1(\omega)$.

Let $j$ and $\ell$ denote the integers such that (8.1.1) and (8.1.2)
hold (with $m=1$).

For a given $k\in\Z$, the following three inequalities are equivalent
by (7.1.4) and (7.1.1): $N_k(\omega) \ge 1$, $\delta^{(1)}_k(\omega)
\not= 0$, and $W^{(1)}_k(\omega) \not= 0$. Hence by (8.1.2),
$$E = \{k\in \{j,j+1,\dots,\ell-1\}: W^{(1)}_k(\omega) \not= 0\}.
\leqno{(8.5.1)}$$

Since $\omega\in \Omega_0$ (by hypothesis), the sequence
$(W^{(1)}_k(\omega)$, $k\in \Z)$ of elements of $\{0,1,\dots,6\}$
satisfies Condition $\cS$ (see Statement 6.3(A)(i) and Definition
2.3(A) again). Hence by (8.1.1), there exist integers
$i(1),i(2),\dots,i(6)$ such that (see also (8.2.1))
$j=i(1)<i(2)<\dots <i(6)<\ell$,
$$W^{(1)}_{i(v)}(\omega) =v\quad \hbox{for each}\quad v\in
\{1,2,\dots,6\},
\leqno{(8.5.2)}$$
and $W^{(1)}_k(\omega) =0$ for all other elements $k\in
\{j,j+1,\dots,\ell-1\}$. Hence by (8.5.1), $E=
\{i(1),i(2),\dots,i(6)\}$. Hence card $E=6$. Also, by (7.1.7),
(7.1.1), and (8.5.2), for each $v\in \{1,2,\dots,6\}$,
$$J(1,i(v))(\omega) = \delta^{(1)}_{i(v)}(\omega)-1 =
W^{(1)}_{i(v)}(\omega) -1 = v-1.$$
All parts of statements (A) and (B) (of Lemma 8.5) have now been
verified for the case $m=1$.

{\it The induction step.} Now suppose $M\in\N$, and statements (A)
and (B) in Lemma 8.5 hold for the case $m=M$. To complete the proof
of Lemma 8.5 by induction, it suffices to prove that all three
statements (A), (B), and (C) hold for the case $m=M+1$.

Suppose $\omega\in \Omega_0$ and $E\in \cE_{M+1}(\omega)$.

Referring to (8.1.1) and (8.1.2), let $j$ and $\ell$ be the integers
such that (for the given $\omega$ and $E$)
$$j<\ell,\quad W^{(M+1)}_j(\omega) = W^{(M+1)}_\ell (\omega) =1,\quad
\hbox{and} \leqno{(8.5.3)}$$
$$W^{(M+1)}_k (\omega) \not= 1 \quad \hbox{for all}\quad k\in
\{j+1,j+2,\dots,\ell-1\}, $$ and
$$E = \{k\in \{j,j+1,\dots,\ell-1\}:N_k(\omega) \ge M+1\}.
\leqno{(8.5.4)}$$

By (8.5.3) and Remark 6.3(B), $W^{(M)}_j (\omega) = W^{(M)}_\ell
(\omega) =1$. Let $a(0)$, $a(1)$, $\dots$, $a(p)$ (where $p$ is a positive
integer) denote the integers such that
$$j=a(0)<a(1)<a(2)<\dots< a(p) =\ell\quad\hbox{and}
\leqno{(8.5.5)}$$
$$\{k\in\{j,j+1,\dots,\ell-1\}: W^{(M)}_k(\omega) =1\} =
\{a(0),a(1),a(2),\dots,a(p-1)\}. \leqno{(8.5.6)}$$

Recall that since $\omega\in \Omega_0$ (by hypothesis), the sequence
$(W^{(M+1)}_k(\omega)$, $k\in\Z)$ of elements of $\{0,1,\dots,6\}$
satisfies Condition $\cS$ (see statement 6.3(A)(i) and Definition
2.3(A) again). Hence by (8.5.3), there exist integers
$i(1),i(2),\dots,i(6)$ with $j=i(1)<i(2)<\dots< i(6)<\ell$ such that
$$\leqalignno{&W^{(M+1)}_{i(u)}(\omega) =u \quad\hbox{for}\quad
u\in\{1,2,\dots,6\} \quad\hbox{and} &(8.5.7) \cr
&W^{(M+1)}_k (\omega) =0 \quad \hbox{for all other $k\in
\{j,j+1,\dots,\ell-1\}$.} \cr}$$
By (8.5.7) and Remark 6.3(B), $W^{(M)}_{i(u)}(\omega) =1$ for each
$u\in \{1,2,\dots,6\}$. Referring to (8.5.5) and (8.5.6), for each
$u\in \{1,2,\dots,6\}$, let $e(u)$ denote the element of
$\{0,1,\dots,p-1\}$ such that
$$i(u) = a(e(u)). \leqno{(8.5.8)}$$
Then by (8.5.5), (8.5.8), and the equality $j=i(1)$ just before
(8.5.7), one has that $a(0) =j = i(1) = a(e(1))$; and hence (see
(8.5.5) and the inequalities right before (8.5.7) again)
$$0= e(1)<e(2)<\dots< e(6) \le p-1. \leqno{(8.5.9)}$$

For each $u\in \{1,2,\dots,6\}$, referring to the inequality
$a(e(u))<a(e(u)+1)$ from (8.5.5) (see (8.5.9)), define the set
$$E_u:= \{k\in{\Z}:\ a(e(u)) \le k <a(e(u)+1) \ \ \hbox{and}
\ \  
N_k(\omega) \ge M\}. \leqno{(8.5.10)}$$
For each $u\in \{1,2,\dots,6\}$, from (8.5.5) and (8.5.6)
(and the sentence right after (8.5.4)), one has
that $W^{(M)}_{a(e(u))}(\omega) = W^{(M)}_{a(e(u)+1)}(\omega) =1$ and
$W^{(M)}_k(\omega)\not= 1$ for all $k\in\Z$ such that
$a(e(u))<k<a(e(u)+1)$. Hence by (8.5.10) and Definition 8.1 with
$m=M$,
$$E_u\in \cE_M(\omega) \quad \hbox{for each}\quad u\in
\{1,2,\dots,6\}. \leqno{(8.5.11)}$$
Hence by the induction assumption of Lemma 8.5(A)(B) for the case
$m=M$, $$\hbox{card }\, E_u =6^M \quad \hbox{for each }\quad
u\in\{1,2,\dots,6\}. \leqno{(8.5.12)}$$

For each $t\in \{0,1,\dots,p-1\}$ (see (8.5.5) and (8.5.6)), one has
the following: For each $k\in\Z$ such that $a(t) \le k<a(t+1)$,
$\Psi(M,k,0)(\omega) = k-a(t)$ by (8.5.5), (8.5.6), (6.2.6), and
Definition 2.2, hence $a(t) = k-\Psi(M,k,0)(\omega)$, and hence by
(7.1.2),
$$\delta^{(M+1)}_k(\omega) = W^{(M+1)}_{k-\Psi(M,k,0)(\omega)}(\omega) =
W^{(M+1)}_{a(t)}(\omega). \leqno{(8.5.13)}$$

For each $u\in \{1,2,\dots,6\}$, $W^{(M+1)}_{a(e(u))}(\omega) =u$ by
(8.5.7) and (8.5.8), and hence by (8.5.13) (and its entire sentence)
one has that
$$\leqalignno{\qquad\qquad\delta^{(M+1)}_k(\omega) &= W^{(M+1)}_{a(e(u))}(\omega) = u \not= 0 &(8.5.14)\cr
&\quad\hbox{for all $k\in\Z$ such that $a(e(u)) \le k< a(e(u)+1)$.}
\cr }$$ 
For each $u\in \{1,2,\dots,6\}$ and each $k\in E_u$,
$N_k(\omega) \ge M$ by (8.5.10), and hence $N_k(\omega) \ge M+1$ by
(8.5.14) and (7.1.4)--(7.1.5). Hence by (8.5.10), for each $u\in
\{1,2,\dots,6\}$,
$$\quad E_u = \left\{k\in {\Z}:\ a(e(u)) \le k< a(e(u)+1) \ \ \hbox{and}\ \ 
N_k(\omega) \ge M+1\right\}. \leqno{(8.5.15)}$$

For any $t\in \{0,1,\dots,p-1\} - \{e(1),e(2),\dots,e(6)\}$ (see
(8.5.8) and its entire sentence), one has that $j\le a(t) <\ell$ by
(8.5.5), $a(t) \not= a(e(u)) = i(u)$ for each
$u \in \{1,2,\dots,6\}$ 
by (8.5.5) and (8.5.8), and
hence $W^{(M+1)}_{a(t)}(\omega) =0$ by (8.5.7). Hence for each $t \in
\{0,1,\dots, p-1\} - \{e(1),e(2),\dots,e(6)\}$ and each $k\in\Z$
such that $a(t) \le k<a(t+1)$, one has that $\delta^{(M+1)}_k(\omega)
=0$ by (8.5.13) and its entire sentence, and hence $N_k(\omega) \le
M$ by (7.1.4). Hence by (8.5.4), (8.5.5), and (8.5.15), the set $E$
has no elements other than the ones in the sets $E_u$, $u\in
\{1,2,\dots,6\}$. In fact by (8.5.4) and (8.5.15), one now has that
$$E = E_1\cup E_2 \cup \dots \cup E_6, \leqno{(8.5.16)}$$
and by (8.5.5), (8.5.9), and (8.5.15) (see also (8.5.12)),
$$E_1<E_2<\dots<E_6 \leqno{(8.5.17)}$$
(see section 2.1(D)). Hence by (8.5.12),
$$\hbox{card}\, E=6^{M+1}. \leqno{(8.5.18)}$$

Thus statement (A) in Lemma 8.5 holds (for the given $\omega$ and $E$)
with $m=M+1$. By (8.5.11), (8.5.16), and (8.5.17), statement (C) in
Lemma 8.5 holds (for the given $\omega$ and $E$) with $m=M+1$. To
complete the induction step and the proof of Lemma 8.5, our task now
is to verify statement (B) (for the given $\omega$ and $E$) with
$m=M+1$.

Now refer to (8.5.11) and our induction assumption of statements (A)
and (B) for the case $m=M$. For each $u\in \{1,2,\dots,6\}$,
referring to (8.5.12), representing the set $E_u$ by
$$\leqalignno{E_u &= \left\{\alpha(u,1),\alpha(u,2),
\alpha(u,3),\dots,\alpha(u,6^M)\right\} &(8.5.19)\cr
&\qquad\hbox{with
$\alpha(u,1)<\alpha(u,2)<\dots<\alpha(u,6^M)$,} \cr }$$
one has that
$$J(M,\alpha(u,v))(\omega) = v-1 \quad \hbox{for all $v\in
\{1,2,\dots,6^M\}$.} \leqno{(8.5.20)}$$

Referring to (8.5.18), represent the set $E$ by
$$E = \left\{ \beta(1),\beta(2),\beta(3),\dots,\beta(6^{M+1})\right\}
\,\,\hbox{with}\quad \beta(1)<\beta(2)<\dots <\beta(6^{M+1}).
\leqno{(8.5.21)}$$ Then for each $u\in \{1,2,\dots,6\}$, by
(8.5.12), (8.5.16), (8.5.17), and (8.5.21), the set $E_u$ contains
precisely the elements $\beta(6^M(u-1)+v)$, $v\in
\{1,2,\dots,6^M\}$. Hence by (8.2.19),
$$\forall\,\, u\in \{1,2,\dots,6\}, \quad \forall\,\, v\in
\{1,2,\dots,6^M\}, \quad \beta\left(6^M(u-1)+v\right) = \alpha(u,v).
\leqno{(8.5.22)}$$

Now suppose $r\in \{1,2,\dots,6^{M+1}\}$. Let $u\in \{1,2,\dots,6\}$
and $v\in \{1,2,\dots,6^M\}$ be such that
$$r = 6^M(u-1)+v. \leqno{(8.5.23)}$$
Then $\beta(r) \in E_u$ by the sentence preceding (8.5.22), hence
$a(e(u)) \le \beta(r) <a(e(u)+1)$ by (8.5.10), and hence
$\delta^{(M+1)}_{\beta(r)}(\omega) = W^{(M+1)}_{a(e(u))} (\omega) =u$
by (8.5.14). Also, $\beta(r)  = \alpha(u,v)$ by (8.5.23) and
(8.5.22). Hence by (7.1.7) (applied twice) and then (8.5.20) and
(8.5.23),
$$\eqalign{
J(M+1,&\beta(r))(\omega) =
\sum^{M+1}_{q=1}6^{q-1}\left(\delta^{(q)}_{\beta(r)} (\omega)
-1\right) \cr
&= 6^M \cdot \left(\delta^{(M+1)}_{\beta(r)}(\omega) -1\right) +
\sum^M_{q=1} 6^{q-1}\left(\delta^{(q)}_{\beta(r)}(\omega) -1\right)
\cr
&= 6^M(u-1) + \sum^M_{q=1} 6^{q-1}
\left(\delta^{(q)}_{\alpha(u,v)}(\omega) -1\right) \cr
&= 6^M (u-1) +J(M,\alpha(u,v))(\omega) \cr
&= 6^M(u-1) + v-1 \cr
&= r-1. \cr}$$

Since $r\in \{1,2,\dots,6^{M+1}\}$ was arbitrary, one has (see
(8.5.21) again) that statement (B) in Lemma 8.5 holds (for the given
$\omega$ and $E$) with $m=M+1$. That completes the induction step and the proof of Lemma 8.5.
\medskip

{\sc Definition 8.6.} Recall from Definitions 8.1 and 8.3 and Lemma
8.5(A) that for each $\omega\in \Omega_0$ (see (6.3.1)), one has that
(i)~card $E=1$ for each $E\in \cD_0(\omega)$, and (ii)~for each
$m\in\N$ and each $E\in \cE_m(\omega)$ (and in particular, for each
$E\in \cD_m(\omega))$, card~$E=6^m$. Recall also Remarks 8.2(D)(E)
and 8.4(B)(C).

(A) For each $\ell\in \{0,1,2,\dots\,\}$ and each set $D\subset \Z$
such that card~$D=6^\ell$, define the set $F_D \subset \Omega_0$ (see
(6.3.1)) as follows:
$$F_D := \{\omega\in \Omega_0:D\in \cD_\ell(\omega)\}.
\leqno{(8.6.1)}$$

(B) For each $m\in\N$ and each set $E\subset \Z$ such that
card~$E=6^m$, define the set $G_E \subset \Omega_0$ as follows:
$$G_E := \{\omega\in \Omega_0:E\in \cE_m(\omega)\}. \leqno{(8.6.2)}$$

{\sc Remark 8.7.} Refer to Lemma 8.5. Suppose $m\in\N$, and $E\subset
\Z$ is a set such that card~$E=6^m$.

(A) The set $G_E$ in (8.6.2) is the set of all $\omega\in \Omega_0$
for which there exist integers $j$ and $\ell$ such that (8.1.1) and
(8.1.2) hold. It follows that $G_E$ is an event (that is, a member of
the $\sigma$-field $\cF$ in our given probability space
$(\Omega,\cF,P))$; and further, by (7.1.5),
$$G_E \in \sigma(\, \overline{W}^{(m)}). \leqno{(8.7.1)}$$
By the same argument, but with the inequality $N_k(\omega) \ge m$ in
(8.1.2) replaced by $N_k(\omega) =m$ (see Remark 8.2(E) and
Definition 8.3(B)), the set $F_E$ in (8.6.1) is an event, and by
(7.1.6), $$F_E \in \sigma(\, \overline{W}^{(m+1)}).
\leqno{(8.7.2)}$$

(B) By Definitions 8.1, 8.3(B), and 8.6 (recall Remark 8.2(E) again),
$$F_E\subset G_E; \leqno{(8.7.3)}$$
$$\forall\,\,\omega\in F_E,\,\,\forall\,\, k\in E,\quad
N_k(\omega)=m;\quad\hbox{and} \leqno{(8.7.4)}$$
$$\forall\,\, \omega\in G_E - F_E, \,\,\forall\,\, k\in
E,\quad N_k(\omega)\ge m+1. \leqno{(8.7.5)}$$

{\sc Remark 8.8.} Equations (8.7.2) and (8.7.4) hold with $m=0$ for
any singleton set $E=\{k\}$ where $k\in \Z$, by (8.6.1), Definition
8.3(A), and (7.1.6).
\medskip

{\sc Definition 8.9.} Refer to Constructions 7.3 and 7.4, to
equations (6.2.6), (7.1.4), and (7.1.7), and also particularly to
Remark 7.2(B).

For each $n\in\N$, define the sequences $Y^{(n)}:= (Y^{(n)}_k$, $k\in
\Z)$ and $X^{(n)}:= (X^{(n)}_k$, $k\in \Z)$ of $\{-1,1\}$--valued
random variables as follows: For each $k\in\Z$ and each
$\omega\in\Omega$,
$$Y^{(n)}_k(\omega) := \cases{ X^{(0)}_k(\omega) &if $N_k(\omega) =0$
\cr
\zeta^{(\ell,\cen)}_{k-\Psi(\ell,k,0)(\omega),J(\ell,k)(\omega)}(\omega)
&if $N_k(\omega) = \ell\in \{1,2,\dots,n-1\}$ \cr
\zeta^{(n,\ord)}_{k-\Psi(n,k,0)(\omega),J(n,k)(\omega)}(\omega) &
if $N_k(\omega) \ge n$, \cr} \leqno{(8.9.1)}$$ and
$$X^{(n)}_k(\omega) := \cases{X^{(0)}_k(\omega) &if $N_k(\omega) =0$
\cr
\zeta^{(\ell,\cen)}_{k-\Psi(\ell,k,0)(\omega),J(\ell,k)(\omega)}(\omega)
&if $N_k(\omega) =\ell\in \{1,2,\dots,n\}$ \cr
\zeta^{(n,\fri)}_{k-\Psi(n,k,0)(\omega),J(n,k)(\omega)}(\omega)
&if $N_k(\omega) \ge n+1$. \cr} \leqno{(8.9.2)}$$ Of course in the
right hand side of (8.9.1), the ``middle'' case is vacuous  (and
should be omitted) if $n=1$.
\medskip

{\sc Remark 8.10.} If $k\in\Z$ and $D= \{k\}$, then by Remark 8.8 (see
 equation (8.7.4)) and (8.9.1) and (8.9.2), for every $\omega\in F_D$,
one has that
$$Y^{(n)}_k(\omega) = X^{(n)}_k(\omega) = X^{(0)}_k(\omega)
\quad\hbox{for all $n\in\N$.} \leqno{(8.10.1)}$$

{\sc Lemma 8.11.} {\it Suppose $m\in\N$, and $E\subset \Z$ is a set such
that
$$\card\, E = 6^m. \leqno{(8.11.1)}$$
Then in the terminology of (2.3) and (2.4), the following statements
hold:

(A) For every $\omega\in G_E$ (see (8.6.2)),
$$Y^{(m)}_E(\omega) = \zeta^{(m,\ord)}_{\min\, E}(\omega)
\leqno{(8.11.2)}$$ and
$$ \leqalignno {
X^{(m)}_E (\omega) = &\zeta^{(m,\cen)}_{\min\, E}(\omega)\cdot
I\left(W^{(m+1)}_{\min\, E} =0\right)(\omega) & (8.11.3) \cr
& \indent +\zeta^{(m,\fri)}_{\min\, E}(\omega) \cdot I\left(W^{(m+1)}_{\min\,
E} \not= 0\right)(\omega). \cr}$$

(B) For every $\omega\in F_E$ (see (8.6.1) and (8.7.3)),}
$$Y^{(n)}_E(\omega) = X^{(n)}_E(\omega) = X^{(m)}_E(\omega) =
\zeta^{(m,\cen)}_{\min\, E}(\omega)\quad \hbox{\it for all $n\ge m+1$.}
\leqno{(8.11.4)}$$
\medskip

{\sc Proof.} We shall prove statements (A) and (B) together.

Suppose $\omega\in G_E$.

Then by (8.11.1) and Definition 8.6(B), $E\in \cE_m(\omega)$.
Referring again to (8.11.1), represent the set $E$ by
$$E = \left\{e(1),e(2),e(3),\dots,e(6^m)\right\}
\quad\hbox{where}\quad e(1)<e(2)<\dots< e(6^m). \leqno{(8.11.5)}$$
Also, referring to Definition 8.1 and Remark 8.2(B), let $j$ and
$\ell$ be the integers such that (8.1.1) and (8.1.2) hold. We shall
refer freely to (8.1.1) and (8.1.2) in the arguments below. By Remark
8.2(A), (8.1.2), and (8.11.5),
$$j = \min \, E = e(1)\quad\hbox{and}\quad \ell>e(6^m).
\leqno{(8.11.6)}$$
For each $k\in E$, by (8.1.1), (8.1.2), (6.2.6), and Definition 2.2,
$\Psi(m,k,0)(\omega) = k-j$. Hence by (8.11.5),
$$\forall\,\, u\in \{1,2,\dots,6^m\}, \quad j = e(u) -
\Psi(m,e(u),0)(\omega). \leqno{(8.11.7)}$$
Also, by (8.11.5) and (8.1.2),
$$\forall\,\, u\in \{1,2,\dots,6^m\}, \quad N_{e(u)}(\omega) \ge m.
\leqno{(8.11.8)}$$

{\it Proof of (8.11.2).} For each $u\in \{1,2,\dots,6^m\}$, by
(8.11.8) and (8.9.1), followed by (8.11.7), (8.11.5), and Lemma
8.5(B) (recall also section 7.3(B)),
$$Y^{(m)}_{e(u)}(\omega) = \zeta^{(m,\ord)}_{e(u) -
\Psi(m,e(u),0)(\omega),J(m,e(u))(\omega)}(\omega) =
\zeta^{(m,\ord)}_{j,u-1}(\omega). \leqno{(8.11.9)}$$
Hence by (8.11.5) and (8.11.6) (see (2.3), (2.4), and section
7.3(B)), one obtains (8.11.2) via
$$Y^{(m)}_E(\omega) = \zeta^{(m,\ord)}_j (\omega) =
\zeta^{(m,\ord)}_{\min\, E}(\omega). \leqno{(8.11.10)}$$

{\it Proof of (8.11.3).} If $W^{(m+1)}_{\min\, E}(\omega) =0$, then
$N_{e(u)}(\omega) = m$ for all $u\in \{1,2,\dots,6^m\}$ by (8.11.5)
and Remark 8.2(E)(i), and following the argument for (8.11.2) (but
using (8.9.2) instead of (8.9.1)), one obtains (8.11.9) and (8.11.10)
with the letters ``$Y$'' and ``ord'' replaced by ``$X$'' and ``cen,''
and in this case (8.11.3) holds. If instead $W^{(m+1)}_{\min\, E}
(\omega) \not= 0$, then $N_{e(u)}(\omega) \ge m+1$ for all $u\in
\{1,2,\dots,6^m\}$ by (8.11.5) and Remark 8.2(E)(ii), and with the
same argument, one obtains (8.11.9) and (8.11.10) with  ``$Y$'' and
``ord'' replaced by ``$X$'' and ``$\fri$,'' and in this case too
equation (8.11.3) holds.

{\it Proof of (8.11.4).} Suppose $\omega\in F_E$ (assumed for
(8.11.4)). Then $E\in \cD_m(\omega)$ by Definition 8.6(A), and hence
$N_{e(u)}(\omega) =m$ for all $u\in \{1,2,\dots,6^m\}$ by (8.11.5)
and Definition 8.3(B). Hence by (8.9.1) and (8.9.2) (applied twice),
for each $n\ge m+1$ and each $u\in \{1,2,\dots,6^m\}$,
$$Y^{(n)}_{e(u)}(\omega) = X^{(n)}_{e(u)}(\omega) =
\zeta^{(m,\cen)}_{e(u) -
\Psi(m,e(u),0)(\omega),J(m,e(u))(\omega)}(\omega) =
X^{(m)}_{e(u)}(\omega).$$
Hence (see (8.11.5) again) for each $n\ge m+1$,
$Y^{(n)}_E(\omega)=X^{(n)}_E(\omega) = X^{(m)}_E(\omega)$. Since
$N_{e(u)}(\omega) =m$ for all $u\in \{1,2,\dots,6^m\}$ (as was noted
above),
 $W^{(m+1)}_{\min \, E}(\omega) =0$ must hold by Remark
8.2(E), and the final equality in (8.11.4) now follows from (8.11.3).
That completes the proof of Lemma 8.11.
\medskip

{\sc Definition 8.12.} Suppose $n\in \N$. Suppose $\cQ :=
\{Q(1),Q(2),\dots, \allowbreak Q(L)\}$ (where $L$ is a positive integer) is a
(finite, nonempty) collection of (pairwise) disjoint subsets of $\Z$
such that
$$\card\, Q(\ell) \in \{1,6,6^2,6^3,\dots,6^n\} \quad \hbox{for all
$\ell\in \{1,2,\dots,L\}$}. \leqno{(8.12.1)}$$ (There is no
assumption of an ``ordering'' of these sets; elements of one set $Q(i)$ may be between elements of another set $Q(j)$.
Also, it is tacitly understood that $L$ can be 1, in which case the phrase
``(pairwise) disjoint'' is meaningless and should be 
omitted.)\ \ Define the
(possibly empty) set
$$\cI := \cI (n,\cQ) := \{\ell\in \{1,\dots,L\}: \card\, Q(\ell) =
6^n\}. \leqno{(8.12.2)}$$
Define the event $H^{(n)}(\cQ)$ by (see (8.6.1)--(8.6.2))
$$H^{(n)}(\cQ):= \left[ \bigcap_{\ell\in \{1,\dots,L\} -\cI}
F_{Q(\ell)}\right]\bigcap \left[\bigcap_{\ell\in \cI}
G_{Q(\ell)}\right]. \leqno{(8.12.3)}$$ Here (if necessary), define
the ``vacuous intersection'' by $\bigcap_{\ell\in \emptyset}(\dots )
:= \Omega$ (the sample space itself).

It is easy to see that this definition of the event $H^{(n)}(\cQ)$
does not depend on the particular order in which the sets in $\cQ$
are labeled (as $Q(\ell)$, $\ell\in \{1,\dots,L\}$).
\medskip

{\sc Lemma 8.13.} {\it Suppose $n\in \N$. In the context of
Definition 8.12, with (8.12.1) satisfied and with $\cI$ and
$H^{(n)}(\cQ)$ defined as in (8.12.2) and (8.12.3), suppose
$$P\left(H^{(n)}(\cQ)\right)>0. \leqno{(8.13.1)}$$
For each $\ell\in \{1,2,\dots,L\}$, define the integer
$$q(\ell) := \min\, Q(\ell). \leqno{(8.13.2)}$$
Then the following statements hold:

(A) One has that $H^{(n)}(\cQ) \in \sigma(\overline{W}^{(n)})$. Also,
for each\ \thinspace $\omega\in H^{(n)}(\cQ)$ and each 
$\ell \in \cI$ (if $\cI$
is nonempty), $W^{(n)}_{q(\ell)}(\omega) =1$.

(B) For each $\omega\in H^{(n)}(\cQ)$ and each $\ell\in
\{1,2,\dots,L \}$ such that $\card\, Q(\ell)=1$ (and hence $Q(\ell) =
\{q(\ell)\}$), one has that
$$X^{(n-1)}_{q(\ell)}(\omega) = Y^{(n)}_{q(\ell)}(\omega) =
X^{(n)}_{q(\ell)}(\omega) = X^{(0)}_{q(\ell)}(\omega).
\leqno{(8.13.3)}$$

(C) For each $\omega\in H^{(n)}(\cQ)$, each $m\in \{1,\dots,n-1\}$
(if $n\ge 2$), and each $\ell\in \{1,\dots,L\}$ such that $\card\,
Q(\ell) =6^m$, one has that (see (2.3) and (2.4))
$$X^{(n-1)}_{Q(\ell)} (\omega) = Y^{(n)}_{Q(\ell)}(\omega) =
X^{(n)}_{Q(\ell)}(\omega) = \zeta^{(m,\cen)}_{q(\ell)}(\omega).
\leqno{(8.13.4)}$$

(D) For any given $\omega\in H^{(n)}(\cQ)$ and any given $\ell\in
\cI$ (see (8.12.2)), representing the set $Q(\ell)$ by $Q(\ell) :=
\{z(1),z(2),\dots,z(6^n)\}$ where $z(1)<z(2)<\dots <z(6^n)$ (and
hence $z(1) = q(\ell))$, and denoting $y(u):= z((u-1)\cdot
6^{n-1}+1)$ for $u\in \{1,2,\dots,6\}$ (and hence $y(1) =z(1)
=q(\ell))$, one has (see section 2.1(C)) that
$$\leqalignno{\qquad X^{(n-1)}_{Q(\ell)} (\omega)
 = \Bigl\langle
&\zeta^{(n-1,\fri)}_{q(\ell)}(\omega),\zeta^{(n-1,\fri)}_{y(2)}(\omega),
\zeta^{(n-1,\fri)}_{y(3)}(\omega), &(8.13.5)\cr &\quad
\zeta^{(n-1,\fri)}_{y(4)}(\omega),
\zeta^{(n-1,\fri)}_{y(5)}(\omega),
\zeta^{(n-1,\fri)}_{y(6)}(\omega)\Bigr\rangle \cr}$$ (with
$\zeta^{(0,\fri)}_k(\omega) := X^{(0)}_k(\omega)$ for all $k\in \Z$
in the case $n=1$),
$$Y^{(n)}_{Q(\ell)}(\omega) = \zeta^{(n,\ord)}_{q(\ell)}(\omega)
\leqno{(8.13.6)}$$
and
$$X^{(n)}_{Q(\ell)}(\omega) = \zeta^{(n,\cen)}_{q(\ell)}(\omega)
\cdot I\left(W^{(n+1)}_{q(\ell)} =0\right)(\omega) +\zeta^{(n,\fri)}_{q(\ell)}
(\omega) \cdot I\left(W^{(n+1)}_{q(\ell)} \not=0 \right)(\omega).
\leqno{(8.13.7)}$$

(E) Refer to (8.13.1) and (8.12.2). The random variables
$I(W^{(n+1)}_{q(\ell)}\not= 0)$, $\ell\in \cI$, the random variables
$X^{(0)}_k$, $k\in\Z$, and the random variables $\zeta^{(n,\ord)}_k$,
$\zeta^{(n,\cen)}_k$, $\zeta^{(n,\fri)}_k$, $n\in\N$, $k\in\Z$ are
conditionally independent given the event $H^{(n)}(\cQ)$. Also, for
each $\ell\in\cI$, conditional on $H^{(n)}(\cQ)$, the random
variable $I(W^{(n+1)}_{q(\ell)}\not=0)$ takes the value 0 resp.\ 1
with probability 5/8 resp.\ 3/8.

(F) The random vectors $X^{(n-1)}_{Q(\ell)}$, $\ell\in \{1,\dots,L\}$
are conditionally independent given~$H^{(n)}(Q)$.

(G) The random vectors $Y^{(n)}_{Q(\ell)}$, $\ell\in
\{1,\dots,L\}$ are conditionally independent given $H^{(n)}(\cQ)$.

(H) The random vectors $X^{(n)}_{Q(\ell)}$, $\ell\in
\{1,\dots,L\}$ are conditionally independent given $H^{(n)}(\cQ)$.

(I) Refer to (8.12.2), Section 2.5(A)(B), and Definition 3.3.  For each $\ell\in \cI$, one has that
(i) $\cL\left(X^{(n-1)}_{Q(\ell)}\, \Bigl| \, 
\cH^{(n)}(\cQ)\right) = \left(\nu^{(n-1)}_{\fri}\right)^{[6]}$ and 

\noindent (ii) $\cL\left(Y^{(n)}_{Q(\ell)}\, \Bigl| \,  
H^{(n)}(\cQ)\right)
= \cL\left(X^{(n)}_{Q(\ell)}\, \Bigl| \,  
H^{(n)}(\cQ)\right) =
\nu^{(n)}_{\ord}$.

(J) Recall from Section 2.1(H) the notation $\union(\cQ):= Q(1) \cup
Q(2) \cup \dots \cup Q(L)$. One has that
$\cL\left(Y^{(n)}_{\union(\cQ)}\, \Bigl| \,  
H^{(n)}(\cQ)\right) =
\cL\left(X^{(n)}_{\union (\cQ)}\, \Bigl| \,  
H^{(n)}(\cQ)\right)$.

(K) Suppose $S$ is a nonempty finite subset of $\Z$ such that $S\cap
\union (\cQ)$ is nonempty and $\card(S\cap Q(\ell)) \le 5$ for all
$\ell\in \{1,\dots,L\}$. Then
$$\cL\left(X^{(n-1)}_{S\,\cap\, \union
(\cQ)} \, \Bigl| \,  
H^{(n)}(\cQ)\right) = \cL\left(Y^{(n)}_{S\,\cap \,\union
(\cQ)}\, \Bigl| \,  
H^{(n)}(\cQ)\right).$$

(L) For any nonempty set $S\subset Q(1)\cup Q(2) \cup \dots \cup
Q(L)$, one has that
$$E\left[ \left(\sum_{k\in S}X_k^{(n-1)}\right)^6 \, \Biggl| \,
H^{(n)}(\cQ)\right] \ge E\left[ \left(\sum_{k\in
S}Y^{(n)}_k\right)^6 \, \Biggl| \, H^{(n)}(\cQ)\right].
\leqno{(8.13.8)}$$

(M) Suppose $\ell\in \cI$ (see (8.12.2)) and $S$ is a set such that
$Q(\ell) \subset S\subset Q(1)\cup Q(2) \cup \dots \cup Q(L)$. Then}
$$\leqalignno{&E\left[\left(\sum_{k\in S} X^{(n-1)}_k\right)^6 \,\Bigg| \,H^{(n)}(\cQ)\right] &(8.13.9) \cr
& \indent \ge E\left[\left(\sum_{k\in S} Y^{(n)}_k\right)^6 \,\Bigg| \,
H^{(n)}(\cQ)\right] +720 \cdot 4^{6(n-1)}. \cr}$$
\medskip

{\sc Proof.}  Let us first prove statement (A). For each $\ell\in
\{1,\dots,L\} -\cI$, letting $m(\ell) \in \{0,1,\dots,n-1\}$ denote
the integer such that card~$Q(\ell) = 6^{m(\ell)}$ (see (8.12.1) and
(8.12.2)), one has that $F_{Q(\ell)} \in
\sigma(\overline{W}^{(m(\ell)+1)}) \subset
\sigma(\overline{W}^{(n)})$ by (8.7.2) (see section 6.3(E) again).
For each $\ell\in \cI$, $G_{Q(\ell)}\in \sigma(\overline{W}^{(n)})$
by (8.7.1) (and (8.12.2)). Hence $H^{(n)}(\cQ) \in
\sigma(\overline{W}^{(n)})$ by (8.12.3). Also, for each $\omega\in
H^{(n)}(\cQ)$ and each $\ell\in \cI$, one has that $\omega\in
G_{Q(\ell)}$ by (8.12.3), hence $Q(\ell) \in \cE_n(\omega)$ by
(8.12.2) and Definition 8.6(B), hence $W^{(n)}_{q(\ell)}(\omega) =1$
by (8.13.2), (8.2.1), and (8.1.1). Both parts of statement (A) have
been proved.

{\it Proof of (B).} For each $\omega\in H^{(n)}(\cQ)$ and each
$\ell\in \{1,\dots,L\}$ such that $\card\, Q(\ell) =1$, one has that
$\omega\in F_{Q(\ell)}$ by (8.12.2) and (8.12.3), and hence
(8.13.3) holds by (8.13.2) and (8.10.1) (regardless of whether $n=1$
or $n\ge 2$).

{\it Proof of (C).} If $\omega\in H^{(n)}(\cQ)$, $m\in
\{1,\dots,n-1\}$, $\ell\in \{1,\dots,L\}$, and $\card\,  Q(\ell)=6^m$,
then $\omega\in F_{Q(\ell)}$ by (8.12.2) and (8.12.3), and hence
(8.13.4) holds by (8.13.2) and (8.11.4) (regardless of whether
$n-1=m$ or $n-1>m$).

{\it Proof of (D).} Suppose $\omega\in H^{(n)}(\cQ)$ and $\ell \in
\cI$. Then $\omega\in G_{Q(\ell)}$ by (8.12.3). Referring to (8.13.2),
one can apply (8.11.2) and (8.11.3) with $m=n$ (see (8.11.1) and
(8.12.2)), and thereby one obtains (8.13.6) and (8.13.7). The
remaining task is to prove (8.13.5).

If $n=1$, and hence $\card\, Q(\ell)=6$, then in the terminology of
statement (D), $Q(\ell) := \{z(1),z(2),\dots,z(6)\}$ where
$z(1)<z(2)<\dots<z(6)$, also $y(u) = z(u)$ for each $u\in
\{1,\dots,6\}$ (and $q(\ell) = y(1) = z(1))$, hence
$$X^{(n-1)}_{Q(\ell)} (\omega) =
\left(X^{(0)}_{q(\ell)}(\omega),X^{(0)}_{y(2)}(\omega),
X^{(0)}_{y(3)}(\omega),\dots,X^{(0)}_{y(6)}(\omega)\right),$$
and thus (8.13.5) holds (see the phrase right after (8.13.5)).

Now suppose instead that $n\ge 2$. Recall again from (8.12.3) that
$\omega\in G_{Q(\ell)}$ (since $\ell\in \cI$). By (8.6.2), $Q(\ell)
\in \cE_n(\omega)$. Hence by (8.1.2),
$$\forall\,\,k\in Q(\ell),\qquad N_k(\omega) \ge n.
\leqno{(8.13.10)}$$
Applying Lemma 8.5(C) (with $m=n$), let $E(1), E(2), \dots,
E(6)$ be sets $\subset \Z$ such that $E(1) \cup E(2)\cup \dots \cup
E(6) = Q(\ell)$, $E(1)<E(2)<\dots<E(6)$, and $E(u) \in
\cE_{n-1}(\omega)$ for each $u\in \{1,\dots,6\}$.  For each $u\in
\{1,\dots,6\}$, card $E(u) = 6^{n-1}$ by Lemma 8.5(A) (applied with
$m=n-1$). It follows that in the terminology of statement (D) (here
in Lemma 8.13), for each $u\in \{1,\dots,6\}$,
$$E(u) = \left\{z((u-1) \cdot 6^{n-1}+1), z((u-1)\cdot 6^{n-1}
+2),\dots,z(u\cdot 6^{n-1})\right\},$$
and hence also $\min E(u) = y(u)$ (recall also that $y(1) =
q(\ell))$. Also (see section 2.1(C))
$$X^{(n-1)}_{Q(\ell)}(\omega) = \left\langle
X^{(n-1)}_{E(1)}(\omega) ,X^{(n-1)}_{E(2)}(\omega)
,\dots,X^{(n-1)}_{E(6)}(\omega)\right\rangle.$$
For each $u\in \{1,\dots,6\}$, since $E(u) \in \cE_{n-1}(\omega)$ (as
was noted above), one has by (8.13.10) and Remark 8.2(E) that
$W^{(n)}_{y(u)}(\omega) \not= 0$ must hold, and hence by (8.11.3),
$X^{(n-1)}_{E(u)}(\omega) = \zeta^{(n-1,\fri)}_{y(u)}(\omega)$.
Equation (8.13.5) follows. That completes the proof of statement (D).

{\it Proof of (E).} By statement (A), $H^{(n)}(\cQ) = H^{(n)}(\cQ)
\bigcap [\,\bigcap_{\ell\in \cI}\{W^{(n)}_{q(\ell)} =1\}]$. Hence by
((8.13.1) and)  Lemma 6.4(D), conditional on the event
$H^{(n)}(\cQ)$, the $\{0,1\}$--valued random variables
$I(W^{(n+1)}_{q(\ell)}\not= 0)$, $\ell\in \cI$ are independent, with
each taking the value 0 resp.\ 1 with probability 5/8 resp.\ 3.8.
Also, the random sequences $\overline{\xi}^{(\infty)}$,
$\overline{\zeta}^{(\infty)}$, and $X^{(0)}$ are independent by
Constructions 7.3 and 7.4; and hence the random sequences
$\overline{W}^{(n+1)}$, $\overline{\zeta}^{(\infty)}$, and $X^{(0)}$
are independent by (6.2.1) (and Definition 2.6). Also, the random
variables $X^{(0)}_k$, $k\in \Z$ and $\zeta^{(n,\ord)}_k$,
$\zeta^{(n,\cen)}_k$, and $\zeta^{(n,\fri)}_k$, $n\in\N$, $k\in \Z$
are all independent of each other, by Constructions 7.3(A) and
7.4(A). Putting these pieces (and the first sentence of statement
(A)) together, one obtains statement (E) by an awkward but trivial
calculation involving conditional probabilities.

{\it Proof of (F), (G), and (H).} We shall just give the argument for
(H). The arguments for (F) and (G) are similar, using (8.13.5) resp.\
(8.13.6) in the places where (in the argument below for (H)) (8.13.7) is used.  (Keep in mind that here in our context
of Definition 8.1.2, the sets $Q(\ell)$, 
$\ell \in \{1,2,\dots,L\}$ are (pairwise) disjoint
(if $L \geq 2$).)

Define the random variables (vectors) $Z_\ell$, $\ell\in
\{1,\dots,L\}$ as follows: For $\omega\in \Omega_0$,
$$Z_\ell(\omega) := \cases{X^{(0)}_{q(\ell)}(\omega) &if card
$Q(\ell) =1$ \cr
\zeta^{(m,\cen)}_{q(\ell)} (\omega) &if card $Q(\ell)=6^m$ where
$m\in \{1,\dots,n-1\}$ \cr
\hbox{[RHS of (8.13.7)]} &if $\ell\in \cI$ \cr} \leqno{(8.13.11)}$$
(Of course in the right hand side of (8.13.11), the ``middle'' part
is vacuous, and should be omitted, if $n=1$.)

By statement (E) (and the trivial fact that
$I(W^{(n+1)}_{q(\ell)}=0) = 1-I(W^{(n+1)}_{q(\ell)} \not= 0)$ for
$\ell\in \cI$), the random variables $Z_1,Z_2,\dots,Z_L$ are
conditionally independent given $H^{(n)}(\cQ)$. Also, by (8.13.3),
(8.13.4), (8.13.7), and (8.13.11), $X^{(n)}_{Q(\ell)}(\omega) =
Z_\ell(\omega)$ for each $\omega\in H^{(n)}(\cQ)$ and each $\ell\in
\{1,\dots,L\}$. Hence (see section 2.5(A)) 
$$\eqalign{&\cL\left(\left(X^{(n)}_{Q(1)},X^{(n)}_{Q(2)},\dots,X^{(n)}_{Q(L)}\right)\, \Big| \,
H^{(n)}(\cQ)\right) = \cL\left( (Z_1,Z_2,\dots,Z_L) \, \Big| \,
H^{(n)}(\cQ)\right) \cr &=\cL\left(Z_1\, \Big| \,
H^{(n)}(\cQ)\right)\times \dots \times \cL\left(Z_L\, \Big| \,
H^{(n)}(\cQ)\right)\cr &= \cL\left(X^{(n)}_{Q(1)}\, \Big| \,
H^{(n)}(\cQ)\right)\times \dots \times \cL\left(X^{(n)}_{Q(L)}
\,\Big| \, H^{(n)}(\cQ)\right). \cr}$$ 
Thus (H) holds.

{\it Proof of (I).} Suppose $\ell\in \cI$. Recall the first sentence
in statement (A), and recall from the proof of (E) that the random
sequences $\overline{W}^{(n+1)}$, $\overline{\zeta}^{(\infty)}$, and
$X^{(0)}$ are independent of each other. By Construction 7.3(A) and
equation (8.13.5) (see also Section 2.5(A)(B)),
$$\eqalign{\cL\left(X^{(n-1)}_{Q(\ell)}\, \Big| \, H^{(n)}(\cQ)\right) &=
\cL\left([\hbox{RHS of (8.13.5)}]\, \Big| \, H^{(n)}(\cQ)\right) \cr
&= \cL([\hbox{RHS of (8.13.5)}]) =
\left(\nu^{(n-1)}_{\fri}\right)^{[6]}.\cr}$$ Thus statement (i) in
(I) holds. The equality $\cL(Y^{(n)}_{Q(\ell)}\mid H^{(n)}(\cQ)) =
\nu^{(n)}_{\ord}$ in statement (ii) holds by a similar argument
using (8.13.6) instead of (8.13.5). To obtain the final equality in
statement (ii), note that for any \break
\noindent $a\in
\{-1,1\}^{\sxtp(n)}$ (see (2.1)), by (8.13.7), statement (E), the
first sentence of (A), Construction 7.3(A), and Remark 3.4,
$$\eqalign{&P\left(X^{(n)}_{Q(\ell)} = a\, \Big| \, H^{(n)}(\cQ)\right) \cr
&=P\left(W^{(n+1)}_{q(\ell)} = 0 \, \Big| \, H^{(n)}(\cQ)\right)
\cdot P\left(X^{(n)}_{Q(\ell)} =a\, \Big| \, \left\{
W^{(n+1)}_{q(\ell)} =0 \right\} \cap H^{(n)}(\cQ)\right) \cr 
&\qquad
+P\left(W^{(n+1)}_{q(\ell)} \not= 0 \, \Big| \,
H^{(n)}(\cQ)\right)\cr 
&\qquad \qquad \cdot P\left(X^{(n)}_{Q(\ell)} =a\,
\Big| \, \left\{W^{(n+1)}_{q(\ell)} \not= 0\right\} \cap
H^{(n)}(\cQ)\right) \cr &= (5/8) \cdot
P\left(\zeta^{(n,\cen)}_{q(\ell)} = a\, \Big| \, \left\{
W^{(n+1)}_{q(\ell)} =0 \right\} \cap H^{(n)}(\cQ)\right) \cr &\qquad
+ (3/8) \cdot P\left(\zeta^{(n,\fri)}_{q(\ell)} =a \, \Big| \,
\left\{W^{(n+1)}_{q(\ell)} \not= 0 \right\}\cap H^{(n)}(\cQ)\right)
\cr &= (5/8) \cdot P\left(\zeta^{(n,\cen)}_{q(\ell)} =a\right) +
(3/8) \cdot P\left(\zeta^{(n,\fri)}_{q(\ell)} =a\right) \cr 
&= (5/8) \cdot \nu^{(n)}_{\cen}(\{a\}) 
+ (3/8) \cdot \nu^{(n)}_{\fri}(\{a\}) =
\nu^{(n)}_{\ord}(\{a\}). \cr}$$ 
Thus the last equality in statement (ii) in
(I) holds. That completes the proof of statement (I).

{\it Proof of (J).} Recall again that the sets
$Q(\ell)$, $\ell \in \{1,2,\dots,L\}$ are (pairwise)
disjoint (if $L \geq 2$). 
By statements (B), (C), and I(ii),
$$\cL\left(Y^{(n)}_{Q(\ell)} \, \Big| \, H^{(n)}(\cQ)\right) =
\cL\left(X^{(n)}_{Q(\ell)}\, \Big| \, H^{(n)} (\cQ)\right)$$ for
each $\ell \in \{1,\dots,L\}$. Hence statement (J) holds by
statements (G) and~(H).

{\it Proof of (K).} By statements (F) and (G), it suffices to prove
that for each $\ell\in \{1,\dots,L\}$ such that $S\cap Q(\ell)$ is
nonempty, one has that (see (2.4) and section 2.5(A) again)
$$\cL\left(X^{(n-1)}_{S\cap Q(\ell)} \, \Big| \, H^{(n)}(\cQ)\right) =
\cL\left( Y^{(n)}_{S\cap Q(\ell)}\, \Big| \, H^{(n)}(\cQ)\right).
$$ For $\ell\not\in \cI$ (see (8.12.2)), that holds trivially by
(8.13.3) and (8.13.4). For $\ell\in \cI$, it holds by statement (I)(i)(ii),
Lemma 3.7(A), and the assumption (in statement (K)) that 
$\card (S\cap
Q(\ell)) \le 5$. Thus (K) holds.

{\it Proof of (L) and (M).} First just suppose $S$ is a nonempty subset of
$Q(1) \cup Q(2) \cup \dots \cup Q(L)$. For a given vector
$(k(1),k(2),\dots,k(6)) \in S^6$ (the coordinates need not be
distinct), if the set $K:= \{k(1),k(2),\dots,k(6)\}$ intersects at
least two of the sets $Q(\ell)$ (and hence card$(K\cap Q(\ell)) \le
5$ for each $\ell\in \{1,\dots,L\})$, then by statements (F), (G),
and (K),
$$E\left[ \prod^6_{u=1} X^{(n-1)}_{k(u)} \, \bigg| \, H^{(n)}(\cQ)\right] =
E\left[ \prod^6_{u=1} Y^{(n)}_{k(u)} \, \bigg| \,
H^{(n)}(\cQ)\right]. \leqno{(8.13.12)}$$ 
Also, by (8.13.3) and
(8.13.4), equation (8.13.12) holds if $k(1),\dots,k(6)$ are elements
of the same set $Q(\ell)$ where $\ell\not\in \cI$. Hence the two
sides of (8.13.12) can differ only in (some) cases where
$k(1),\dots,k(6)$ are elements of the same set $Q(\ell)$ where
$\ell\in \cI$. Hence, using the notation $T\uparrow 6:= T\times
T\times T\times T\times T\times T$ for a given set $T$, one has that
$$\leqalignno{& &(8.13.13)\cr
&E\Bigg[ \left(\sum_{k\in S} X^{(n-1)}_k \right)^6 
\, \Bigg| \,
H^{(n)}(\cQ)\Bigg] - E\Bigg[\left( \sum_{k\in S} Y^{(n)}_k\right)^6
\, \Bigg| \, H^{(n)}(\cQ)\Bigg] \cr &= \sum_{(k(1),\dots,k(6)) \in
S\uparrow 6}([\hbox{LHS of (8.13.12)}] - [\hbox{RHS of (8.13.12)}])
\cr &=\sum_{\ell\in \cI} \sum_{(k(1),\dots,k(6)) \in (S\cap Q(\ell))
\uparrow 6}([\hbox{LHS of (8.13.12)}] - [\hbox{RHS of (8.13.12)}])
\cr 
&= \sum_{\ell\in \cI} \Bigg[E\Bigg[\left(\sum_{k\in S \cap
Q(\ell)}\!\!\!\!X^{(n-1)}_k\right)^6 \, \Bigg| \,
H^{(n)}(\cQ)\Bigg] \cr 
& \indent \indent \indent - E\Bigg[\left(\sum_{k\in S\cap
Q(\ell)}\!\!\!\! Y^{(n)}_k\right)^6 \, \Bigg| \,
H^{(n)}(\cQ)\Bigg]\, \Bigg]. \cr}$$ 
Hence by statement I(i)(ii) and
Lemma 3.7(B), equation (8.13.8) holds. Now if also $S\supset Q(\ell)$
for some $\ell\in \cI$, then by (8.13.13), statement I(i)(ii), and
Lemma 3.7(B)(C), equation (8.13.9) holds. That completes the proofs
of (L) and (M) and of Lemma~8.13.
\medskip

{\sc Definition 8.14.} Suppose $S$ is a nonempty finite subset of
$\Z$, and $n$ is a positive integer. A ``class $\cC(n)$ covering of
$S$'' is a family $\cQ := \{Q(1),Q(2),\dots,Q(L)\}$ of finitely many
subsets of $\Z$ (where $L$ is the number of sets in the family $\cQ$)
with the following four properties:

(i) For each $\ell\in \{1,\dots,L\}$, card~$Q(\ell)\in
\{1,6,6^2,\dots,6^n\}$.

(ii) The sets $Q(1),Q(2),\dots,Q(L)$ are (pairwise) disjoint (if
$L\ge 2$).

(iii) $S\subset Q(1)\cup Q(2)\cup \dots \cup Q(L)$.

(iv) For each $\ell\in \{1,\dots,L\}$, the set $S\cap Q(\ell)$ is
nonempty.
\medskip

{\sc Lemma 8.15.} {\it Suppose $S$ is a nonempty finite subset of $\Z$,
$n$ is a positive integer, and $\omega\in \Omega_0$. Then there
exists exactly one class $\cC(n)$ covering $\cQ$ of $S$ such that
(see (8.12.2)--(8.12.3)) $\omega\in H^{(n)}(\cQ)$. That covering
$\cQ$ is the family of all sets
$$D\in [\cD_0 (\omega) \cup \cD_1(\omega)\cup \dots \cup \cD_{n-1}
(\omega)] \cup \cE_n(\omega)$$
such that the set $D\cap S$ is nonempty.}

{\sc Proof.}  Suppose $S$, $n$ and $\omega$ are as in the statement
of Lemma 8.15. Referring to Definitions 8.1 and 8.3, define the
family $\cP_n$ of subsets of $\Z$ by
$$\cP_n := [\cD_0(\omega) \cup \cD_1(\omega) \cup \dots \cup \cD_{n-1}
(\omega) ] \cup \cE_n(\omega). \leqno{(8.15.1)}$$
Then by Remark 8.4(C), $\cP_n$ is a partition of the set $\Z$ itself
into countably many nonempty finite sets.

Let $\cQ^*$ denote the family of all sets $D\in \cP_n$ such that the
set $D\cap S$ is nonempty. Since the set $S$ is finite and the
members of $\cP_n$ (and hence of $\cQ^*$) are (pairwise) disjoint,
the family $\cQ^*$ has only finitely many sets. By (8.15.1), Lemma
8.5(A), and Definition 8.3(A)(B), the family $\cQ^*$ satisfies
property (i) in Definition 8.14. Properties (ii), (iii), and (iv) in
Definition 8.14 hold for $\cQ^*$ as trivial consequences of the
definition of $\cQ^*$ and the fact that $\cP^n$ is a partition of
$\Z$. Thus $\cQ^*$ is (by Definition 8.14) a class $\cC(n)$ covering
of~$S$.

Our next task is to show that $\omega\in H^{(n)}(\cQ^*)$. Refer to
property (i) in Definition 8.14 and refer to Lemma 8.5(A) again. For
each $m\in \{0,1,\dots,$ $n-1\}$ and each $D\in \cQ^*$ such that
$\card\, D=6^m$, one must have that $D\in \cD_m(\omega)$ by (8.15.1)
(since $D\in \cP_n$ but $D$ cannot belong to $\cD_u(\omega)$ or even
$\cE_u(\omega)$ for any $u\not= m$ by Lemma 8.5(A)), and hence
$\omega\in F_D$ by (8.6.1). For each $D\in \cQ^*$ such that card
$D=6^n$, one similarly has that $D\in \cE_n(\omega)$ by (8.15.1) and
hence $\omega\in G_D$ by (8.6.2). Hence by (8.12.2)--(8.12.3),
$\omega\in H^{(n)}(\cQ^*)$ (which is simply the intersection of the
various events $F_D$ and $G_D$ mentioned above).

By the definition of $\cQ^*$ and its properties identified above, all
of Lemma 8.15 has been proved except for uniqueness.

Suppose $\cQ$ is any class $\cC(n)$ covering of the set $S$ such that
$\omega\in H^{(n)}(\cQ)$. To prove uniqueness in Lemma 8.15, our task
is to show that $\cQ = \cQ^*$. Refer again to property (i) in
Definition 8.14. For each $m\in \{0,1,\dots,n-1\}$ and each $D\in
\cQ$ such that $\card\, D=6^m$, one has that $\omega\in F_D$ (since
$H^{(n)}(\cQ) \subset F_D$ by (8.12.2)--(8.12.3)), hence $D\in
\cD_n(\omega)$ by (8.6.1), hence $D\in \cP_n$ by (8.15.1). For each
$D\in \cQ$ such that $\card\,D=6^n$, one similarly has that $\omega\in
G_D$ (again by (8.12.2)--(8.12.3)), hence $D\in \cE_n(\omega)$ by
(8.6.2), hence $D\in \cP_n$ by (8.15.1). Hence $\cQ\subset \cP_n$.

If $D\in \cP_n$ and the set $D\cap S$ is empty, then (by property
(iv) in Definition 8.14)) $D$ cannot belong to $\cQ$. If $D\in \cP_n$
and (instead) the set $D\cap S$ is nonempty, then (by property (iii)
in Definition 8.14), $D$ must belong to $\cQ$ (since no other set in
$\cP_n$ will contain a given $k\in D\cap S$). Hence $\cQ = \cQ^*$ (by
the definition of $\cQ^*$). That completes the proof of uniqueness,
and of Lemma~8.15.
\medskip

{\sc Lemma 8.16.} {\it Suppose $n\in \N$. Then the following
four statements hold:

(A) The random sequences $Y^{(n)}$ and $X^{(n)}$ have the same
distribution (on $\{-1,1\}^{\Z}$).

(B) For any nonempty set $S\subset \Z$ such that $\card\, S\le 5$,
the random vectors $X^{(n-1)}_S$ and $Y^{(n)}_S$ (see (2.4)) have the
same distribution (on $\{-1,1\}^{\card \, S}$).

(C) For any nonempty finite set $S\subset\Z$, $E(\sum_{k\in
S}X^{(n-1)}_k)^6 \ge$  \break \noindent$E(\sum_{k\in
S}Y^{(n)}_k)^6$.

(D) For any integer $m\ge 6 \cdot 16^n$,}
$$E\left(\sum^m_{k=1}X^{(n-1)}_k\right)^6 \ge  E\left(\sum^m_{k=1}Y^{(n)}_k\right)^6 + 360
\cdot 4^{6(n-1)}.$$

{\sc Proof.} Before addressing any of the four statements (A)---(D)
in Lemma 8.16, let us present some arguments that will be common to
the proofs of all four statements.

Let $n\in \N$ be arbitrary but fixed.

Suppose $S$ is a nonempty finite subset of $\Z$.

Since there are only countably many finite subsets of $\Z$, there
exist only countably many finite families of finite subsets of $\Z$,
and hence (see Definition 8.14) only countably many class $\cC(n)$
coverings $\cQ$ of $S$. For each such $\cQ$, $H^{(n)}(\cQ) \subset
\Omega_0$ (again see (6.3.1)--(6.3.2)) by (8.6.1), (8.6.2), and
(8.12.2)--(8.12.3). By Lemma 8.15, the events $H^{(n)}(\cQ)$, for
class $\cC(n)$ coverings $\cQ$ of $S$, form a countable partition of
$\Omega_0$. (Some of those events $H^{(n)}(\cQ)$ may be empty.)

In the sums below, the index $\cQ$ ranges over all class $\cC(n)$
coverings of $S$ such that $P(H^{(n)}(\cQ))>0$. By a simple argument
(see section 2.5(A) and equation (2.4)),
$$\cL\left( Y^{(n)}_S\right) = \sum_{\cQ} \cL\left(Y^{(n)}_S \, \Big| \, 
H^{(n)}(\cQ)\right) \cdot P\left(H^{(n)}(\cQ)\right),
\leqno{(8.16.1)}$$
and the analogous statements holds with $Y^{(n)}_S$ replaced by $X^{(n-1)}_S$ and by $X^{(n)}_S$.

For each such covering $\cQ$ (with $P(H^{(n)}(\cQ))>0$), one has that
$\cL(Y^{(n)}_S \mid H^{(n)}(\cQ)) = \cL(X^{(n)}_S \mid H^{(n)}(\cQ))$
by Lemma 8.13(J) and property (iii) in Definition 8.14.  Hence by
(8.16.1) and its analog for $X^{(n)}_S$, one has that $\cL(Y^{(n)}_S)
= \cL(X^{(n)}_S)$. Since $S$ was an arbitrary nonempty finite subset
of $\Z$, statement (A) in Lemma 8.16 follows.

Next, in the case where $\card\, S\le 5$, an exactly analogous argument
holds with $X^{(n)}_S$ replaced by $X^{(n-1)}_S$, using Lemma
8.13(K). Thus statement (B) in Lemma 8.16 holds.

Next (regardless of the (finite) cardinality of $S$),
$$\leqalignno{&&(8.16.2) \cr
&E\left(\sum_{k\in S}X^{(n-1)}_k\right)^6 - E\left(\sum_{k\in
S}Y^{(n)}_k\right)^6 \cr &= \sum_{\cQ} E\left[\left(\sum_{k\in
S}X^{(n-1)}_k\right)^6\!\! - \left(\sum_{k\in S} Y^{(n)}_k\right)^6
\, \Bigg| \, H^{(n)}(\cQ)\right] \cdot P(H^{(n)}(\cQ)). \cr }$$ Hence
by Lemma 8.13(L), statement (C) in Lemma 8.16 holds.

{\it Proof of statement (D).} Suppose (for our given $n\in\N$) that
$m\ge 6 \cdot 16^n$. Now let $S:= \{1,2,\dots,m\}$.

Let $B$ denote the set (event) of all $\omega\in \Omega_0$ such that
there exist distinct integers $i,j\in \{1,2,\dots,6 \cdot 16^n\}$
such that $W^{(n)}_i(\omega) = W^{(n)}_j(\omega) =1$.

Suppose $\omega\in B$. Then by Definition 8.1 and a trivial argument,
there exists a set $E\in \cE_n(\omega)$ such that $E\subset
\{1,2,\dots,6\cdot 16^n\}$ (and hence $E\subset \{1,2,\dots,m\})$.
Now card $E=6^n$ by Lemma 8.5(A). Also, by Lemma 8.15, $E$ is a
member of the unique class $\cC(n)$ covering $\cQ$ of
$\{1,2,\dots,m\}$ such that $\omega\in H^{(n)}(\cQ)$.

Hence $B\subset \bigcup_{\cQ\in \cJ}H^{(n)}(\cQ)$ where $\cJ$ denotes
the family of all class $\cC(n)$ coverings $\cQ$ of $\{1,2,\dots,m\}$
such that $\cQ$ contains at least one member $E$ such that $E\subset
\{1,2,\dots, m\}$ and card $E=6^n$. Let $\cJ +$ denote the family of
all $\cQ\in \cJ$ such that $P(H^{(n)}(\cQ))>0$. Now by Lemma 6.7,
$P(B) \ge 1/2$. Hence  $\sum_{\cQ\in \cJ+} P(H^{(n)}(\cQ))\ge 1/2$.

Now by (8.16.2), together with Lemma 8.13(M) (for $\cQ\in \cJ+$) and
Lemma 8.13(L) (for the other $\cQ$'s such that $P(H^{(n)}(\cQ))>0)$,
$$\hbox{[LHS of (8.16.2)]} \ge \sum_{\cQ\in \cJ+} 720 \cdot
4^{6(n-1)} \cdot P(H^{(n)}(\cQ)) \ge 360 \cdot 4^{6(n-1)}.$$
Thus statement (D) in Lemma 8.16 holds. That completes the proof.
\medskip

{\sc Lemma 8.17.} {\it (A) For every $n\in \N$ and every nonempty
set $S\subset \Z$ with $\card \, S\le 5$, one has that (see (2.4)
and section 2.5(A))
$$\cL(X^{(n)}_S) = \cL(X^{(0)}_S).$$

(B) If $m$ and $N$ are positive integers such that $m\ge 6 \cdot
16^N$, then for every $n\ge N$,}
$$E\left( \sum^m_{k=1}X^{(n)}_k\right)^6 \le 15 m^3 - 360 \cdot
4^{6(N-1)}.$$
\medskip

{\sc Proof.}  Statement (A) follows from Lemma 8.16(A)(B) and
induction.

{\it Proof of (B).} Suppose $m$ and $N$ are as in statement (B).
Refer to Construction 7.4(A). Obviously $E(X^{(0)}_0)^\ell =1$ resp.\
0 if $\ell$ is an even resp.\ odd integer. Hence by a well known,
elementary calculation (the ``6th moment'' analog of the argument,
involving 4th moments, in [1, proof of Theorem 6.1]), one that
$E(\sum^m_{k=1} X^{(0)}_k)^6 \le 15 m^3$. Also, by Lemma 8.16 (A)(C),
the sequence of numbers $E(\sum^m_{k=1} X^{(n)}_k)^6$, $n\in
\{0,1,2,\dots\,\}$ is nonincreasing. Also, by Lemma 8.16(A)(D),
$$E\left(\sum^m_{k=1}X^{(N)}_k\right)^6 \le E\left(\sum^m_{k=1}
X^{(N-1)}_k\right)^6 - 360 \cdot 4^{6(N-1)}.$$
Combining these three preceding sentences, one obtains statement (B)
in Lemma 8.17.

\bigskip

\centerline {\bf 9.  Proof of properties (A), (B), (E), and (F) in Theorem 1.1}
\bigskip
This section contains arguments somewhat reminiscent of ones in [6]; but because of some nontrivial differences, the arguments here will be given in full. 

By Lemma 7.6 and Remark 7.7, the random sequence $X := (X_k$, $k\in
\Z)$ defined in (7.5.1) satisfies strict stationarity and Property
(C) in Theorem 1.1. Here in Section 9, it will be shown that the
sequence $X$ satisfies properties (A), (B), (E), and (F) in
Theorem~1.1. Finally in Section 10, property (D) (the triviality of
the double tail $\sigma$-field) will be verified, and the proof of
Theorem~1.1 will then be complete.

We start with some preliminary arguments. By (7.5.1), (8.9.1),
(8.9.2), and Constructions 7.3 and 7.4, the random variables $X_k$, $X_k^{(0)}$, 
$X^{(n)}_k$, $Y^{(n)}_k$, $k\in \Z$, $n\in \N$ all take their values in the set $\{-1,1\}$.

For each $k\in \Z$, $P(N_k=\infty) =0$ by Lemma 7.8. For each $k\in
\Z$ and each $\omega\in \Omega$ such that $N_k(\omega)<\infty$, one
has by (7.5.1) and (8.9.2) that $X_k(\omega) = X^{(n)}_k(\omega)$ for
all $n\ge N_k(\omega)$. Thus trivially,
$$\forall\, k\in {\Z},\quad X^{(n)}_k \rightarrow X_k \,\,\,\, \hbox{a.s.\ \ as
$n\to \infty$.}
\leqno{(9.1)}$$

Now let us prove properties (A), (B), (E), and (F) in Theorem~1.1.

{\it Proof of (A) and (B).} Suppose $S$ is any nonempty subset of
$\Z$ such that $\card\, S\le 5$. By (9.1) (see also (2.4)), $X^{(n)}_S
\rightarrow X_S$ a.s.\ as $n\to \infty$. As an elementary
consequence, $X^{(n)}_S \Rightarrow X_S$ (convergence in
distribution) as $n\to \infty$. Hence by Lemma 8.17(A) (see also Section 2.5(A)), $\cL(X_S ) = \cL(X^{(0)}_S)$.

That last equality has two consequences: First, taking 
$S = \{k\}$
for an arbitrary $k\in \Z$ and applying (7.4.1), one obtains property
(A) in Theorem~1.1. Second, taking arbitrary subsets 
$S\subset \Z$
with $\card\, S =5$, and applying the independence of the $X^{(0)}_k$'s
in Construction 7.4(A), one obtains property (B) in Theorem~1.1.

{\it Proof of property (E).} We shall now use the notations in (1.1).
Let $M$ be an arbitrary but fixed integer such that $M\ge 96$. To
prove property (E), it suffices to show that
$$E\left[ M^{-1/2}S(X,M)\right]^6 \le 15-16^{-6}.
\leqno{(9.2)}$$

Let $N$ be the positive integer such that
$$6 \cdot 16^N \le M< 6 \cdot 16^{N+1}. \leqno{(9.3)}$$
By (9.2) and Lemma 8.17(B),
$$\forall\,\, n\ge N, \quad E\left[S(X^{(n)},M)\right]^6 \le 15 M^3 -
360 \cdot 4^{6(N-1)}.
\leqno{(9.4)}$$

Now for each $n\in\N$, $[S(X^{(n)},M)]^6 \le M^6$ (since each
$X^{(n)}_k$ takes its values in $\{-1,1\}$). Also, as a trivial
consequence of (9.1), $[S(X^{(n)},M)]^6 \to [S(X,M)]^6$ a.s.\ as
$n\to \infty$. Hence by dominated convergence, $E[S(X^{(n)},M)]^6$
 $\to E[S(X,M)]^6$ as $n\to \infty$. Hence by (9.4), $E[S(X,M)]^6 \le
[$RHS of (9.4)$]$. Dividing both sides of this by $M^3$ and then
applying the second inequality in (9.3) and simple arithmetic, one
obtains
$$\eqalign{[\hbox{LHS of (9.2)}] &\le 15 - 360\cdot 4^{6(N-1)}/M^3
\cr \le \,\,&15 - [360 \cdot 16^{3(N-1)}]/[216 \cdot 16^{3(N+1)}] =
15 - (360/216) \cdot 16^{-6}. \cr}$$ Thus (9.2) holds, and property
(E) is proved.

{\it Proof of property (F).}  This is analogous to
a corresponding argument in [6]. However, it is short and
is also at the core of this construction. 
It will be given in full here.
  
By property (A) (proved above), $EX_0
=0$ and $EX^2_0 = EX^4_0 =1$. Hence by strict stationarity, property
(B), and the well known, elementary calculation in [1, proof of
Theorem 6.1] (which requires only 4--tuplewise independence), one has
that for each $n\in \N$ (again see (1.1)),
$$E\left(n^{-1/2}S_n\right) =0, \leqno{(9.5)}$$
$$E\left(n^{-1/2}S_n\right)^2 =1, \quad \hbox{and}
\leqno{(9.6)}$$
$$E\left(n^{-1/2}S_n\right)^4 \le 3.
\leqno{(9.7)}$$
By (say) (9.6) and Chebyshev's inequality, the family of
distributions of the normalized partial sums $(n^{-1/2} S_n$, $n\in
\N)$ is tight.

Now for the proof of property (F), suppose $Q$ is an infinite subset
of $\N$. Because of tightness, there exists an infinite set $T\subset
Q$ and a probability measure $\mu$ on $(\R,\cR)$ (both $T$ and $\mu$
henceforth fixed) such that
$$n^{-1/2}S_n \Rightarrow \mu \quad \hbox{as $n\to \infty$, $n\in
T$.} \leqno{(9.8)}$$
To complete the proof of property (F), our task is to show that $\mu$
is neither degenerate nor normal.

Because of (9.7), (9.8), and [1, p.\ 338, the Corollary], one has by
(9.5) and (9.6) that
$$\int_{x\in \R} x\mu(dx) =0 \quad\hbox{and}\quad \int_{x\in \R}x^2
\mu(dx) =1. \leqno{(9.9)}$$
Hence the probability measure $\mu$ has positive variance and is
therefore nondegenerate. Our task now is to prove that $\mu$ is not
normal.

If $\mu$ were normal, then by (9.9) it would have to be the $N(0,1)$
distribution, and hence by a well known calculation (e.g.\ [1, p.\
275, Example 21.1]) one would have $\int_{x\in \R} x^6 \mu(dx) = 15$.
By (9.8) and [1, p.\ 334, Corollary 1, and p.\ 338, Theorem 25.11],
one would then have
$$\liminf_{n\to \infty} E\left(n^{-1/2} S_n\right)^6 \ge 15.$$
But that contradicts property (E) (proved above). Hence $\mu$ cannot
be normal. That completes the proof of property (F).
\bigskip

\centerline {\bf 10.   Proof of property (D) in Theorem 1.1}
\bigskip 
This is the final
piece in the proof of Theorem~1.1. This argument will be divided into
thirteen ``steps'' (including some ``lemmas''), numbered 10.1, 10.2, etc.
\medskip

{\sc Step 10.1.} Let $F'$ be an arbitrary but fixed event such that
$$F' \in \cT_{\rm double} (X).
\leqno{(10.1.1)}$$ To prove property (D), i.e.\ to show that
$\cT_{\double}(X)$ is trivial, our task is to show that $P(F')=0$
or~1.

Suppose instead that
$$0<P(F')<1. \leqno{(10.1.2)}$$
Our task is to produce a contradiction. That task will consume the
rest of Section 10 here.
\medskip

{\sc Step 10.2.} By the assumption of (10.1.2), one has that
$[P(F')]^2<P(F')$. Let $\varepsilon\in (0,1)$ be fixed such that
$$[P(F') +\varepsilon]\cdot [P(F')+3\varepsilon]< P(F') - 4\varepsilon.
\leqno{(10.2.1)}$$

By (10.1.1), (7.4.2) (and its entire paragraph), 
(7.5.1), Lemma 7.6, and Definition 2.6, one
trivially has that $F'\in \sigma(X) \,\,\dot\subset\,\,
\sigma(\eta)$.  (One can do better, but this will suffice.)  Hence
from Constructions 6.1, 7.3, and 7.4 and a standard
measure--theoretic argument (again see (7.4.2)), there 
exists a positive integer $N'$ and
an event
$$\leqalignno{
F'' \in\,\, &\sigma\left(X^{(0)}_k, \, -N'\le k\le N'\right)
&(10.2.2)\cr &\vee \sigma\left(\xi^{(n)}_k, \,\, 1\le n\le
N',\,\,-N'\le k\le N'\right) \cr
&\vee
\sigma\left(\zeta^{(n)}_k,\,\, 1\le n\le N',\,\, -N'\le k\le
N'\right) \cr}$$ (see (7.3.1)) such that
$$P\left(F''\,\, \triangle\,\, F'\right) \le \varepsilon
\leqno{(10.2.3)}$$
(where $\Delta$ denotes symmetric difference).

By (10.2.3) and a trivial argument,
$$P(F'') \le P(F')+\varepsilon.
\leqno{(10.2.4)}$$

Increasing $N'$ if necessary, assume further that
$$N' \ge 1+ 6\cdot (8/5)^{30} \cdot(8/3)^6. \leqno{(10.2.5)}$$
(The positive integer $N'$ and the event $F''$, satisfying 
(10.2.2)--(10.2.5), are now taken as ``fixed.'') \medskip

{\sc Step 10.3.} Refer to section 6.3(A).  Random variables defined
in this step will be defined on $\Omega_0$ and left undefined on
$\Omega^c_0$.

For each $n\in \N$, we shall define a random positive integer $L_n$
and a sequence $(H(n,1),H(n,2), H(n,3),\dots\,)$ of random positive
integers such that for each $\omega\in \Omega_0$,
$$L_n(\omega) <H(n,1)(\omega)<H(n,2)(\omega) <H(n,3)(\omega)< \dots
\leqno{(10.3.1)}$$
and
$$\leqalignno{\qquad \qquad\bigl\{ k\in {\Z}: k>L_n(\omega&) \,\, \hbox{and}\,\,
W^{(n)}_k(\omega) \,=1\bigr\} &(10.3.2)\cr & = \{
H(n,1)(\omega),H(n,2)(\omega),H(n,3)(\omega),\dots  \}, \cr}$$ and
also (see section 6.1(C))
$$\sigma(L_n,H(n,1),H(N,2),H(n,3),\dots\,)\,\,\dot\subset\,\,
\sigma\left(\overline{\xi}^{(n)}\right). \leqno{(10.3.3)}$$
This will be done in such a way that
$$EL_1 \le N' +6\cdot (8/5)^{30} \cdot (8/3)^6. \leqno{(10.3.4)}$$
Also, for each $n\ge 2$, we shall define a random positive integer
$\Phi_n$, satisfying
$$\sigma(\Phi_n) \,\,\dot\subset \,\,
\sigma\left(\overline{\xi}^{(n)}\right) \quad \hbox{and}\quad E\Phi_n
\le 6 \cdot(8/5)^{30} \cdot (8/3)^6. \leqno{(10.3.5)}$$
The definition will be recursive, and is as follows:

To start with $n=1$, refer to Statement 6.3(A)(ii) (and Section \hfil\break 
2.4(A)(B)(C)(D)). 
Define the positive integer-valued random variable $L_1$ as
follows: For every $\omega\in \Omega_0$,
$$\leqalignno{\qquad\qquad L_1(\omega) := \min \bigl\{k\in {\Z}: k\ge N'+6 \,\, \hbox{and}\,
\bigl[ \xi^{(1)t}_{k-5}&(\omega)\, \mid \xi^{(1)t}_{k-4}(\omega)\mid
&(10.3.6) \cr & \cdots \mid \xi^{(1)t}_k(\omega)\bigr] = I_6\bigr\}
\cr}$$ 
(where as usual, the superscript $t$ denotes transpose).
By Lemma 5.2(C) (see Section 6.1(B)), 
equation (10.3.4) holds. Now define the
random positive integers $H(1,i)$, $i\in\N$ uniquely by (10.3.1) and
(10.3.2) (for $n=1$) for $\omega\in \Omega_0$. Trivially by (10.3.6), (10.3.1)--(10.3.2), and (6.2.1) (and 
Remark 2.6(B)),
equation (10.3.3) holds for $n=1$.

Now suppose $n\ge 2$, and the positive integer--valued random
variables $L_{n-1}$ and $H(n-1,i)$, $i\in \N$ are already defined,
satisfying (10.3.1), (10.3.2), and (10.3.3) with $n$ replaced by
$n-1$.

Referring to (10.3.1)--(10.3.2) and Statement 6.3(A)(iii) (with $n$
replaced by $n-1$), define the positive integer-valued random
variable $\Phi_n$ as follows: For each $\omega\in \Omega_0$,
$$\leqalignno{\qquad \qquad\Phi_n(\omega) :=  \min \Bigl\{k\ge 6:
\Bigl[& \xi^{(n)t}_{H(n-1,k-5)(\omega)}(\omega)\mid
\xi^{(n)t}_{H(n-1,k-4)(\omega)}(\omega)\mid &(10.3.7) \cr
&\qquad\quad\quad\quad \,\,\cdots\mid
\xi^{(n)t}_{H(n-1,k)(\omega)}(\omega)\Bigr] =I_6\Bigr\}. \cr}
$$ Next,
define the positive integer-valued random variable $L_n$ as follows:
For all $\omega\in \Omega_0$,
$$L_n(\omega) := H(n-1,\Phi_n(\omega))(\omega). \leqno{(10.3.8)}$$
Now define the random positive integers $H(n,i)$, $i\in \N$
(uniquely) by (10.3.1) and (10.3.2) (for the given $n$) for
$\omega\in \Omega_0$.

By (10.3.3) (with $n$ replaced by $n-1$) and (10.3.7) (for the given
$n$), the first part $(\sigma(\Phi_n) \,\, \dot\subset\,\,
\sigma(\overline{\xi}^{(n)}))$ of (10.3.5) holds; hence by (10.3.3)
(with $n$ replaced by $n-1$) and (10.3.8) (for the given $n$),
$\sigma(L_n) \,\,\dot\subset\,\, \sigma(\overline{\xi}^{(n)})$; and
hence (10.3.3) holds for the given $n$ by (10.3.1)--(10.3.2) and
(6.2.1). To complete this recursive definition, all that remains is
to verify the second part $(E\Phi_n \le \dots\,)$ in (10.3.5).

For that purpose, we shall apply Lemma 5.2(D), with the random
integers $\kappa(j)$, $j\in\Z$ there defined by $\kappa(j) :=
H(n-1,j)$ for $j\ge 1$ and (just as a frivolous formality) $\kappa(j)
:= j$ (a constant random variable) for $j\le 0$ (see (10.3.1)). Then
$\sigma(\kappa(j))$, $j\in
{\Z})\,\,\dot\subset\,\,\sigma(\overline{\xi}^{(n-1)})$ by (10.3.3)
with $n$ replaced by $n-1$. Since the sequences $\xi^{(n)}$ and
$\overline{\xi}^{(n-1)}$ are independent (and $\xi^{(n)}$ has the
appropriate distribution), all conclusions of Lemma 5.2(D) apply to
the sequence $(\xi^{(n)}_{\kappa(j)}$, $j\in \Z)$. Applying Lemma 5.2(D)(iv) to that sequence, 
one obtains that the random variable
$\Phi_n$ (see (10.3.7)) satisfies the second part of (10.3.5). That
completes this recursive definition.

By (10.3.6), (10.3.1), (10.3.7) (which gives $\Phi_n(\omega) \ge 6$ for $n\ge
2$ and $\omega\in \Omega_0)$, (10.3.8), and induction, one has that
for all $\omega\in \Omega_0$,
$$\leqalignno{\qquad N'+6 \le L_1(\omega) <H(1,1)(\omega)
&<L_2(\omega)<H(2,1)(\omega) &(10.3.9) \cr &<L_3(\omega)<
H(3,1)(\omega)<\dots\,. \cr}$$

{\sc Step 10.4.} The ``universe'' of ``scaffolding'' in Sections 6,
7, and 8 was based on the following array of independent random variables from 
Constructions 6.1, 7.3 (recall (7.3.1)), and 7.4:
$$ X^{(0)}_k,\,\, k\in{\Z};\,\, \hbox{and}\,\,
\xi^{(n)}_k,\zeta^{(n,\ord)}_k,
\zeta^{(n,\cen)}_k,\zeta^{(n,\fri)}_k,\,\, (n,k)\in\N\times\Z.\quad
\leqno{(10.4.1)}$$ Here in Step 10.4, we shall construct an
``alternate universe'' based on an array
$$X^{*(0)}_k,\, k\in {\Z}; \quad\hbox{and}\quad
\xi^{*(n)}_k,\zeta^{*(n,\ord)}_k, \zeta^{*(n,\cen)}_k,
\zeta^{*(n,\fri)}_k, \,\, (n,k)\in\N\times\Z. \leqno{(10.4.2)}$$
To create this ``alternate universe,'' we shall keep the ``original''
random variables in the array (10.4.1), except that the random
variables in (10.4.1) that are mentioned in the right hand side of
(10.2.2) will be replaced by ``independent copies.''

(A) Refer to the positive integer $N'$ in (10.2.2) and (10.2.5). Let
$$\leqalignno{\qquad X^{*(0)}_k,
 k\in \{-&N',-N'+1,-N'+2,\dots,N'\}\,\,\,\hbox{and} &(10.4.3)\cr
 & \xi^{*(n)}_k,
\zeta^{*(n,\ord)}_k,\zeta^{*(n,\cen)}_k,\zeta^{*(n,\fri)}_k,  \cr &
n\in \{1,2,\dots,N'\},\quad k\in \{-N',-N'+1,\dots,N'\} \cr}$$ be an
array of independent random variables, with this array being
independent of the entire array (10.4.1) (and hence independent of
the entire collection of random variables studied in sections 6, 7,
8, and 9), such that for each $k\in \{-N',-N'+1,\dots,N'\}$ and
$n\in \{1,2,\dots,N'\}$, the random variable $X^{*(0)}_k$ resp.\
$\xi^{*(n)}_k$ resp.\ $\zeta^{*(n,\ord)}_k$ resp.\
$\zeta^{*(n,\cen)}_k$ resp.\ $\zeta^{*(n,\fri)}_k$ takes its values
in the space $\{-1,1\}$ resp.\ $\{0,1\}^6$ resp.\
$\{-1,1\}^{\sxtp(n)}$ (see (2.1) again) resp.\ $\{-1,1\}^{\sxtp(n)}$
resp.\ $\{-1,1\}^{\sxtp(n)}$ and has the same distribution as the
random variable $X^{(0)}_k$ resp.\ $\xi^{(n)}_k$ resp.\
$\zeta^{(n,\ord)}_k$ resp.\ $\zeta^{(n,\cen)}_k$ resp.\
$\zeta^{(n,\fri)}_k$.

(B) For each $k\in {\Z} - \{-N',-N'+1,\dots,N'\}$, define the
$\{-1,1\}$--valued random variable $X^{*(0)}_k$ by $X^{*(0)}_k :=
X^{(0)}_k$ (that is, for all $\omega\in \Omega$, $X^{*(0)}_k(\omega)
:= X^{(0)}_k(\omega))$.

(C) For each $(n,k) \in ({\N}\times {\Z}) -
(\{1,2,\dots,N'\}\times\{-N',-N'+1,\dots,N'\})$, define the
$\{0,1\}^6$--valued random variable $\xi^{*(n)}_k$ by $\xi^{*(n)}_k
:= \xi^{(n)}_k$, and define the $\{-1,1\}^{\sxtp(n)}$--valued random
variables (see (2.1)) $\zeta^{*(n,\ord)}_k$, $\zeta^{*(n,\cen)}_k$,
and $\zeta^{*(n,\fri)}_k$ by $\zeta^{*(n,\ord)}_k :=
\zeta^{(n,\ord)}_k$, $\zeta^{*(n,\cen)}_k:= \zeta^{(n,\cen)}_k$, and
$\zeta^{*(n,\fri)}_k := \zeta^{(n,\fri)}_k$. That completes the
definition of the array (10.4.2).

(D) Recall from Constructions 6.1, 7.3, and 7.4 that the random
variables in the array (10.4.1) are independent of each other. It
follows from the conditions in (A), (B), and (C) above that the
random variables in the array (10.4.2) are independent of each other,
and that in fact the two arrays (10.4.1) and (10.4.2) have the same
distribution (on the appropriate space).

(E) Further, the array (10.4.2) is independent of the entire
$\sigma$--field in the right hand side of (10.2.2). Hence by (10.2.2)
itself,
$$\hbox{the event $F''$ is independent of the array (10.4.2).}
\leqno{(10.4.4)}$$

(F) Starting with (A), (B), (C), and (D) above, we construct analogs
of all of the random variables constructed in Sections 6, 7, and 8,
using the array (10.4.2) in place of the array (10.4.1). The
notations will be the same except that (as in (10.4.2)) an asterisk
will be inserted after each ``main symbol.'' The general pattern is
that for a given random variable $Z$ of the form $Z:= f$(the array
(10.4.1)), where $f$ is an appropriate measurable function, the
analogous random variable $Z^*$ will be given by $Z^* =f$(the array
(10.4.2)) (with the same function $f$). Thus the analog of (6.1.3)
is
$$\overline{\xi}^{*(n)}_k :=
\left(\xi^{*(1)}_k,\xi^{*(2)}_k,\dots,\xi^{*(n)}_k\right),$$
the analog of (6.2.5) is
$$W^{*(1)}_k := g_{\spaced}\left(\xi^{*(1)}_k, \xi^{*(1)}_{k-1},
\xi^{*(1)}_{k-2},\dots\,\right),$$
the analog of (6.2.6) is
$$\Psi^*(n,k,j) := \psi_j\left(W^{*(n)}_k,
W^{*(n)}_{k-1},W^{*(n)}_{k-2},\dots\,\right),$$
and so on. Along the way, corresponding analogs of the events
$\Omega^{(n)}_{\good}$ and $\Omega_0$ in Section 6.2 and equation
(6.3.1) are defined and denoted $\Omega^{*(n)}_{\good}$ and
$\Omega^*_0$.

Of course because of (B) and (C) above, some of the new
(``asterisk'') random variables defined in this way will be exactly
the same as the original counterparts in Sections 6, 7, and~8.

(G) Now define the event
$$\Omega' := \Omega_0 \cap \Omega^*_0. \leqno{(10.4.5)}$$
By (6.3.2) and its ``asterisk'' counterpart $(P(\Omega^*_0) =1)$,
one has that
$$P(\Omega')=1. \leqno{(10.4.6)}$$
Also (for example), for every $\omega\in \Omega'$, 
Statements 6.3(A)(i)(ii)(iii)
and their ``asterisk'' counterparts all hold.

(H) The random variables formulated in Step 10.3 
($L_n$ and $H(n,i)$ for $n,i \in \N$, and
$\Phi_n$ for $n \geq 2$), are defined 
(and positive integer-valued)
at every $\omega \in \Omega'$, by (10.4.5) and the second
sentence in Step 10.3 itself.
For those random variables, we will not need to refer
to, and will therefore not formally define, ``asterisk'' counterparts.    
\medskip

{\sc Step 10.5.} (A) Refer to Step 10.4(B)(C). The portion of the
array (10.4.2) involving indices $k\le -N'-1$ (and any $n\in\N$)
coincides with the corresponding portion of the array (10.4.1).

(B) In particular, for any $k\le -N'-1$, $\xi^{*(1)}_k = \xi^{(1)}_k$
(that is, $\xi^{*(1)}_k(\omega) = \xi^{(1)}_k(\omega)$ for all
$\omega\in \Omega$). Hence $W^{*(1)}_k = W^{(1)}_k$ for every $k\le -N'
-1$, since by (6.2.5) (in both arrays (10.4.1) and (10.4.2)), for any
$k\le -N'-1$ and any $\omega\in\Omega$,
$$\leqalignno{\qquad \quad W^{*(1)}_k(\omega) &= g_{\spaced}
\left(\xi^{*(1)}_k(\omega),\xi^{*(2)}_k(\omega),
\xi^{*(3)}_k(\omega),\dots\,\right) &(10.5.1) \cr
&= g_{\spaced} \left(\xi^{(1)}_k(\omega), \xi^{(2)}_k(\omega),
\xi^{(3)}_k(\omega),\dots\,\right) = W^{(1)}_k(\omega). \cr}$$
By the same type of argument, using (6.2.6) and (6.2.7) (in both
arrays (10.4.1) and (10.4.2)) together with induction on $n$, one can
show the following two facts together: that $\Psi^*(n,k,j) =
\Psi(n,k,j)$ for each $n\in \N$, each $k\le -N' -1$, and each $j\ge
0$, and that
$$\forall\, n\in{\N}, \,\, \forall\, k\le -N'-1, \quad W^{*(n)}_k =
W^{(n)}_k. \leqno{(10.5.2)}$$

(C) Continuing to apply (A) above in the same way as in (B), one
obtains for any $k\le -N'-1$ and $n\in \N$ the following equalities
of random variables: $\overline{\xi}^{*(n)}_k =
\overline{\xi}^{(n)}_k$ (see (6.1.3)), $\overline{W}^{*(n)}_k =
\overline{W}^{(n)}_k$ (see (6.3.6)), $\delta^{*(n)}_k =
\delta^{(n)}_k$ (see (7.1.1) and (7.1.2)), $N^*_k = N_k$ (see
(7.1.4)), $J^*(n,k) = J(n,k)$ (see (7.1.7)), $\zeta^{*(n)}_k =
\zeta^{(n)}_k$ (see (7.3.1)), $\overline{\zeta}^{*(n)}_k =
\overline{\zeta}^{(n)}_k$ (see (7.3.2)), and finally, from (7.5.1),
$$\forall\, k\le -N'-1,\qquad X^*_k = X_k. \leqno{(10.5.3)}$$

{\sc Lemma 10.6.} {\it Refer to Step 10.3 (including equation
(10.3.9)).  Refer also to (10.4.5)--(10.4.6) and 
Step 10.4(H). 
Suppose $n\in\N$. Then for every 
$\omega\in\Omega'$, the following statements hold:}
$$\forall\,\, k\ge L_n(\omega) +1, \quad W^{*(n)}_k(\omega) =
W^{(n)}_k(\omega); \leqno{(10.6.1)}$$
$$\forall\,\, k\ge H(n,1)(\omega), \quad \delta^{*(n)}_k(\omega) =
\delta^{(n)}_k(\omega); \leqno{(10.6.2)}$$
$$\forall \,\, k\ge H(n,1)(\omega), \quad \Psi^*(n,k,0)(\omega) =
\Psi(n,k,0)(\omega); \,\,\hbox{and} \leqno{(10.6.3)}$$
$$\forall\,\, k\ge H(n,1)(\omega), \quad k-\Psi^*(n,k,0)(\omega) =
k-\Psi(n,k,0)(\omega) \ge H(n,1)(\omega). \leqno{(10.6.4)}$$

{\sc Proof.} Throughout this proof, let $\omega \in \Omega'$
be arbitrary but fixed.  (Again recall Step 10.4(G)(H).) 

For each $n\in\N$, let the function 
$\varphi_n: (\{0,1\})^n \to
\{0,1,2,3,4,5,6\}$ be as in Lemma 4.7.

We shall first verify (10.6.1) and (10.6.2) for $n=1$. For each $i \geq L_1(\omega) - 5$, one has that
$\xi_i^{*(1)}(\omega) = \xi_i^{(1)}(\omega)$
by (10.3.9) and Step 10.4(C). 
Hence for each integer $k \geq L_1(\omega) + 1$, 
by the ``asterisk'' counterpart of (6.2.5), 
then two applications of (10.3.6) and Lemma 4.7
together, followed by (6.2.5) itself, with
$L_1(\omega)$ written here as $L(1)(\omega)$,
$$\eqalign{W^{*(1)}_k(\omega) &= g_{\spaced}
\left(\xi^{*(1)}_k(\omega),\xi^{*(1)}_{k-1}(\omega),\xi^{*(1)}_{k-2}(\omega),
\dots\, \right) \cr &=
\varphi_{k-L(1)(\omega)}\left(\xi^{*(1)}_{L(1)(\omega)+1}(\omega),\xi^{*(1)}_{L(1)(\omega)
+2} (\omega),\dots,\xi^{*(1)}_{k}(\omega)\right) \cr
&=\varphi_{k-L(1)(\omega)}\left(\xi^{(1)}_{L(1)(\omega)+1}(\omega),\xi^{(1)}_{L(1)(\omega)+2}(\omega),\dots,\xi^{(1)}_k(\omega)\right)
\cr &=
g_{\spaced}\left(\xi^{(1)}_k(\omega),\xi^{(1)}_{k-1}(\omega),\xi^{(1)}_{k-2}(\omega),\dots\,\right)
\cr &=W^{(1)}_k(\omega). \cr}$$ 
Thus (10.6.1) holds for $n=1$. 
By (10.6.1) (for $n = 1$) and (7.1.1) (and its ``asterisk'' counterpart),
$\delta^{*(1)}_k(\omega) = \delta^{(1)}_k(\omega)$ for all
$k\ge L_1(\omega) + 1$ (and in particular
for all $k\ge H(1,1)(\omega)$ by (10.3.9)). Thus (10.6.2) holds for $n=1$.

Now we use induction. Suppose $N\in\N$, and
suppose (for the given $\omega \in \Omega'$) that (10.6.1) 
and (10.6.2) hold for $n=N$. For the induction step,
we shall show (for the given $\omega\in \Omega'$) that
(10.6.3) and (10.6.4) hold for $n=N$, and then show that (10.6.1) and (10.6.2)
hold for $n=N+1$. This actually suffices to prove Lemma 10.6 by
induction.

{\it Verification of (10.6.3)--(10.6.4) for $n=N$.} Suppose $k\ge
H(N,1)(\omega)$. Referring to (10.3.1) and (10.3.2), let $m$ denote
the positive integer such that $H(N,m)(\omega) \le
k<H(N,m+1)(\omega)$. By the inductive assumption of (10.6.1) for
$n=N$, equation (10.3.2) holds for $n=N$ with $W^{(N)}_k(\omega)$
replaced by $W^{*(N)}_k(\omega)$. Hence by (6.2.6) and its
``asterisk'' analog (and Definition 2.2(A)), together with (10.3.1)
and (10.3.2),
$$\eqalign{
\Psi^*(N,k,0)(\omega) &=
\psi_0\left(W^{*(N)}_k(\omega),W^{*(N)}_{k-1}(\omega),W^{*(N)}_{k-2}(\omega),\dots\,\right)
\cr &= k-H(N,m)(\omega) \cr &=
\psi_0\left(W^{(N)}_k(\omega),W^{(N)}_{k-1}(\omega),W^{(N)}_{k-2}(\omega)
,\dots\,\right)\cr &= \Psi(N,k,0)(\omega), \cr}$$ and hence also
(see (10.3.1) again)
$$k-\Psi^*(N,k,0)(\omega) = k-\Psi(N,k,0)(\omega) = H(N,m)(\omega)
\ge H(N,1)(\omega).$$
Thus (10.6.3) and (10.6.4) hold for $n=N$.

{\it Verification of (10.6.1) for} $n=N+1$.  By (10.3.8) (see (10.3.7)) $L_{N+1}(\omega) =
H(N,\Phi_{N+1}(\omega))(\omega)$. 
If $k\ge L_{N+1}(\omega)+1$ but
$k$ is not one of the integers $H(N,m)(\omega)$, $m\ge
\Phi_{N+1}(\omega)+1$, then 
$W^{(N)}_k(\omega) \not= 1$ by ((10.3.1) and) (10.3.2)
for $n = N$, hence $W_k^{*(N)} = W_k^{(N)} \not= 1$ by 
((10.3.9) and) the inductive assumption of (10.6.1)
for $n = N$, and hence  
$W^{*(N+1)}_k(\omega) = W^{(N+1)}_k(\omega) =0$ by both
Remark 6.3(B) and its `` asterisk'' counterpart.
Hence (see (10.3.1) again) to complete the proof of (10.6.1) for $n = N+1$,
what remains is to show that for each integer 
$m \geq \Phi_{N+1}(\omega) + 1$, one has that
$W_{H(N,m)(\omega)}^{*(N+1)}(\omega)
= W_{H(N,m)(\omega)}^{(N+1)}(\omega)$. 

Suppose $m\ge \Phi_{N+1}(\omega) +1$.  To obtain the
desired equality, using the function $\varphi_u$ in
Lemma 4.7 for the positive integer $u := m - \Phi_{N+1}(\omega)$, using the convention in
Notation 2.1(G), and writing the positive integer 
$\Phi_{N+1}(\omega)$ as
$\Phi(N+1)(\omega)$, we shall show that
$$\leqalignno{&&(10.6.5) \cr
&\ W^{*(N+1)}_{H(N,m)(\omega)}(\omega)  \cr 
&= g_{\spaced}
\left(\xi^{*(N+1)}_j,\, \, j\in \{i\le H(N,m)(\omega):
W^{*(N)}_i(\omega) =1\}\right) \cr &= 
\varphi_u
\left(\xi^{*(N+1)}_{H(N,\Phi(N+1)(\omega)+1)(\omega)}(\omega),
\xi^{*(N+1)}_{H(N,\Phi(N+1)(\omega)+2)(\omega)}(\omega),\dots,\xi^{*(N+1)}_{H(N,m)(\omega)}(\omega)\right)
\cr &=
\varphi_u
\left(\xi^{(N+1)}_{H(N,\Phi(N+1)(\omega)+1)(\omega)}(\omega),\xi^{(N+1)}_{H(N,\Phi(N+1)(\omega)+2)(\omega)}(\omega),\dots,
\xi^{(N+1)}_{H(N,m)(\omega)}(\omega)\right) \cr &= g_{\spaced}
\left(\xi^{(N+1)}_{j},\,\, j\in\{i\le
H(N,m)(\omega):W^{(N)}_i(\omega) =1\}\right) \cr &=
W^{(N+1)}_{H(N,m)(\omega)}(\omega). \cr}$$
To verify this, first note that by (10.3.9) and Step 10.4(C),
the equality 
$\xi_k^{*(N+1)}(\omega) = \xi_k^{(N+1)}(\omega)$ 
holds for
each $k \geq H(N,1)(\omega)$, {\it and in particular it holds for each $k$ of the form $k = H(N,\ell)(\omega)$ for
$\ell \geq \Phi_{N+1}(\omega) - 5$} by (10.3.1) (with $n=N$)
and (10.3.7) (with $n = N+1$).  That trivially yields the
third equality in (10.6.5).
Next, by (10.3.1)--(10.3.2) and the inductive assumption of (10.6.1) for $n = N$, 
$W_k^{*(N)}(\omega) = W_k^{(N)}(\omega) = 1$
for $k \in \{ H(N,1)(\omega), \allowbreak H(N,2)(\omega),\  
H(N,3)(\omega),\ \dots\}$, and
$W_k^{*(N)}(\omega) = W_k^{(N)}(\omega) \not= 1$ 
for all other $k \geq H(N,1)(\omega)$.
This fact has a couple of consequences:  
First (just consider that fact for $k = H(N,m)(\omega)$),  
by (6.2.14) and (6.2.9)--(6.2.10) and their
``asterisk'' counterparts (for $n = N$), the first and fifth 
(i.e.\ last) equalities in (10.6.5) hold.
Then, by (10.3.7) with $n = N+1$ (combined with the observation in italics in
the first sentence after (10.6.5)), one has from Lemma 4.7
that the second and fourth equalities in (10.6.5) also hold.   
Thus (10.6.5) holds.  That completes the proof of
(10.6.1) for $n = N+1$. 

{\it Verification of (10.6.2) for} $n=N+1$. Now suppose $k\ge
H(N+1,1)(\omega)$. Referring to (10.3.1) for $n=N+1$, let $m\in\N$
be such that $H(N+1,m)(\omega) \le k<H(N+1,m+1)(\omega)$. Now
$W^{*(N+1)}_{H(N+1,m)(\omega)}(\omega) =
W^{(N+1)}_{H(N+1,m)(\omega)} (\omega) =1$ by (10.3.2) (for $n=N+1$)
and (10.6.1) for $n=N+1$ (just verified above), and hence
$W^{*(N)}_{H(N+1,m)(\omega)}(\omega) =1$ and
$W^{(N)}_{H(N+1,m)(\omega)}(\omega) =1$ by Remark 6.3(B) (and its
``asterisk'' counterpart). Hence $\Psi(N,k,0)(\omega) \le
k-H(N+1,m)(\omega)$ by (6.2.6) (and Definition 2.2(A)), hence
$k-\Psi(N,k,0)(\omega) \ge H(N+1,m)(\omega) \ge L_{N+1}(\omega) +1$
(see (10.3.2) with $n=N+1$), hence
$$W^{*(N+1)}_{k-\Psi(N,k,0)(\omega)}(\omega) =
W^{(N+1)}_{k-\Psi(N,k,0)(\omega)}(\omega) \leqno{(10.6.6)}$$ by
(10.6.1) for $n=N+1$ (verified above). Also (since $k\ge
H(N+1,m)(\omega) \ge H(N+1,1)(\omega)>H(N,1)(\omega)$ by (10.3.9))
$\Psi^*(N,k,0)(\omega) = \Psi(N,k,0)(\omega)$ by (10.6.3) for $n=N$
(verified above). Substituting that into the left hand side 
of (10.6.6) and then applying (7.1.2) and its ``asterisk'' counterpart,
one obtains $\delta^{*(N+1)}_k(\omega) = \delta^{(N+1)}_k(\omega)$.
Thus (10.6.2) has been verified for $n=N+1$. That completes the
induction step and the proof of Lemma~10.6.
\medskip

{\sc Step 10.7.} For each $n\in\N$, each $\omega\in\Omega'$ (see Step
10.4(G)), and each set $S\subset \R$, define the sets
$$G^{(n)}(S)(\omega) := \left\{k\in {\Z}: k\in S\ \
\hbox{and}\ \ 
W^{(n)}_k(\omega) =1\right\} \leqno{(10.7.1)}$$
and
$$G^{*(n)}(S)(\omega) := \left\{k\in {\Z}: k\in S \ \ \hbox{and}\ \ 
W^{*(n)}_k(\omega) =1\right\}. \leqno{(10.7.2)}$$

For each $n\in\N$, define the random sets $\Gamma_n$ and $\Gamma^*_n$
as follows: For each $\omega\in\Omega'$ (see Step 10.4(G) again),
$$\leqalignno{\qquad\qquad\Gamma_n(\omega) :=& G^{(n)}([-N',L_n(\omega)])(\omega) &(10.7.3)\cr &=  \, \Bigl\{
k\in{\Z}: -N' \le k\le L_n(\omega)\,\, 
\hbox{and}\,\, W^{(n)}_k(\omega) =1\Bigr\} \cr}$$ 
and
$$\leqalignno{\qquad\qquad\Gamma^*_n(\omega) :=&  G^{*(n)} ([-N',L_n(\omega)])(\omega) &(10.7.4) \cr &= \, 
\Bigl\{ k\in{\Z}: -N' \le k\le L_n(\omega)\,\, 
\hbox{and}\,\, W^{*(n)}_k(\omega) =1\Bigr\}. \cr}$$

For each $n\in\N$, by (10.7.3), (10.7.4), (10.3.3), and (6.2.1) and
its ``asterisk'' counterpart, one has that
$$\sigma\left(\Gamma_n, \Gamma^*_n\right)\,\,\dot\subset\,\,
\sigma\left(\overline{\xi}^{(n)},\overline{\xi}^{*(n)}\right).
\leqno{(10.7.5)}$$

Also, for $n\in\N$ and $\omega\in\Omega'$, by (6.2.9) and (6.2.10), one
has (see (10.7.1)) the following reformulation of (6.2.14) 
in the convention of Notation 2.1(G):
$$\leqalignno{\qquad\quad&\hbox{If $k\in\Z$ is such that $W^{(n)}_k(\omega) =1$,} &(10.7.6)\cr
&\hbox{then $W^{(n+1)}_k(\omega) = g_{\spaced}\left(\xi^{(n+1)}_j(\omega),
\,j\in G^{(n)}((-\infty,k])(\omega)\right)$.}  \cr}$$ Analogously
(see (10.7.2)), for $n\in\N$ and $\omega\in\Omega'$,
$$\leqalignno{\qquad\qquad&\hbox{If $k\in\Z$ is such that $W^{*(n)}_k(\omega) =1$,} &(10.7.7)\cr
&\hbox{then $W^{*(n+1)}_k(\omega) =
g_{\spaced}\left(\xi^{*(n+1)}_j(\omega),\, j\in
G^{*(n)}((-\infty,k])(\omega)\right)$.}  \cr}$$

{\sc Lemma 10.8.} {\it For each $n\in\N$, }
$$E(\card\, \Gamma_n)\le 3N' \quad\hbox{\it and} \leqno{(10.8.1)}$$
$$E(\card\, \Gamma^*_n) \le 3N'. \leqno{(10.8.2)}$$

{\sc Proof.} It suffices to give the argument for (10.8.1). The
argument for (10.8.2) is exactly analogous.

First, for each $\omega\in\Omega'$, by (10.7.3), $\Gamma_1(\omega)
\subset \{ -N', -N'+1,\dots,L_1(\omega)\}$ and hence $\card\,
\Gamma_1(\omega) \le 1+N' +L_1(\omega)$. Hence by (10.3.4) and
(10.2.5),
$$\leqalignno{
E(\card\, \Gamma_1) &\le 1+N'+EL_1 &(10.8.3) \cr
&\le 1+N' +N' +6\cdot (8/5)^{30} \cdot (8/3)^{6} \cr
&\le 3N'. \cr}$$

Next, for each $n\ge 2$ and each $\omega\in\Omega'$, by (10.7.3),
(10.7.1), Remark 6.3(C), (10.3.9), (10.3.1), (10.3.2), and (10.3.8),
$$\eqalign{
\card\, \Gamma_n(\omega)
 &\le 1+(1/6) \cdot \card\,
G^{(n-1)}([-N',L_n(\omega)])(\omega) \cr &= 1+(1/6)\cdot \card\,
G^{(n-1)}([-N', L_{n-1}(\omega)])(\omega) \cr &\qquad + (1/6) \cdot
\card\, G^{(n-1)}((L_{n-1}(\omega), L_n(\omega)])(\omega) \cr &=
1+(1/6)\cdot \card\, \Gamma_{n-1}(\omega) \cr &\qquad +(1/6) \cdot
\card \bigl\{H(n-1,1)(\omega), H(n-1,2)(\omega), \cr
&\qquad\qquad\qquad\qquad\qquad\qquad\dots,
H(n-1,\Phi_n(\omega))(\omega)\bigr\} \cr &=1+(1/6) \cdot \card\,
\Gamma_{n-1}(\omega) + (1/6) \cdot \Phi_n(\omega). \cr}$$ Hence for
each $n\ge 2$, by (10.3.5) and (10.2.5),
$$\eqalign{
E(\card\, \Gamma_n) &\le 1+(1/6) \cdot E(\card\, \Gamma_{n-1})+(1/6)
\cdot E\Phi_n \cr
&\le 1+(1/6) \cdot E(\card\, \Gamma_{n-1}) +N'. \cr}$$
Hence for any $n\in\N$ such that $E(\card\,\Gamma_{n-1})\le 3N'$, one
trivially has that $E(\card\, \Gamma_n) \le 3N'$. Hence by (10.8.3)
and induction, (10.8.1) holds for all $n\in\N$. That completes the
proof of Lemma~10.8.
\medskip

{\sc Step 10.9.} Recall the number $\varepsilon \in (0,1)$ fixed in the sentence containing (10.2.1).
Let $M$ be a positive integer sufficiently large
that 
$$6N'/M \le \varepsilon. \leqno{(10.9.1)}$$

By Lemma 10.8, $E(\card(\Gamma_n\cup \Gamma^*_n)) \le 6N'$ for each
positive integer $n$. Hence by Fatou's Lemma,
$$E\left(\liminf_{n\to\infty}\card(\Gamma_n\cup \Gamma^*_n)\right) \le
6N'.$$
Hence by Markov's inequality and (10.9.1),
$$P\left(\liminf_{n\to\infty} \card(\Gamma_n\cup \Gamma^*_n) \ge
M\right) \le 6N'/M\le \varepsilon.  \leqno{(10.9.2)}$$

Next, define the number $\theta\in(0,1)$ by
$$\theta:= (5/8)^{6M}. \leqno{(10.9.3)}$$

Also, for each $n\ge N'$, define the event
$$B_n:= \left\{\card(\Gamma_n\cup \Gamma^*_n)\le M\right\};
\leqno{(10.9.4)}$$
and for each $n\ge N'+1$, define the event
$$D_n := \left\{\omega\in\Omega': \xi^{(n)}_k(\omega) = (0,0,0,0,0,0)
\quad\hbox{for all $k\in \Gamma_{n-1}(\omega) \cup
\Gamma^*_{n-1}(\omega)$}\right\}. \leqno{(10.9.5)}$$

By (10.7.5), (10.9.4), and (10.9.5),
$$\forall\,\,n\ge N',\quad B_n\, \dot\in\, 
\sigma\left(\overline{\xi}^{(n)},\overline{\xi}^{*(n)}\right);
\quad\hbox{and} \leqno{(10.9.6)}$$
$$\forall\,\, n \ge N'+1,\quad
D_n\, \dot\in\, \sigma\left(\overline{\xi}^{(n)},\overline{\xi}^{*(n)}\right).
\leqno{(10.9.7)}$$

{\sc Lemma 10.10.} {\it Refer to (10.9.3), (10.9.4), and (10.9.5).
Suppose $n\ge N'$, $A \in
\sigma(\overline{\xi}^{(n)},\overline{\xi}^{*(n)})$, and 
$P(A\cap B_n)>0$. Then 
$P(D_{n+1}\mid A\cap B_n) \ge \theta$.}
\medskip

{\sc Proof.} By (10.9.4), the event $B_n$ can be partitioned into an
at most countable collection of events of the form $\{\Gamma_n \cap
\Gamma^*_n =S\}$ for sets $S\subset \Z$ such that 
$\card\, S\le M$
(including the empty set). Accordingly, the event $A\cap B_n$ is
partitioned into events $A\cap \{\Gamma_n\cap \Gamma^*_n = S\}$ for
such~$S$.

For any such $S$ such that $P(A\cap \{\Gamma_n\cap \Gamma^*_n =
S\})>0$, one has by (10.7.5) (and the hypothesis of Lemma 10.10
here) that $A \cap \{\Gamma_n\cap \Gamma^*_n =S\} \dot\in
\sigma(\overline{\xi}^n,\overline{\xi}^{*(n)})$ and is therefore
independent of the sequence $\xi^{(n+1)}$
(recall again Construction 6.1(A)(C)), 
and hence by (10.9.5),
Construction 6.1(A) (including (6.1.1) there), and (10.9.3),
$$\eqalign{
&P(D_{n+1}\mid A\cap \{\Gamma_n\cap \Gamma^*_n\} =S) \cr
&= P\left(\xi^{(n+1)}_k = (0,0,0,0,0,0) \,\,\forall\,\,k\in S\mid
A\cap \left\{\Gamma_n\cap \Gamma^*_n\right\} = S\right) \cr
&= P\left(\xi^{(n+1)}_k = (0,0,0,0,0,0)\,\,\forall\,\, k\in S\right)
\cr
&= ((5/8)^6)^{\card\, S} \cr
&\ge (5/8)^{6M} =\theta. \cr}$$
Lemma 10.10 now follows from Remark 2.8.
\medskip

{\sc Lemma 10.11.} $P\left(\bigcup^\infty_{n=N'+1} 
D_n\right) \ge 1-2\varepsilon$.
\medskip

   Here again of course the number $\varepsilon \in (0,1)$
is as fixed in the sentence containing (10.2.1).
One can do better than Lemma 10.11, but that will not be
necessary.
\medskip

{\sc Proof.}  Suppose Lemma 10.11 is false.  That is,
suppose instead that 
$$ s := P\Bigl(\bigcap_{n=N'+1}^{\infty} D_n^c\Bigl)
> 2 \varepsilon.  \leqno (10.11.1) $$
We shall aim for a contradiction.

   Referring to (10.9.3) and (10.11.1), let $\Theta$ be an 
integer such that
$$ \Theta \geq N'+1 \quad {\rm and} \quad
P\Bigl(\bigcap_{n=N'+1}^\Theta D_n^c \Bigl) \leq 
s + \varepsilon \theta /2.
\leqno (10.11.2)
$$
In the arguments below, we shall use more compact
notations for the events in (10.11.1) and (10.11.2):
$$  E := \bigcap_{n=N'+1}^{\infty} D_n^c
\quad {\rm and} \quad   
E' := \bigcap_{n=N'+1}^\Theta D_n^c.  \leqno (10.11.3) $$ 

   Referring to (10.11.2), (10.9.4) and (10.9.5), define the events 
$A_j,\ j \in \{\Theta+1, \Theta+2, \Theta+3, \dots\} \cup \{\infty\}$ as follows:
$$ \leqalignno{
& \indent \indent \indent A_{\Theta+1} := B_{\Theta+1}; & (10.11.4) \cr
& \forall j \in \{\Theta+2, \Theta+3, \Theta+4, \dots\}, \quad
A_j := B_j \bigcap \Bigl(\bigcap_{i=\Theta+1}^{j-1}B_i^c
\Bigl); \quad {\rm and} \cr
&A_{\infty} := \Bigl(\bigcap_{i=\Theta+1}^{\infty}B_i^c
\Bigl). \cr
}$$
These events $A_{\Theta+1}, A_{\Theta+2}, A_{\Theta+3}, \dots$ and 
$A_\infty$ form a partition of $\Omega$.
Also, by (10.11.4) and (10.9.4),  
$ A_\infty\, \dot\subset\,  
\{ \liminf_{n \to \infty} \card(\Gamma_n \cup \Gamma^*_n)
\geq M+1\}$, and hence by (10.9.2), 
$P(A_\infty) \leq \varepsilon$.  Hence
$P(E' \cap A_\infty) \leq \varepsilon$.
Now $P(E') \geq P(E) > 2\varepsilon$ by (10.11.3) and (10.11.1).  Hence
$$ P(E' \cap A_\infty^c) > \varepsilon.  \leqno (10.11.5) $$

   Now suppose 
$j \in \{\Theta+1, \Theta+2, \Theta+3, \dots\}$
is such that $P(E' \cap A_j) > 0$.
(Then $j \geq N' + 2$ by (10.11.2).)\ \ 
By (10.9.6), (10.9.7), (10.11.3), and (10.11.4)
(and Constructions 6.1(A)(C)), 
$ E' \cap A_j \cap B_j = E' \cap A_j\, \dot\in\,  
\sigma (\bar\xi^{(j)}, \bar\xi^{*(j)}) $.
Hence by Lemma 10.10,
$P(D_{j+1}|E' \cap A_j) \geq \theta$.  Of course by
(10.11.3) the events $D_{j+1}$ and $E$ are disjoint.
It follows that $P(E|E' \cap A_j) \leq 1 - \theta$.

   Now by the sentence right after (10.11.4), the events
$E'\cap A_j$, 
$j \in \{\Theta+1, \Theta+2, \Theta+3, \dots\}$
form a partition of the event $E' \cap A_\infty^c$.
Hence by the calculations in the preceding paragraph, 
together with Remark 2.8, one has that
$P(E|E' \cap A_\infty^c) \leq 1 - \theta$.
Since $E \subset E'$ (recall (10.11.3)), this can be 
rewritten as
$ P(E \cap A_\infty^c) \leq  
(1-\theta) P(E' \cap A_\infty^c)$.  Combining that
with (10.11.5), one now has
$$ P(E' \cap A_\infty^c) - P(E \cap A_\infty^c)
\geq \theta \cdot P(E' \cap A_\infty^c) \geq \theta 
\varepsilon. \leqno (10.11.6) $$
Since (again) $E \subset E'$, one also has that
$P(E' \cap A_\infty) - P(E \cap A_\infty) \geq 0$.
Combining that with (10.11.6), one obtains 
$$ P(E') - P(E) \geq \theta \varepsilon. \leqno (10.11.7)$$ 

   However, by (10.11.1). (10.11.2). and (10.11.3),
$ P(E') - P(E) \leq s + \varepsilon \theta / 2 - s
= \varepsilon \theta/2$, which contradicts (10.11.7).
Hence Lemma 10.11 must hold after all.
\medskip

{\sc Lemma 10.12.} {\it For every integer $N\ge N'$ and every
$\omega\in D_{N+1}$, there exists a positive integer $K=K(\omega)$
such that $X^*_k(\omega) = X_k(\omega)$ for all $k\ge K$.}
\medskip

{\it Proof.} Throughout this proof, let $N$ and $\omega$ be fixed
such that
$$N\ge N' \quad\hbox{and}\quad \omega\in D_{N+1}. \leqno{(10.12.1)}$$
Of course $\omega\in \Omega'$ by (10.9.5); and all 
properties of $\omega$ mentioned in the last sentence
of Step 10.4(G) will be tacitly taken for granted
here.

Suppose $k\in \{-N',-N'+1,\dots, L_N(\omega)\}$ is such that
$W^{*(N)}_k(\omega) =1$. Then $k\in \Gamma^*_N(\omega)$ by (10.7.4),
hence $\xi^{(N+1)}_k(\omega) = (0,0,0,0,0,0)$ by (10.12.1) and
(10.9.5), and hence also $\xi^{*(N+1)}_k(\omega) = (0,0,0,0,0,0)$ by (10.12.1) and 
Step 10.4(C). Also (trivially) by (10.7.2), $k = \max
G^{*(N)}((-\infty,k])(\omega)$, and hence (in the convention of Notation 2.1(G))
$$W^{*(N+1)}_k(\omega) = g_{\spaced} 
\left( \xi^{*(N+1)}_j(\omega), \, j\in
G^{*(N)} ((-\infty,k])(\omega)\right) = 0$$
by (10.7.7) and Lemma~4.8 (equation (4.8.1)).

We have shown that (under the hypothesis of Lemma 10.12)
$W^{*(N+1)}_k(\omega) \allowbreak =0$ for every $k\in
\{-N',-N'+1,\dots,L_N(\omega)\}$ such that $W^{*(N)}_k(\omega) =1$.
For all other values of $k\in \{-N',-N'+1,\dots,L_N(\omega)\}$,
$W^{*(N+1)}_k(\omega) =0$ by (the ``asterisk'' analog of) Remark
6.3(B). Hence $W^{*(N+1)}_k(\omega) =0$ for all $k\in\Z$ such that
$-N'\le k\le L_N(\omega)$. By an exactly analogous argument
(using (10.7.3) and (10.7.1) instead of (10.7.4) and (10.7.2)),
$W^{(N+1)}_k(\omega) =0$ for all such $k$. That is,
$$\forall\, k\in \{-N',-N'+1,\dots,L_N(\omega)\}, \quad
W^{*(N+1)}_k(\omega) = W^{(N+1)}_k(\omega) =0. \leqno{(10.12.2)}$$

Next, suppose $u\in\N$. Then $W^{*(N)}_{H(N,u)(\omega)} (\omega) =
W^{(N)}_{H(N,u)(\omega)}(\omega) =1$ by Lemma 10.6 (equation
(10.6.1)), (10.3.1), and (10.3.2); and hence by (10.7.7),
$$W^{*(N+1)}_{H(N,u)(\omega)}(\omega) =
g_{\spaced}\left(\xi^{*(N+1)}_j (\omega), \, j\in
G^{*(N)}((-\infty,H(N,u)(\omega)])(\omega)\right).
\leqno{(10.12.3)}$$ Again, for every $k\in
G^{*(N)}((-\infty,H(N,u)(\omega)])(\omega)$ such that $-N'\le k\le
L_N(\omega)$, one has that $\xi^{*(N+1)}_k (\omega) =
(0,0,0,0,0,0)$, by (10.7.2) and the arguments in the
paragraph after that of (10.12.1). Referring to (10.3.1)--(10.3.2)
and applying Lemma 4.8 (equation (4.8.2)), one can delete those particular elements
$k$ from the set $G^{*(N)}((-\infty,H(N,u)(\omega)])(\omega)$
without changing the right hand side of (10.12.3). That is, one
obtains (see (10.7.2) and (10.3.1)--(10.3.2) again)
$$\leqalignno{
&&(10.12.4)\cr
 W^{*(N+1)}_{H(N,u)}(\omega) &=
g_{\spaced}\Bigl(\xi^{*(N+1)}_j(\omega), \, j\in G^{*(N)}
((-\infty,-N'-1])(\omega)  \cr 
&\qquad\qquad \qquad \quad\ \  \cup
G^{*(N)}([L_N(\omega) +1,H(N,u)(\omega)])(\omega)\Bigr)  \cr
&=g_{\spaced}\Bigl(\xi^{*(N+1)}_j(\omega),\, j\in
G^{*(N)}((-\infty,-N'-1])(\omega) \cr &\qquad\qquad \quad\cup
\{H(N,1)(\omega),H(N,2)(\omega),\dots,H(N,u)(\omega)\}\Bigr)  \cr}$$
Now by (10.12.1) and Step 10.4(C), $\xi^{*(N+1)}_j(\omega) =
\xi^{(N+1)}_j(\omega)$ for every $j\in\Z$. Also,
$$G^{*(N)}((-\infty,-N'-1])(\omega) =
G^{(N)}((-\infty,-N'-1])(\omega)$$ by (10.7.1), (10.7.2), and
(10.5.2). Hence by (10.12.4), 
$$\leqalignno{
&&(10.12.5)\cr
 W^{*(N+1)}_{H(N,u)(\omega)}(\omega) = \,
&g_{\spaced}\Bigl(\xi^{(N+1)}_j(\omega),\, j\in
G^{(N)}((-\infty,-N'-1])(\omega)  \cr & \indent \cup
\{H(N,1)(\omega),H(N,2)(\omega),\dots,H(N,u)(\omega)\}\Bigr). \cr}$$
By an exactly analogous argument, one obtains
$$W^{(N+1)}_{H(N,u)(\omega)}(\omega) = [\hbox{RHS of (10.12.5)}],$$
and hence by (10.12.5) itself,
$$W^{*(N+1)}_{H(N,u)(\omega)} (\omega)  
= W^{(N+1)}_{H(N,u)(\omega)}(\omega).
\leqno{(10.12.6)}$$

Equation (10.12.6) was obtained for arbitrary $u\in\N$. For all other
$k\ge L_N(\omega)+1$ besides
$\{H(N,1)(\omega),H(N,2)(\omega),\dots,\,\}$, one has that
$W^{*(N)}_k(\omega) = W^{(N)}_k(\omega) \not= 1$
by Lemma 10.6 (equation (10.6.1)) and (10.3.2), and hence
$W^{*(N+1)}_k(\omega) = W^{(N+1)}_k(\omega) =0$ by Remark 6.3(B) and
its ``asterisk'' counterpart. Combining that with (10.12.6), one now
has that
$$\forall\,k\ge L_N(\omega) +1,\quad W^{*(N+1)}_k (\omega) =
W^{(N+1)}_k(\omega). \leqno{(10.12.7)}$$

Now by (10.12.2), (10.12.7), and (10.5.2), one now has that
$$\forall\, k\in {\Z},\quad W^{*(N+1)}_k(\omega) = W^{(N+1)}_k(\omega).
\leqno{(10.12.8)}$$

Now recall from Step 10.4(C) that $\xi^{*(n)}_k(\omega) =
\xi^{(n)}_k(\omega)$ for all $n\ge N'+1$ and all $k\in\Z$. Recall (10.12.1) again.  Starting
with (10.12.8) and using (6.2.6) and (6.2.7) (and their ``asterisk''
counterparts) and induction on $n$, one has that
$$\forall\, n\ge N+1,\,\,\forall\, k\in{\Z},\quad W^{*(n)}_k(\omega) =
W^{(n)}_k(\omega); \,\, \hbox{and} \leqno{(10.12.9)}$$
$$\, \forall\, n\ge N+1, \, \forall\, k\in{\Z},\, \forall\, j\ge 0,
\,\, \Psi^*(n,k,j)(\omega) = \Psi(n,k,j)(\omega).
\leqno{(10.12.10)}$$

Now recall (10.3.9). Our candidate for the integer $K=K(\omega)$ in
Lemma 10.12 will be $H(N,1)(\omega)$.

By (10.3.9) and Step 10.4(B), one has that
$$\forall\, k\ge H(N,1)(\omega), \quad X^{*(0)}_k(\omega) =
X^{(0)}_k(\omega). \leqno{(10.12.11)}$$

Now keep in mind that $H(N,1)(\omega) \ge H(n,1)(\omega)$ for all
$n\in\{1,\dots,N\}$ by (10.3.9). By (6.3.5) and Lemma 10.6 (equation
(10.6.4), applied twice), one has that
$$\leqalignno{\forall\, n\in\{1,\dots,N\}, \qquad \forall\,
k&\ge H(N,1)(\omega), &(10.12.12) \cr
k-\Psi^*(n,k,0)(\omega) = k-\Psi(n,k,0)(\omega) &\ge k-\Psi(N,k,0)(\omega)\cr
& \ge H(N,1)(\omega). \cr}$$
By Lemma 10.6 (equation (10.6.2)),
$$\forall\, n\in \{1,\dots,N\},\,\, \forall\, k\ge H(N,1)(\omega),
\quad \delta^{*(n)}_k(\omega) = \delta^{(n)}_k(\omega).
\leqno{(10.12.13)}$$
Also, by (7.1.2) and its ``asterisk'' counterpart, equation
(10.12.12) for $n=N$, and equation (10.12.8),
$$\forall\, k\ge H(N,1)(\omega),\quad \delta^{*(N+1)}_k (\omega) =
\delta^{(N+1)}_k(\omega). \leqno{(10.12.14)}$$
Next, recall from Step 10.4(C) (and (7.3.1) and its
``asterisk'' counterpart) that $\zeta^{*(n)}_k(\omega) =
\zeta^{(n)}_k(\omega)$ for all $n\in\N$ and all $k\ge N'+1$.  Since $H(N,1)(\omega) \ge N'+1$ by
(10.3.9), one has by (10.12.12) that
$$\leqalignno{\qquad\qquad\qquad\forall\, n\in \{1,\dots,N\}, \,\,\, &\quad \forall\,
k\ge H(N,1)(\omega), &(10.12.15) \cr
\zeta^{*(n)}_{k-\Psi^*(n,k,0)(\omega)}(\omega) &=
\zeta^{*(n)}_{k-\Psi(n,k,0)(\omega)}(\omega) =
\zeta^{(n)}_{k-\Psi(n,k,0)(\omega)}(\omega). \cr}$$

Next by Step 10.4(C) and (10.12.1), one has that for every $n\ge N+1$
(in fact every $n\ge N'+1)$ and every $k\in\Z$,
$\zeta^{*(n)}_k(\omega) = \zeta^{(n)}_k(\omega)$. Hence by
(10.12.10), for every $n\ge N+1$ and every $k\in\Z$,
$\zeta^{*(n)}_{k-\Psi^*(n,k,0)(\omega)}(\omega) =
\zeta^{(n)}_{k-\Psi(n,k,0)(\omega)}(\omega)$. Also, by (10.12.9),
(10.12.10), and (7.1.2) and its ``asterisk'' counterpart, for every
$n\ge N+2$ and every $k\in\Z$, $\delta^{*(n)}_k(\omega) =
\delta^{(n)}_k(\omega)$. Combining these facts with (10.12.13),
(10.12.14),  and (10.12.15), one now has that
$$\forall\, n\in{\N},\,\,\forall\, k\ge H(N,1)(\omega),\quad
\delta_k^{*(n)}(\omega) = \delta^{(n)}_k(\omega); \quad \hbox{and}
\leqno{(10.12.16)}$$
$$\forall\, n\in{\N},\, \forall\, k\ge H(N,1)(\omega), \quad
\zeta^{*(n)}_{k-\Psi^*(n,k,0)(\omega)}(\omega) =
\zeta^{(n)}_{k-\Psi(n,k,0)(\omega)}(\omega). \leqno{(10.12.17)}$$

By (7.1.7) and (7.1.4) and their ``asterisk'' counterparts, together
with (10.12.16), one obtains
$$\forall\,n\in {\N},\, \forall\, k \ge H(N,1)(\omega), \quad
J^*(n,k)(\omega) = J(n,k)(\omega); \quad \hbox{and}
\leqno{(10.12.18)}$$
$$\forall\, k\ge H(N,1)(\omega), \quad
N^*_k(\omega) = N_k(\omega). \leqno{(10.12.19)}$$

Now by (7.5.1) (and its ``asterisk'' counterpart), 
(10.12.11), (10.12.17) (recall (7.3.1) and its ``asterisk''
counterpart),  
(10.12.18), and (10.12.19), 
one has that $X^*_k(\omega) = X_k(\omega)$
for all $k\ge H(N,1)(\omega)$. Thus Lemma 10.12 holds with $K=
K(\omega) = H(N,1)(\omega)$.
\medskip

{\sc Step 10.13.} By Lemma 10.12 and Lemma 10.11,
$$P\left(\bigcup_{J\in\N} \left\{X^*_k = X_k\,\,\forall\, k\ge
J\right\}\right) \ge P\left(\bigcup^\infty_{n=N'+1} D_n\right) \ge
1-2\varepsilon. \leqno{(10.13.1)}$$
Also (trivially) for each $J\in\N$,
$$\bigl\{X^*_k =X_k \,\,\forall\, k\ge J\bigr\}\subset \bigl\{X^*_k =
X_k\,\, \forall\, k\ge J+1\bigr\}.$$
Hence by (10.13.1), $\lim_{J\to\infty} P(X^*_k =X_k$ $\forall\, k\ge J)
\ge 1-2\varepsilon$. Accordingly, let $J'$ be a positive integer such
that
$$J'\ge N'+1\quad \hbox{and} \leqno{(10.13.2)}$$
$$P\left(X^*_k =X_k\,\,\forall\, k\ge J'\right) \ge 1-3\varepsilon.
\leqno{(10.13.3)}$$

Now recall (10.1.1). One has that $F'\in \sigma(X_k$, $|k| \ge J')$.
Thus (see e.g.\ [1, Theorem 20.1(i), trivially extended],
and recall property (A) in Theorem 1.1) 
there exists a Borel set $B'\subset \{-1,1\}^{\N}$ such that
$$F' = \left\{(X_{J'},X_{-J'}, X_{J'+1},X_{-J'-1}, X_{J'+2},
X_{-J'-2},\dots\,)\in B'\right\} \leqno{(10.13.4)}$$
Define the (``asterisk counterpart'') event
$$F''' := \left\{(X^*_{J'},X^*_{-J'}, X^*_{J'+1},X^*_{-J'-1}, X^*_{J'+2},
X^*_{-J'-2},\dots\,)\in B'\right\} \leqno{(10.13.5)}$$

By (10.13.2), (10.13.3), and (10.5.3),
$$P\left(\hbox{$X^*_k =X_k$ for all}\,\,
k\in{\Z} \,\,\hbox{such that}\,\, |k|\ge J'\right) \ge
1-3\varepsilon.$$ Hence by (10.13.4), (10.13.5), and a simple
argument,
$$P(F'\, \triangle \, F''') \le 3\varepsilon \leqno{(10.13.6)}$$
(where $\triangle$ denotes symmetric difference). Hence by a simple
standard argument,
$$P(F''') \le P(F') +3\varepsilon. \leqno{(10.13.7)}$$

Refer to Construction 7.4(B).  
By Lemma 7.6 (and Remark 2.6(B)), 
one has that $\sigma(X) \dot\subset \sigma(\eta)$, the
$\sigma$-field generated by the array (10.4.1). The ``asterisk''
counterpart is $\sigma(X^*) \dot\subset \sigma(\eta^*)$, the
$\sigma$-field generated by the array (10.4.2). Hence by (10.4.4) and (10.13.5),
$$\hbox{the events $F''$ and $F'''$ are independent.}
\leqno{(10.13.8)}$$

Now by (10.2.3), (10.13.6), and a standard elementary argument (see
e.g.\ [5, Vol.\ 1, Appendix, Section A053(V)]),
$$\eqalign{
P((F'' \cap F''')\, \triangle\,  F') &= P((F'' \cap F''')\, \triangle \, (F' \cap F'))
\cr
&\le P(F''\,  \triangle\, F') +P(F'''\, \triangle \, F') \le 4\varepsilon. \cr}$$
Hence by a standard simple argument, $P(F'' \cap F''') \ge P(F')
-4\varepsilon$. Hence by (10.13.8), (10.2.4), and (10.13.7),
$$\eqalign{
P(F') -4\varepsilon &\le P(F''\cap F''') = P(F'')\cdot P(F''') \cr
&\le [P(F') +\varepsilon]\cdot [P(F') +3\varepsilon]. \cr}$$
However, that contradicts (10.2.1). Hence (10.1.2) must be false, and \hfil\break
$\cT_{\double}(X)$ is trivial after all. That completes the proof of
Property (D) in Theorem 1.1, and of Theorem 1.1 itself.

\bigskip
\bigskip
\bigskip
{\bf Acknowledgements.} The author thanks Jon Aaronson and Benjamin
Weiss for their encouragement and helpful comments and for pointing
out reference [16]; and the author also thanks Ron Peled for point
out reference [20].

\vfill\eject
 \centerline{\bf References}
\smallskip

\refs{[1]} P.\  Billingsley, {\it Probability and Measure}, 3rd ed.\ 
(Wiley, New York, 1995).

\refs{[2]} R.C. Bradley, A bilaterally deterministic $\rho$--mixing
stationary random sequence, {\it Trans.\ Amer.\ Math.\ Soc.} 294 (1986), 55-66.

\refs{[3]} R.C.\  Bradley, A stationary, pairwise independent,
absolutely regular sequence for which the central limit theorem
fails, {\it Probab.\ Th.\ Rel.\ Fields} 81 (1989), 1-10.

\refs{[4]} R.C.\ Bradley, On a stationary, triple-wise independent,
absolutely regular counterexample to the central limit theorem, {\it
Rocky Mountain J.\ Math.} 37 (2007), 25-44.

\refs{[5]} R.C.\ Bradley, {\it Introduction to Strong Mixing Conditions,}
Volumes 1, 2, and 3 (Kendrick Press, Heber City, Utah, 2007).

\refs{[6]} R.C.\ Bradley and A.R.\ Pruss, A strictly
stationary, $N$-tuplewise independent counterexample to the
central limit theorem,
{\it Stochastic Process.\ Appl.} 119 (2009), 3300-3318.

\refs{[7]} P.J.\ Brockwell and R.A.\ Davis, {\it Time Series: Theory
and Methods,} 2nd ed.\ (Springer, New York, 1991).

\refs{[8]} R.M.\ Burton, M.\ Denker, and M.\ Smorodinsky, Finite
state bilaterally deterministic strongly mixing processes, {\it Israel J.\ Math.} 95 (1996), 115-133.

\refs{[9]} R.\ Cogburn, Asymptotic properties of stationary
sequences, {\it Univ.\ Calif.\ Publ.\ Statist.} 3 (1960), 99-146.

\refs{[10]} J.A.\ Cuesta and C.\ Matran, On the asymptotic
behavior of sums of pairwise independent random variables,
{\it Statist.\ Probab.\ Letters\/} 11 (1991) 201-210.
(Correction Ibid.\ 12 (1991) 183.)

\refs{[11]} Yu.A.\ Davydov, Mixing conditions for Markov chains, {\it
Theory Probab.\ Appl.} 18 (1973), 312-328.

\refs{[12]} H.\ Dehling, M.\ Denker, and W.\ Philipp, Central limit
theorems for mixing sequences of random variables under minimal
conditions, {\it Ann.\ Probab.} 14 (1986), 1359-1370.

\refs{[13]} M.\ Denker, Uniform integrability and the central limit
theorem for strongly mixing processes, In: {\it Dependence in
Probability and Statistics} (E.\ Eberlein and M.S. Taqqu, eds.), {\it
Progress in Probability and Statistics,} Volume 11, pp.\ 269-274
(Birkh\"auser, Boston, 1986).

\refs{[14]} N.\ Etemadi, An elementary proof of the strong law of
large numbers, {\it Z.\ Wahrsch.\ verw.\ Gebiete} 55
(1981), 119-122.

\refs{[15]} W.\ Feller, {\it An Introduction to Probability Theory
and its Applications}, Volume~1, 3rd\ ed. (Wiley, New York, 1968).

\refs{[16]} L.\ Flaminio, Mixing $k$--fold independent processes of
zero entropy, {\it Proc. Amer.\ Math.\ Soc.} 118 (1993), 1263-1269.

\refs{[17]} N.\ Friedman and D.S.\ Ornstein, On isomorphism of weak
Bernoulli transformations, {\it Advances in Math.} 5 (1970),
365-394.

\refs{[18]} N.\ Herrndorf, Stationary strongly mixing sequences not
satisfying the central limit theorem, {\it Ann. Probab.} 11 (1983),
809-813.

\refs{[19]} S. Janson, Some pairwise independent sequences for which
the central limit theorem fails, {\it Stochastics} 23 (1988),
439-448.

\refs{[20]} F.S.\ MacWilliams and N.J.A.\ Sloane, {\it The Theory of
Error--Correcting Codes} (North Holland, Amsterdam, 1977).

\refs{[21]} F.\ Merlev\`ede and M.\ Peligrad, The functional central
limit theorem under the strong mixing condition, {\it Ann.\ Probab.}
28 (2000), 1336-1352.

\refs{[22]} T.\ Mori and K.\ Yoshihara, A note on the central limit
theorem for stationary strong-mixing sequences, Yokohama Math.\ J.\
34 (1986), 143-146.

\refs{[23]} D.S.\ Ornstein, Two Bernoulli shifts with infinite
entropy are isomorphic, {\it Advances in Math.} 5 (1970), 339-348.

\refs{[24]} D.S.\ Ornstein, Factors of Bernoulli shifts are
Bernoulli shifts, {\it Advances in Math.} 5 (1970), 349-364.

\refs{[25]} D.S.\ Ornstein and B.\ Weiss, Every transformation is
bilaterally deterministic, {\it Israel J.\ Math.} 21 (1975),
154-158.

\refs{[26]} K.\ Petersen, {\it Ergodic Theory,} first paperback
edition (Cambridge University Press, Cambridge, 1989).

\refs{[27]} J.\ Pitman, {\it Probability} (Springer, New York, 1993).

\refs{[28]} A.R.\ Pruss, A bounded $N$-tuplewise independent and
identically distributed counterexample to the CLT, {\it Probab.\
Th.\ Rel.\ Fields} 111 (1998), 323-332.

\refs{[29]} M.\ Rosenblatt, A central limit theorem and a strong
mixing condition, {\it Proc.\ Nat.\ Acad.\ Sci.\ USA} 42 (1956), 43-47.

\refs{[30]} M.\ Rosenblatt, Stationary Markov chains and
independent random variables,  
{\it J.\ Math.\ Mech.} 9 (1960), 945-949.

\refs{[31]} P.C.\ Shields, {\it The Ergodic Theory of Discrete Sample
Paths} (American Mathematical Society, Providence, R.I., 1996).

\refs{[32]} M.\ Smorodinsky, A partition on a Bernoulli shift which
is not weakly Bernoulli, {\it Math.\ Systems Th.} 5 (1971), 201-203.

\refs{[33]} V.A.\ Volkonskii and Yu.A.\ Rozanov, Some limit theorems
for random functions~I, {\it Theory Probab.\ Appl.} 4 (1959),
178-197.

\bye